\documentclass{article}

 \usepackage{amssymb}					
 \usepackage{amsmath}					
 \usepackage{amsthm}					
 \usepackage{latexsym}					
 \usepackage{subfigure}
 \usepackage{epsfig}					
 \usepackage{pstricks,pst-node,pst-text,pst-tree}
 \usepackage{amscd}

 \newtheorem{theorem}{Theorem}[section]
 
 \newtheorem{corollary}[theorem]{Corollary}
 \newtheorem{lemma}[theorem]{Lemma}
 \newtheorem{proposition}[theorem]{Proposition}

 \theoremstyle{definition}\newtheorem{definition}[theorem]{Definition}

 \theoremstyle{remark}\newtheorem{remark}[theorem]{Remark}
 
 \newtheorem{example}[theorem]{Example}
 \newtheorem{warning}[theorem]{Warning}

 \numberwithin{equation}{section}

 \newcommand{\para}{\medskip \indent}
 \newcommand{\one}{\ensuremath{(\mathrm{i})}}
 \newcommand{\two}{\ensuremath{(\mathrm{ii})}}
 \newcommand{\three}{\ensuremath{(\mathrm{iii})}}
 
 \newcommand{\plusone}{\ensuremath{(+1)}}
 \newcommand{\minusone}{\ensuremath{(-1)}}
 \newcommand{\plusminusone}{\ensuremath{(\pm 1)}}

 \newcommand{\Amp}{\ensuremath{\operatorname{Amp}}}
 \newcommand{\Aut}{\ensuremath{\operatorname{Aut}}}
 \newcommand{\Band}{\operatorname{Band}}
 \newcommand{\C}{\ensuremath{\mathbb{C}}}
 \newcommand{\Cycle}{\operatorname{Cycle}}
 
 \newcommand{\GL}{\ensuremath{\operatorname{GL}}}
 \newcommand{\GExt}{\ensuremath{\operatorname{\hbox{$G$}-Ext}}}
 \newcommand{\GHom}{\ensuremath{\operatorname{\hbox{$G$}-Hom}}}
 \newcommand{\Ext}{\ensuremath{\operatorname{Ext}}}
 
 \newcommand{\Hom}{\ensuremath{\mbox{Hom}}}
 \newcommand{\Irr}{\ensuremath{\operatorname{Irr}}}
 \newcommand{\0}{\ensuremath{\operatorname{0}}} 
 \newcommand{\I}{\ensuremath{\operatorname{I}}}
 \newcommand{\II}{\ensuremath{\operatorname{I\!I}}}
 \newcommand{\III}{\ensuremath{\operatorname{I\!I\!I}}} 
 \newcommand{\LL}{\mathbf{L}}
 
 \newcommand{\Paut}{\ensuremath{\operatorname{\mathbb{P}Aut}}}
 \newcommand{\Pic}{\ensuremath{\operatorname{Pic}}}
 \newcommand{\Q}{\ensuremath{\mathbb{Q}}}
 \newcommand{\Quot}{\operatorname{Quot}}
 \newcommand{\R}{\ensuremath{\mathbb{R}}} 
 \newcommand{\Ro}{\ensuremath{R_{1}}}
 \newcommand{\Rt}{\ensuremath{R_{2}}}
 \newcommand{\RR}{\mathbf{R}}
 
 \newcommand{\SL}{\ensuremath{\operatorname{SL}}}
 
 \newcommand{\Spec}{\ensuremath{\operatorname{Spec}}}
 \newcommand{\Supp}{\ensuremath{\operatorname{Supp}}}

 \newcommand{\Z}{\ensuremath{\mathbb{Z}}}

 \newcommand{\ch}{\ensuremath{\operatorname{ch}}}
 \newcommand{\codim}{\operatorname{codim}} 
 \newcommand{\degree}{\ensuremath{\operatorname{deg}}}
 \newcommand{\diag}{\ensuremath{\operatorname{diag}}}

 \newcommand{\ellzero}{\ensuremath{\operatorname{\ell_{0}}}}
 \newcommand{\ev}{\ensuremath{\operatorname{ev}}}
 \newcommand{\gen}{\text{gen}}
 \newcommand{\git}{\ensuremath{\operatorname{/\!\!/}}}
 \newcommand{\ghilb}{\ensuremath{G}\operatorname{-Hilb}}
 \newcommand{\ind}{\operatorname{ind}}
 \newcommand{\ltensor}{\overset{\LL}{\otimes}}
 \newcommand{\mc}{\mathcal}
 \newcommand{\owe}{\ensuremath{\mathcal{O}}}

 \newcommand{\res}{\ensuremath{\operatorname{p}}}
 \newcommand{\socle}{\ensuremath{\operatorname{socle}}}
 \newcommand{\st}{\ensuremath{\operatorname{\bigm{|}}}}

 \newcommand{\Ghilb}[1]{\ensuremath{G\operatorname{-Hilb}(\mathbb{C}^{#1})}}
 
 \newcommand{\taut}[1]{\ensuremath{\mathcal{R}_{#1}}}
 \newcommand{\uni}[1]{\ensuremath{\mathcal{U}_{#1}}}
 \newcommand{\FM}[1]{\ensuremath{\Phi_{#1}}}
 \newcommand{\fm}[1]{\ensuremath{\varphi_{#1}}}
 \newcommand{\lb}[1]{\ensuremath{L_{#1}}}
 \newcommand{\fmstar}[1]{\ensuremath{\varphi_{#1}^{*}}}
 \newcommand{\cont}[1]{\ensuremath{\operatorname{cont_{#1}}}}
 \newcommand{\Dom}[1]{\text{Fill}({#1})}

\makeatletter
 \newlength{\typesize}
\makeatother

\newlength{\vvoff}
\newlength{\hhoff}

\newcommand{\locateoffcenter}[1]{%
\addtolength{\vvoff}{-0.25\typesize}%
\raisebox{\vvoff}{\hspace{\hhoff}\makebox(0,0){\smash{#1}}}
}
\newcommand{\object}[1]{%
\setlength{\vvoff}{0pt}%
\setlength{\hhoff}{0pt}%
\locateoffcenter{#1}
}

\newcommand{\swlabel}[1]{%
\setlength{\vvoff}{-0.5\typesize}%
\setlength{\hhoff}{0.75\typesize}%
\locateoffcenter{#1}
}
\newcommand{\nwlabel}[1]{%
\setlength{\vvoff}{-0.5\typesize}%
\setlength{\hhoff}{-0.75\typesize}%
\locateoffcenter{#1}
}
\newcommand{\selabel}[1]{%
\setlength{\vvoff}{0.5\typesize}%
\setlength{\hhoff}{0.75\typesize}%
\locateoffcenter{#1}
}
\newcommand{\nelabel}[1]{%
\setlength{\vvoff}{0.5\typesize}%
\setlength{\hhoff}{-0.75\typesize}%
\locateoffcenter{#1}
}

 \title{Flops of \protect\(\ghilb\protect\) and equivalences of derived categories by variation of GIT quotient}
 \author{Alastair Craw and Akira Ishii}
 \date{}

\begin{document}
 
 \maketitle

 \begin{abstract}
 For a finite subgroup \(G \subset \SL(3,\C)\),  Bridgeland, King and Reid proved that the moduli space of \(G\)-clusters is a crepant resolution of the quotient \(\C^{3}/G\).  This paper considers the moduli spaces \(\mc{M}_{\theta}\), introduced by Kronheimer and further studied by Sardo Infirri,  which coincide with \(\ghilb\) for a particular choice of GIT parameter \(\theta\). For \(G\) Abelian,  we prove that every projective crepant resolution of \(\C^{3}/G\) is isomorphic to \(\mc{M}_{\theta}\) for some parameter \(\theta\).  The key step is the description of GIT chambers in terms of the \(K\)-theory of the moduli space via the appropriate Fourier-Mukai transform. We also uncover explicit equivalences between the derived categories of moduli \(\mc{M}_{\theta}\) for parameters lying in adjacent GIT chambers. 
 
 
 \end{abstract}

 \section{Introduction}
 \label{sec:intro}
 For a finite subgroup \(G\subset \SL(3,\C)\),  write \(\ghilb\) for the Hilbert scheme parametrising \emph{\(G\)-clusters} in \(\C^{3}\), i.e.,  the scheme parametrising \(G\)-invariant subschemes \(Z\subset\C^{3}\) of dimension zero with global sections \(H^{0}(\mathcal{O}_{Z})\) isomorphic as a \(\C[G]\)-module to the regular representation of \(G\).   Nakamura~\cite{Nakamura:ago} proved that \(\ghilb\) is a crepant resolution of \(\C^{3}/G\) when \(G\) is Abelian and conjectured that the same holds for all finite subgroups \(G\subset \SL(3,\C)\).   Bridgeland,  King and Reid~\cite{Bridgeland:mim} subsequently proved the conjecture by establishing an equivalence of derived categories \(\Phi\colon D(Y)\to D^{G}(\C^{3})\),  where \(D(Y)\) and \(D^{G}(\C^{3})\) are the bounded derived categories of coherent sheaves on \(Y = \ghilb\) and \(G\)-equivariant coherent sheaves on \(\C^{3}\) respectively.

 This article generalises the notion of \(G\)-cluster.  For a finite subgroup \(G\subset \GL(n,\C)\),  a \emph{\(G\)-constellation} is a \(G\)-equivariant coherent sheaf \(F\) on \(\C^{n}\) with global sections \(H^{0}(F)\) isomorphic as a \(\C[G]\)-module to the regular representation \(R\) of \(G\).   Set
 \[
 \Theta \! := \{\, \theta\in \Hom_{\Z}(R(G),\Q) \st \theta(R) = 0\,\},
 \]
 where \(R(G)\) is the representation ring of \(G\).  For \(\theta\in \Theta\),  a \(G\)-constellation \(F\) is said to be \emph{\(\theta\)-stable} (or \emph{\(\theta\)-semistable}) if every proper \(G\)-equivariant coherent subsheaf \(0 \subset E \subset F\) satisfies \(\theta(E) > 0\) (or \(\theta(E) \geq 0\)).  We construct moduli spaces \(\mc{M}_{\theta}\) of \(\theta\)-stable \(G\)-constellations and \(\overline{\mc{M}_{\theta}}\) of \(\theta\)-semistable \(G\)-constellations using geometric invariant theory (GIT).  For generic \(\theta\), the moduli \(\mc{M}_{\theta}\) depends only upon the (open) GIT chamber \(C\subset \Theta\) containing \(\theta\in \Theta\),  so we write \(\mc{M}_{C}\) in place of \(\mc{M}_{\theta}\) for any \(\theta\subset C\).  

 These moduli spaces have appeared before under the guise of moduli of representations of the McKay quiver.  The hyperkahler quotient construction by Kronheimer~\cite{Kronheimer:hkq} for \(G\subset \SL(2,\C)\) was generalised using geometric invariant theory (GIT) by Sardo Infirri~\cite{Sardo-Infirri:ros}.  
The notion of stability for quiver representations was introduced by King~\cite{King:mor},  where it was shown that quiver stability coincides with the notion of stability occuring in the GIT construction of \cite{Sardo-Infirri:ros}.  Thus, \(\mc{M}_{C}\) and \(\overline{\mc{M}_{C}}\) are examples of the type introduced by King.  More recently,  these moduli have appeared in the physics literature as moduli of \(D\)0-branes on the orbifold \(\C^{n}/G\),  see Douglas et.~al.~\cite{Douglas:ale,Douglas:ord}. 


 Our interest in these moduli spaces begins with the observation by Ito and Nakajima~\cite{Ito:mc3} that \(\ghilb = \mc{M}_{C_{0}}\) for some chamber \(C_{0}\subset \Theta\),   i.e.,  there is a crepant resolution \(\mc{M}_{C_{0}}\to\C^{3}/G\) for \(G\subset \SL(3,\C)\).  More generally,  the method of \cite{Bridgeland:mim} shows that for every chamber \(C\subset \Theta\) there is an equivalence of derived categories \(\Phi_{C}\colon D(\mc{M}_{C})\to D^{G}(\C^{3})\) and a crepant resolution \(\tau\colon \mc{M}_{C}\to\C^{3}/G\) (see \S\ref{sec:bkr}). It is therefore natural to ask whether every (projective) crepant resolution may be realised as a moduli space \(\mc{M}_{C}\) for some chamber \(C\).  The main result of this paper answers this question affirmatively in the Abelian case:

 \begin{theorem}
 \label{thm:main}
 For a finite Abelian subgroup \(G\subset \SL(3,\C)\),  suppose \(Y\to\C^{3}/G\) is a projective crepant resolution.  Then \(Y \cong \mc{M}_C\) for some chamber \(C\subset \Theta\).
 \end{theorem}
 
 Since every projective crepant resolution is obtained by a finite sequence of flops from \(\mc{M}_{C_{0}} = \ghilb\),  it is enough to show that if \(Y\cong \mathcal{M}_C\) for some chamber \(C\) then for any flop \(Y'\) of \(Y\) there is a chamber \(C'\) (not necessarily adjacent to \(C\)) such that \(\mathcal{M}_{C'} \cong Y'\).

 We now review the key elements in more detail.  The first step is to understand the walls of chambers in \(\Theta\) in terms of Fourier--Mukai transforms.  The derived equivalence \(\FM{C}\) induces a \(\Z\)-linear isomorphism \(\varphi_{C}\colon K_{0}(\mc{M}_{C})\to R(G)\) between the Grothendieck group of coherent sheaves supported on \(\tau^{-1}(0)\) and the representation ring.  The pairing between \(K(\mc{M}_{C})\) and \(K_{0}(\mc{M}_{C})\) is perfect so \(\varphi_{C}^{*}\colon \Hom_{\Z}(R(G),\Q) \to K(\mc{M}_{C})_{\Q}\) is an isomorphism.  The restriction of \(\varphi_{C}^{*}\) to the subspace \(\Theta\subset \Hom_{\Z}(R(G),\Q)\) has image equal to \(F^{1}\subset K(\mc{M}_{C})_{\Q}\),  where \(F^{i}\subset K(\mc{M}_{C})_{\Q}\) denotes the subspace spanned by sheaves with support of codimension at least \(i\).  The quotient \(F^{1}/F^{2}\) is isomorphic to \(\Pic(\mc{M}_{C})_{\Q}\),  so we have three sides of the following commutative diagram:
 \[
 \begin{CD}
 \Theta @>{\fmstar{C}}>> F^1 \\
 @V{L_{C}}VV              @VV{\res}V  \\
 \Pic(Y)_{\Q}  @<{\sim}<< F^1/F^2
 \end{CD}
 \]
 The map \(L_{C}\) that completes the diagram (see \S\ref{sec:classificationofwalls} for the construction) sends a parameter \(\theta\in C\) to the ample fractional line bundle \(\owe_{\mc{M}_{\theta}}(1)\) on \(\mc{M}_{\theta} = Y\) canonically derived in the GIT construction of \(\mc{M}_{\theta}\).  In particular,  \(L_{C}(C)\) lies inside the ample cone \(\Amp(Y)\subset \Pic(Y)_{\Q}\).  The inequalities defining the ample cone lift to give inequalities of the form
 \[
 \theta(\fm{C}(\owe_{\ell})) > 0 \quad \text{for exceptional curves }\ell\subset Y
 \]
 that are satisfied by each \(\theta\in C\).  Some of these inequalities may be redundant,  otherwise they determine walls of `type \I' or `type \III' (this means that the wall determines a birational contraction of type \I\ or \III\ from \(\mc{M}_{C}\),  see Definition~\ref{defn:type}).  In contrast to the two-dimensional case,  there are additional walls of \(C\) that do not arise as lifts of walls from the ample cone.  These walls can also be described in terms of the Fourier-Mukai transform using the \emph{tautological bundles} \(\{\taut{\rho}\}\) (see \S2.1 for the definition).  Theorem~\ref{thm:chamberstrong} asserts that \(C\) is determined by the inequalities coming from the ample cone together with inequalities of the form
 \[
 \theta(\fm{C}(\taut{\rho}^{-1}\otimes \omega_{D})) < 0\quad\text{and}\quad\theta(\fm{C}(\taut{\rho}^{-1}\vert_D)) > 0
 \]
 for every compact divisor \(D\) and irreducible representation \(\rho\).  Many of these inequalities are redundant,  but those that are not determine walls of `type \0'. 

 The second main step is to understand how the moduli \(\mc{M}_{C}\),  their tautological bundles \(\taut{C}\) and the Fourier--Mukai transforms \(\FM{C}\) change as \(\theta\) passes through a wall from \(C\) to an adjacent chamber \(C'\).  The description of chambers in terms of Fourier--Mukai transforms leads to the following result (see \S\ref{sec:twists}): 
 \begin{theorem}
 \label{thm:walls}
 Let \(C\) and \(C'\) be adjacent chambers in \(\Theta\) with \(W = \overline{C}\cap\overline{C'}\).  Then the relation between the moduli spaces \(\mc{M}_{C},\mc{M}_{C'}\) and the explicit form of the derived equivalence 
 \[
 \FM{C'}^{-1}\circ \FM{C}\colon D(\mc{M}_{C}) \overset{\sim}{\longrightarrow} D(\mc{M}_{C'})
 \]
 is determined by the type of the wall \(W\):
 \begin{description}
 \item[Type \0\,:] \(\mc{M}_{C}\) is isomorphic to \(\mc{M}_{C'}\) and the derived equivalence \(\FM{C'}^{-1}\circ \FM{C}\) is a Seidel--Thomas twist~\cite{Seidel:bga} (up to a tensor product with a line bundle).
 \item[Type \I\,:] \(\mc{M}_{C}\) is a flop of \(\mc{M}_{C'}\) and \(\FM{C'}^{-1}\circ \FM{C}\) is the induced derived equivalence described by Bondal--Orlov~\cite[\S3]{Bondal:sod}.
 \item[Type \III:] \(\mc{M}_{C}\) is isomorphic to \(\mc{M}_{C'}\) and \(\FM{C'}^{-1}\circ \FM{C}\) is an \(EZ\)-transformation of Horja--Szendr\H{o}i~\cite{Horja:dca,Szendroi:ffm} (up to a tensor product with a line bundle).
 \end{description}
 \end{theorem}

 As Horja~\cite[\S4.1]{Horja:dca} observes,  both types of derived self-equivalence listed here are EZ-transformations (up to a tensor product with a line bundle).
 It would be interesting to discover the extent to which these transformations determine the group of self-equivalences of the derived category of \(\mc{M}_{C}\) (together with shifts,  automorphisms of \(\mc{M}_{C}\) and tensoring by line bundles).

 To prove Theorem~\ref{thm:main} we start from a chamber \(C\) and change the parameter so that the projection \(\res(\fmstar{C}(\theta))\in F^{1}/F^{2}\) moves towards the desired boundary of the ample cone of \(Y = \mc{M}_{C}\).  If \(C\) does not have the appropriate wall of type \I\ we pass into an adjacent chamber \(C_{1}\) separated from \(C\) by a wall of type \0.  A property of Seidel--Thomas twists reveals that \(\fmstar{C}(C)\) and \(\fmstar{C_{1}}(C_{1})\) are adjacent cones in \(F^{1}\),  modulo \(F^{2}\).  We proceed in this way towards the desired wall of the ample cone.  The main result of \S\ref{sec:maintheorem} establishes that after crossing finitely many walls of type \0\ we reach a chamber \(C_{n}\) defining \(Y\) with the desired wall of type \I. Thus,  by crossing the wall we induce the desired flop \(Y = \mc{M}_{C} \dashrightarrow \mc{M}_{C'} = Y'\).

 One consequence of Theorem~\ref{thm:main} is that every crepant resolution of \(\C^{3}/G\) comes armed with tautological bundles,  so \(\ghilb\) is no longer distinguished in this sense.  This should be a comfort to 3-folders:  no one minimal model should be privileged above all others.  Nevertheless,   the chamber \(C_{0}\) defining \(\ghilb = \mc{M}_{C_{0}}\) enjoys several very nice properties as we show in \S\ref{sec:ghilb}.
 For example,  every birational contraction of type \I\ or \III\ from \(\ghilb\) is achieved by crossing a wall of \(C_{0}\);  this is false in general,  see Example~\ref{ex:3stepflop}.  Moreover,  an inequality of the form \(\theta(\fm{C_{0}}(\taut{\rho}^{-1}\vert_{D})) > 0\) is a type \0\ wall of \(C_{0}\) if and only if \(D\) is compact irreducible divisor and \(\rho\) marks \(D\) according to Reid's recipe~\cite{Reid:mc,Craw:emc}.

 The commutative square described above exists without the assumption that \(G\) is Abelian,  so our approach may generalise to all finite subgroups \(G\subset \SL(3,\C)\).  However,  we use the Abelian assumption at several stages:
 \begin{itemize}
 \item The decomposition \(R = \Ro \oplus \Rt\) of the regular representation into sub and quotient representations is uniquely determined by a wall.
  \item The reducedness of the fibre \(\tau^{-1}(0)\) in Lemma~\ref{lemma:quotinequal} and Corollary~\ref{coro:quotinequal},  and the reducedness in Lemma~\ref{lemma:reduced}. 
 \item The rigidity result of \S\ref{sec:repsquiver}. This is essential for describing walls of type \0\ in terms of Fourier--Mukai transforms.
 \end{itemize}

 \noindent \textbf{Conventions} If \(W\) and \(\rho\) are representations of \(G\) with \(\rho\) irreducible,  we write \(\rho\subseteq W\) whenever \(\Hom_{G}(\rho,W)\neq 0\).  The terms `\(G\)-sheaf' or `\(G\)-equivariant sheaf' on \(\C^{n}\) are shorthand for a sheaf on \(\C^{n}\) with a compatible \(G\)-action as in \cite[\S4]{Bridgeland:mim}.

 \bigskip

 \noindent \textbf{Acknowledgements} The authors wish to thank H.~Nakajima and M.~Reid for introducing them to this problem.  Sincere thanks also to A.~King and T.~Logvinenko for many useful discussions, and to the Isaac Newton Institute for their hospitality while part of this paper was written.


 \section{Moduli of \protect\(G\)-constellations}

 This section discusses the moduli spaces \(\mc{M}_{\theta}\) of \(\theta\)-stable \(G\)-constellations for a finite subgroup \(G\subset \GL(n,\C)\),  generalising the moduli space of \(G\)-clusters.
 We also introduce the Fourier-Mukai transform \(\FM{\theta}\colon D(\mc{M}_{\theta}) \to D^G(\C^n)\) that plays a key role in this paper.
 In the special case \(G \subset \SL(n, \C)\) with \(n=2,3\) and for a generic choice of the GIT parameter \(\theta\), the method of Bridgeland, King and Reid~\cite{Bridgeland:mim} shows that \(\mc{M}_{\theta}\) is a crepant resolution of \(\C^{n}/G\) and \(\FM{\theta}\) is an equivalence.  
 This generalisation was well known to those authors.

 \subsection{Construction of the moduli \protect\(\mc{M}_{\theta}\protect\)}
 \label{sec:construction}
 For a finite subgroup \(G\subset \GL(n,\C)\),  write \(\Irr(G)\) for the set of equivalence classes of irreducible representations of \(G\).  Let \(R=\oplus_{\rho \in \Irr(G)} R_{\rho} \otimes \rho\) denote the regular representation of \(G\) and \(R(G) = \oplus_{\rho\in \Irr(G)}\Z\cdot \rho\) the representation ring.

 \begin{definition}
 A \emph{\(G\)-constellation} is a \(G\)-equivariant coherent sheaf \(F\) on \(\C^{n}\) such that \(H^{0}(F)\) is isomorphic as a \(\C[G]\)-module to \(R\).  Set 
 \[
 \Theta \! := \{\, \theta\in \Hom_{\Z}(R(G),\Q)\st \theta(R) = 0\,\}.
 \]
 For \(\theta\in \Theta\), a \(G\)-constellation is said to be \emph{\(\theta\)-stable} if every proper \(G\)-equivariant coherent subsheaf \(0 \subset E \subset F\) satisfies \(\theta(E) > 0\),  i.e., \(\theta(H^{0}(E)) > 0 = \theta(H^{0}(F))\).  The notion of \emph{\(\theta\)-semistable} is the same with \(\geq\) replacing \(>\).
 \end{definition}

 We consider moduli spaces \(\mc{M}_{\theta}\) of \(\theta\)-stable \(G\)-constellations and \(\overline{\mc{M}_{\theta}}\) of \(\theta\)-semistable \(G\)-constellations. 
The above definition of stability is a direct
translation of that for quiver representations introduced by
King~\cite{King:mor} into the language of sheaves on \(\C^n\) via the
observation by Ito--Nakajima~\cite[\S3]{Ito:mc3}.  We now recall the GIT construction of the moduli spaces (compare \cite{Sardo-Infirri:ros}).  It is convenient to set \(V = \C^{n}\).  Consider the affine scheme
 \[
\mc N=\{\, B \in \Hom_{\C[G]}(R, V \otimes R)\st B\wedge B=0 \,\},
 \]
 where \(B \wedge B\) lies in \(\Hom_{\C[G]}(R,\wedge^2 V  \otimes R)\).  Each map \(B\in \mc{N}\) determines an action of \(V^{*}\) on \(R\) and the condition \(B \wedge B=0\) ensures this action is commutative.  In this way, \(B\) endows \(R\) with a \(\C[V]\)-module structure,  where \(\C[V] = \oplus_{k=0}^{\infty} \text{Sym}^k(V^*)\) is the polynomial ring of functions on \(V\).  Thus \(B\) determines a \(G\)-constellation on \(V\).

 Denote by \(\Aut_{\C[G]}(R) \subset \GL(R)\) the group of \(G\)-equivariant automorphisms of \(R\).  It decomposes as \(\Aut_{\C[G]}(R)= \prod_{\rho \in \Irr(G)} \GL(R_{\rho})\).  This group acts on \(\mc{N}\) by conjugation and the diagonal scalar subgroup \(\C^{*}\) acts trivially,  leaving a faithful action of 
 \[
 \Paut_{\C[G]}(R):= \Aut_{\C[G]}(R)/\C^{*}
 \]
 on \(\mc{N}\) by conjugation. Let \(\theta\in \Hom_{\Z}(R(G),\Z)\) satisfy \(\theta(R) = 0\). Then \(\theta\) determines a character \(\chi_{\theta}\) of \(\Paut_{\C[G]}(R)\) mapping \([(g_{\rho})_{\rho}]\) to \(\prod_{\rho} \det(g_{\rho})^{\theta(\rho)}\in\C^{*}\). This character determines a linearisation of the trivial line bundle on \(\mc{N}\).  We denote by
 \[
 \mc{M}_{\theta}\! := \mc{N}_{\theta}/\Paut_{\C[G]}(R)\quad\text{and}\quad\overline{\mc{M}_{\theta}} \! := \overline{\mc N_{\theta}}\git\Paut_{\C[G]}(R)
 \]
 the geometric and categorical quotients of the open subsets \(\mc{N}_{\theta}\) and \(\overline{\mc{N}_{\theta}}\) of \(\chi_{\theta}\)-stable and \(\chi_{\theta}\)-semistable points of \(\mc{N}\).
 King's result (which holds in a more general context) asserts that \(\mc{M}_{\theta}\) and \(\overline{\mc{M}_{\theta}}\) are equivalently the moduli spaces of \(\theta\)-stable and \(\theta\)-semistable \(G\)-constellations respectively.
 There is an isomorphism \(\mathcal{M}_{k\theta}\cong \mathcal{M}_{\theta}\) for \(k\in \Z_{> 0}\),  so \(\mathcal{M}_{\theta}\) is well defined even for parameters  \(\theta\in \Hom_{\Z}(R(G),\Q)\),  i.e.,  for \(\theta\in \Theta\).

 To let \(R \otimes \mc{O}_{\mc{N}}\) and the universal homomorphism \(R \otimes \mc O_{\mc N} \rightarrow V \otimes R \otimes \mc O_{\mc N}\) on \(\mc{N}\) descend to \(\mc{M}_{\theta}\),  we have to determine an equivariant \(\Paut_{\C[G]}(R)\)-action on \(R \otimes \mc{O}_{\mc{N}}\).  Define a subgroup 
 \[
 \Aut'_{\C[G]}(R) = \{\,\textstyle{ \prod_{\rho \in \Irr(G)} (g_{\rho}) \in \Aut_{\C[G]}(R) \st g_{\rho_0}=1} \,\}
 \]
 of \(\Aut_{\C[G]}(R)\),  where $\rho_0$ is the trivial representation and hence \(\dim R_{\rho_0} = 1\).
 Then the projection gives an isomorphism \(\Aut'_{\C[G]}(R) \to \Paut_{\C[G]}(R)\).
 On the other hand, \(\Aut'_{\C[G]}(R)\) has a natural equivariant action on \(R \otimes \mc O_{\mc N}\).
 Descent theory of coherent sheaves for a faithfully flat morphism (see \cite[VIII]{Grothendieck:sga1}) implies that \(R \otimes \owe_{\mathcal{N}}\) and \(R_{\rho} \otimes \owe_{\mathcal{N}}\) descend to locally free sheaves \(\mc{R}\!:= \mc{R}_{\theta}\) and \(\taut{\rho}\!:= (\mc{R}_{\theta})_{\rho}\) respectively on \(\mc{M}_{\theta}\) such that \(\mc{R} = \oplus_{\rho \in \Irr(G)} \mc{R}_{\rho} \otimes \rho\).
 Moreover, we have a homomorphism \(\mc{R} \to V \otimes \mc{R}\) determining the universal \(G\)-constellation \(\mc{U}_{\theta}\) on \(\mc M_{\theta}\times \C^{n}\).
 Note by the definition of \(\Aut'_{\C[G]}(R)\) that the line bundle \(\mc{R}_{\rho_{0}}\) is trivial.

 \subsection{The morphism to \protect\(X = \C^{n}/G\protect\)}
 A parameter \(\theta \in \Theta\) is \emph{generic} if every \(\theta\)-semistable \(G\)-constellation is \(\theta\)-stable.

 \begin{proposition}
 \label{prop:maptox}
 The moduli spaces and the singularity \(X = \C^{n}/G\) are related in the following way.
 \begin{enumerate}
 \item There is a closed immersion \(X \to \overline{\mc M_0}\) such that an orbit \(G\cdot x\) is mapped to the class of \(G\)-constellations whose support is \(G\cdot x\).
This makes $X$ an irreducible component of $\overline{\mc M_0}$.
 \item There is a morphism \(\tau\colon \mc{M}_{\theta} \to X\) that associates a \(G\)-constellation to its support. If \(\theta\) is generic then this morphism is projective.
 \end{enumerate}
 \end{proposition}
 \begin{proof}
 The morphism \(G \times \C^n \to \C^n \times \C^n\) sending \((g, x)\) to \((g\cdot x, x)\) is \(G\)-equivariant with respect to the action of \(G\) on the first factor, so the second projection makes \(\C[G] \otimes \owe_{\C^n}\) into a flat family of \(G\)-constellations parametrised by \(\C^n\).
 Since every \(G\)-constellation is 0-semistable
 and since $\overline{\mc{M}_{0}}$ {\it corepresents} the
 moduli functor (see Huybrechts and Lehn~\cite{Huy}), this family gives rise to a morphism \(\C^n \to \overline{\mc M_0}\) which is \(G\)-equivariant with respect to the trivial action of \(G\) on \(\overline{\mc{M}_{0}}\),  and hence to a morphism \(X \to \overline{\mc{M}_{0}}\).
 To see that this is a closed immersion, consider the morphism from \(\overline{\mc M_0}\) to the \(G\)-fixed part \((S^{\#G}(\C^n))^G\) of the symmetric product of \(\C^n\), similar to \cite[Example 4.3.6]{Huy}.
Now,  \(X\) is the Chow quotient of \(\C^n\) by \(G\) in the sense of Kapranov~\cite[Theorem 0.4.3]{Kapranov:cqg},  so the composite morphism \(X \to (S^{\#G}(\C^n))^G\),  and hence also \(X \to \overline{\mc M_0}\),  is a closed immersion.
 
 Since \(X\) is affine,  the morphism \(\tau\) in the second assertion is determined by a ring homomorphism \(H^0(\owe_X) \to H^0(\owe_{\mc{M}_{\theta}})\) that we now construct.  
 Given a \(G\)-invariant polynomial \(f\),  multiplication by \(\pi_{\C^n}^*f\) determines an endomophism of the universal \(G\)-constellation \(\mc{U}_{\theta}\) on \(\mc{M}_{\theta} \times \C^{n}\).
 Since \(\mc{U}_{\theta}\) is a flat family of \(\theta\)-stable \(G\)-constellations,  the base change theorem for \(\GExt\) gives rise to an isomorphism \(\pi_{\mc M_{\theta}*}(G\text{-}\mc Hom_{\owe_{\mc M_{\theta} \times\C^n}}(\mc{U}_{\theta}, \mc{U}_{\theta})) \cong \owe_{\mc M_{\theta}}\). 
 It follows that \(\pi_{\C^n}^*f\in \GHom_{\owe_{\mc M_{\theta} \times\C^n}}(\mc{U}_{\theta},\mc{U}_{\theta})\) determines a function on \(\mc{M}_{\theta}\) as required.
  If \(\pi_{\C^n}^{*}f\) vanishes on the support of a \(G\)-constellation \(F\) then multiplication by \(f\) determines a noninvertible endomorphism of \(F\). But this must be zero by the stability of \(F\), so \([F] \in \mc{M}_{\theta}\) is mapped to the support of \(F\).
 Thus we obtain the set-theoretic description of \(\tau\).
 Moreover, the diagram
 \[
 \begin{CD}
 \mc M_{\theta} @>>> \overline{\mc M_{\theta}}\\
 @VVV                                   @VVV \\
 X  @>>> \overline{\mc M_0}
 \end{CD}
 \]
 is set-theoretically commutative.  The right-hand vertical arrow in this diagram is a projective morphism (see, for example,  King~\cite[\S2]{King:mor}).  If \(\theta\) is generic then the top horizontal arrow is an isomorphism,  so \(\mc{M}_{\theta}\) is projective over \(\overline{\mc{M}_{0}}\) and hence over the irreducible component \(X\) via the morphism \(\tau\).
 \end{proof}

 \begin{remark}
 If \(G\) acts freely outside the origin then the morphisms \(X \to \overline{\mc{M}_0}\) and \(\overline{\mc{M}_0} \to (S^{\#G}(\C^n))^G\) are bijective as maps of sets.
 Otherwise, \(\overline{\mc{M}_0}\) has an irreducible component other than \(X\) because the support of a \(G\)-constellation can be more than one orbit.
 Furthermore, \((S^{\#G}(\C^n))^G\) can have more components consisting of the supports of \(G\)-sheaves that are not \(G\)-constellations.
 \end{remark}

 \subsection{The moduli space of \protect\(G\protect\)-clusters}
 \label{sec:chamberghilb}
 Recall that a \emph{\(G\)-cluster} is a \(G\)-invariant subscheme \(Z\subset\C^{n}\) of dimension zero with global sections \(H^{0}(\mathcal{O}_{Z})\) isomorphic as a \(\C[G]\)-module to the regular representation \(R\) of \(G\).  
 Write \(\Ghilb{n}\),  or simply \(\ghilb\),  for the moduli space of \(G\)-clusters.  

 Ito and Nakajima~\cite[\S3]{Ito:mc3} observed that \(\ghilb\) is isomorphic to \(\mathcal{M}_{\theta}\) for parameters \(\theta\) in the cone
 \begin{equation}
 \label{eqn:Theta+}
 \Theta_{+} \!:= \{\theta \in \Theta \st \theta(\rho) > 0 \text{ if } \rho \neq \rho_{0} \}.
 \end{equation}
 To see this,  let \(F\) be a \(\theta\)-stable \(G\)-constellation for \(\theta\in \Theta_{+}\). 
 No proper submodule of \(F\) contains the trivial representation so \(F\) is a cyclic \(\owe_{\C^{n}}\)-module with generator \(\rho_{0}\),  hence \(F\cong \owe_{Z}\) for some \(Z\in \ghilb\) as required.   Note also that the fibre of the bundle \(\taut{} = \pi_{*}\owe_{\mathcal{Z}}\) over a point \(Z\in \ghilb\) is \(H^{0}(\owe_{Z})\).  
 \begin{remark}
 In general,  \(\ghilb = \mathcal{M}_{\theta}\) for a set (in fact a cone) of parameters that \emph{strictly} contains \(\Theta_{+}\),  see Example~\ref{ex:11.1.2.8}
 \end{remark}

 \subsection{Fourier--Mukai transforms}
 \label{sec:bkr}
 Set \(Y = \mathcal{M}_{\theta}\) for some \(\theta \in \Theta\).
 The morphism introduced in Proposition~\ref{prop:maptox} fits into the following commutative diagram:

 \begin{equation*}\label{maindiag}
 \setlength{\unitlength}{36pt}
 \begin{picture}(2.2,2)(0,0)
 \put(1,0){\object{$X$}}
 \put(0,1){\object{$Y$}}
 \put(2,1){\object{$\C^n$}}
 \put(1,2){\object{$Y\times \C^n$}}
 \put(0.25,0.75){\vector(1,-1){0.5}}
 \put(0.5,0.5){\nwlabel{$\tau$}}
 \put(1.75,0.75){\vector(-1,-1){0.5}}
 \put(1.5,0.5){\swlabel{$\pi$}}
 \put(0.75,1.75){\vector(-1,-1){0.5}}
 \put(0.5,1.5){\nelabel{$\pi_{Y}$}}
 \put(1.25,1.75){\vector(1,-1){0.5}}
 \put(1.5,1.5){\selabel{$\pi_{\C^n}$}}
 \end{picture}
 \end{equation*}

 The tautological bundle \(\taut{}\) on \(Y\) satisfies \(\taut{} = {\pi_Y}_* \mathcal{U}_{\theta} = \oplus_{\rho} (\taut{\rho} \otimes \rho)\),  where \(\mathcal{U}_{\theta}\) is the universal \(G\)-constellation on \(Y\times \C^{n}\).  Let \(D(Y)\) and \(D^{G}(\C^n)\) denote the bounded derived categories of coherent sheaves on \(Y\) and \(G\)-equivariant coherent sheaves on \(\C^n\) respectively.  Consider the functor \(\FM{\theta} \colon D(Y) \to D^{G}(\C^n)\) defined by
 \[
 \Phi_{\theta}(-) = {\RR}{\pi_{\C^n}}_*(\mathcal{U}_{\theta}\otimes {\pi_Y}^*(- \otimes \rho_0)).
 \]
 We repeatedly use the following formula to calculate \(\FM{\theta}\) (we sometimes abuse notation and omit $\pi_*$ from the left hand side):
 \begin{equation}
 \label{eqn:FMi}
 \pi_* \FM{\theta}^{i}(-) \cong R^i \tau_*(- \otimes \taut{})
 = \bigoplus_{\rho\in \Irr(G)} H^{i}(-\otimes \taut{\rho})\otimes \rho.
 \end{equation}
 The method of Bridgeland, King and Reid~\cite{Bridgeland:mim} generalises to prove the
 following result.

 \begin{theorem}
 Let \(G\subset \SL(3,\C)\) be a finite subgroup.  If \(\theta\) is generic then \(\tau\colon \mc M_{\theta} \to X\) is a crepant resolution and \(\FM{\theta}\) is an equivalence of triangulated categories.
 \end{theorem}

 The (quasi-)inverse \(\Phi_{\theta}^{-1} \colon D^{G}({\C^3}) \to D(Y)\) given by
 \[
 \FM{\theta}^{-1}(-) = \big{[}{{\RR}\pi_Y}_*(\mathcal{U}_{\theta}^{\vee}[3] \overset{\LL}{\otimes} {\pi_{\C^3}}^*(-))\big{]}^G,
 \]
 where \(\mathcal{U}_{\theta}^{\vee}\) denotes the derived dual \(\RR\mc{H}om_{\owe_{Y\times {\C^3}}}(\mathcal{U}_{\theta},\owe_{Y\times {\C^3}})\).  
Then
 \[
 \FM{\theta}^{-1}(\rho \otimes \owe_{\C^3})
 = \big{[}(\RR {\pi_Y}_*(\mc U_{\theta}^{\vee}[3])) \otimes \rho\big{]}^G
 \cong (\taut{}^{\vee} \otimes \rho)^G = \taut{\rho}^{\vee}.
 \]
 Equivalently, we have \(\FM{\theta}(\taut{\rho}^{\vee}) \cong \rho \otimes \owe_{\C^3}\).

 The restriction of \(\FM{\theta}\) to the full subcategory of \(D(Y)\) consisting of objects supported on the subscheme \(\tau^{-1}(\pi(0))\) induces an isomorphism between the Grothendieck group \(K_{0}(Y)\) of coherent sheaves supported on \(\tau^{-1}(\pi(0))\) and the Grothendieck group \(K_{0}^{G}(\C^{3})\) of \(G\)-equivariant coherent sheaves on \(\C^{3}\) supported at the origin.  Thus,  just as in \cite[\S9]{Bridgeland:mim},  we obtain an isomorphism
 \begin{equation}
 \label{eqn:fmktheory}
 \fm{\theta}\colon K_0(Y)\longrightarrow K_{0}^{G}(\C^{3}).
 \end{equation}
 The equivalence \(\FM{\theta}\) also induces an isomorphism \(K(Y) \to K^G(\C^3)\).  Since \(\FM{\theta}(\taut{\rho}^{\vee}) \cong \rho \otimes \owe_{\C^3}\),  the tautological bundles \(\taut{\rho}\) form a \(\Z\)-basis of \(K(Y)\).  The notation \(\fm{\theta}\) (and \(\fm{C}\) later,  see \S\ref{sec:ThetaFM}) will always denote the isomorphism (\ref{eqn:fmktheory}).

 \subsection{Symmetries of \protect\(\Theta\protect\)}
 \label{sec:symmetries}
 The parameter space \(\Theta\) and the moduli spaces admit the following symmetries.

 \begin{lemma}
 \label{lemma:symmetry}
 \begin{enumerate}
 \item[\one] 
 Fix a one-dimensional representation \(\sigma \in \Irr(G)\) and define \(\theta'\in \Theta\) by \(\theta'(\rho)\! := \theta(\sigma\rho)\). 
 Then \(\mathcal{M}_{\theta'} \cong \mathcal{M}_{\theta}\), and the pushforward of the universal \(G\)-constellation to \(\mc{M}_{\theta'}\) is \(\oplus_{\rho} (\taut{\rho\sigma} \otimes \taut{\sigma}^{-1} \otimes \rho)\).
 The Fourier--Mukai functors are related by
 \(\FM{\theta'}(\alpha)=\sigma \otimes \FM{\theta}(\taut{\sigma}^{-1}\otimes \alpha)\).
 \item[\two] 
 Define \(\theta'\in \Theta\) by \(\theta'(\rho)\! := -\theta(\rho^{*})\).  
 Then \(\mathcal{M}_{\theta'} \cong \mathcal{M}_{\theta}\) and the pushforward of the universal \(G\)-constellation to \(\mc{M}_{\theta'}\) is \(\oplus_{\rho} \taut{\rho^*}^{\vee} \otimes \rho\).
 The Fourier--Mukai functors are related by \(\FM{\theta'}(\alpha)=\FM{\theta}(\alpha^{\vee})^{\vee}\).
 \end{enumerate}
 \begin{proof}
 For \one, the isomorphism is induced by the map \(F \mapsto F \otimes_{\C} \sigma^*\).
 The universal \(G\)-constellation for \(\theta'\) is \(\mc{U}_{\theta} \otimes \pi_{Y}^*(\taut{\sigma}^{-1}) \otimes_{\C} \sigma^*\),  where \(\mc{U}_{\theta}\) is the universal \(G\)-constellation for \(\theta\) (the term \(\taut{\sigma}^{-1}\) appears because of our normalisation of the universal family).
 This implies the statements for the tautological bundles and the Fourier--Mukai transforms in \one.  As for \two,  the isomorphism is induced by the map \(F \mapsto F^{\vee}[3] = \mc{E}xt^3_{\owe_{\C^3}}(F, \owe_{\C^3})\).
 The universal \(G\)-constellation for \(\theta'\) is \(\mc{U}_{\theta}^{\vee}[3] = \mc{E}xt^3_{\owe_{Y \times \C^3}}(\mc{U}_{\theta}, \owe_{\mc M_{\theta} \times \C^3})\).
 Grothendieck duality for \(\pi_{Y}\) gives \({\pi_{Y}}_*(\mc{U}_{\theta}^{\vee}[3]) \cong ({\pi_{Y}}_*\mc{U}_{\theta})^{\vee}\).  This establishes the statement for the tautological bundles.
 Grothendieck duality for \(\pi_{\C^3}\) gives
 \[
\FM{\theta'}(\alpha)=\RR{\pi_{\C^3}}_*
(\mc U_{\theta}^{\vee}[3] \otimes \pi_{Y}^*(\alpha \otimes \rho_0))
\cong
\big{(}\RR{\pi_{\C^3}}_*
(\mc U_{\theta} \otimes \pi_{Y}^*(\alpha^{\vee} \otimes \rho_0))\big{)}^{\vee}.
 \]
 This proves the final statement.
 \end{proof}
 \end{lemma}


 \section{Classifying walls in \protect\(\Theta\protect\)} 
 \label{sec:classify-walls}
 This section defines and classifies the walls of the parameter space \(\Theta\) for a finite subgroup \(G\subset \SL(3,\C)\).   The contraction morphism determined by a wall crossing and the unstable locus of a wall crossing are defined under the additional assumption that \(G\) is Abelian.

 \subsection{Chambers and walls in \protect\(\Theta\protect\)}
 
 \begin{lemma}
 The subset \(\Theta^{\gen}\subset \Theta\) of generic parameters is open and dense.
 It is the disjoint union of finitely many open convex polyhedral cones in \(\Theta\).
 \end{lemma}
 \begin{proof}
 Let \(\theta\) be generic and put
 \begin{eqnarray*}
 C & = & \{\, \eta \in \Theta \st \text{every \(\theta\)-stable \(G\)-constellation is \(\eta\)-stable} \,\} \\
   & = & \left\{\eta \in \Theta  \left|\begin{array}{l} \eta(F)>0 \text{ for every nontrivial subsheaf } \\ F \text{ of every } \theta\text{-stable } G\text{-constellation} \end{array}\right.\right\}.
 \end{eqnarray*}
 This is a convex polyhedral cone. We claim that every \(\eta \in C\) is generic and every \(\eta\)-stable \(G\)-constellation is \(\theta\)-stable.  Indeed,  if \(\eta \in C\), then \(\mc{M}_{\eta}\) contains \(\mc{M}_{\theta}\) as an open set by definition.  Then \(\mc{M}_{\theta}\) must be a connected component of \(\mc{M}_{\eta}\), so the arguments in \cite[\S8]{Bridgeland:mim} imply that \(\mc{M}_{\theta} = \overline{\mc{M}_{\eta}}\).
 \end{proof}
 
 \begin{definition}
 A \emph{chamber} in \(\Theta\) is a convex polyhedral cone \(C\) as in the above lemma.
 A codimension-one face \(W\) of \(\overline{C}\) is a \emph{wall} of \(C\).
 Write \(\mc{M}_{C}\), \(\taut{C}\) and \(\mc{U}_{C}\) for the moduli space \(\mc{M}_{\theta}\),  the tautological bundle \(\taut{\theta}\) and the universal \(G\)-constellation \(\mc{U}_{\theta}\) respectively,  for any \(\theta \in C\).
 \end{definition}
 
 \subsection{Classification of walls}
 \label{sec:classificationofwalls}
 Let \(\theta, \theta' \in \Theta\) be generic parameters and \(C, C'\) the chambers containing them.  Assume that if we put \(\theta_{t} = \frac{1}{2}(1-t)\theta + \frac{1}{2}(1+t)\theta'\),  then \(\theta_{t} \in C'\) when \(0 < t \leq 1\) and \(\theta_{t} \in C\) when \(-1 \leq t < 0\).   We are interested in the case where \(C\) and \(C'\) are adjacent chambers that share a wall.  To simplify notation,  we write
\[
 \taut{} = \bigoplus_{\rho \in \Irr(G)} \taut{\rho}\otimes \rho \quad \text{and}\quad \taut{}' = \bigoplus_{\rho \in \Irr(G)} \taut{\rho}'\otimes \rho
 \]
 for \(\taut{} = \taut{C}\) and \(\taut{}' = \taut{C'}\).  For any \(\eta \in \Theta\), the (fractional) line bundles
 \[
 \lb{C}(\eta):=\bigotimes_{\rho \in \Irr(G)} \left(\det \taut{\rho}\right)^{\eta(\rho) }
 \quad\text{and}\quad
 \lb{C'}(\eta):=\bigotimes_{\rho \in \Irr(G)} \left(\det {\taut{\rho}'}\right)^{\eta(\rho)}
   \]
 are the descendents of the restrictions of \(\owe_{\mc{N}}\) linearised by \(\eta\).
 In particular,  \(\lb{C}(\theta)\) and \(\lb{C'}(\theta')\) are the ample line bundles \(\owe_{\mc{M}_{\theta}}(1)\) and \(\owe_{\mc{M}_{\theta'}}(1)\) respectively,  canonically constructed by GIT.
 There is also a (fractional) ample line bundle \(L_{0} = \owe_{\overline{\mc{M}_{\theta_0}}}(1)\) on \(\overline{\mc{M}_{\theta_0}}\).  

 Since the tautological bundles \(\taut{\rho}\) generate \(K(\mc{M}_{C})\),  their first Chern classes generate \(\Pic(\mc{M}_{C})\) so the map \(L_{C}\colon \Theta\to\Pic(\mc{M}_{C})_{\Q}\) is surjective.  The following lemma holds in a more general situation,  see Thaddeus~\cite{Thaddeus:git}.  

 \begin{lemma}
 \label{lemma:GIT}
 Write \(f\colon \mc{M}_{C} \to \overline{\mc M_{\theta_0}}\) and \(f'\colon \mc{M}_{C'} \to \overline{\mc M_{\theta_0}}\) for the canonical morphisms.  Then \(\lb{C}(\theta_0) \cong f^* L_0\) and \(\lb{C}(\theta) \otimes \lb{C}(\theta') \cong f^{*} L_{0}^{\otimes 2}\),  so \(\lb{C}(\theta')^{-1}\) is \(f\)-ample.  
 Moreover, if \(\lb{C}(\theta_0)\) is ample on \(\mc{M}_{C}\) then \(f,f'\) are injective and \(\mc{M}_{C} \cong \mc M_{C'}\).
 \end{lemma}
 \begin{proof}
 The first assertions follow from uniqueness of descendents for a faithfully flat morphism.    Assume \(\lb{C}(\theta_0)\) is ample,  so \(\lb{C}(\theta')\) is ample.  Then \(\lb{C}(\theta')^{-1}\) is also $f$-ample and hence $f$ must be finite.
 The fibre consists of $G$-constellations in the same S-equivalence class (see \cite{Huy} for S-equivalence).
 Thus it is connected and therefore must be a single point.  Since \(\mc{M}_{C'}\) is also a crepant resolution,  we obtain the result.
 \end{proof}

 \begin{definition}
 \label{defn:type}
 If \(\theta_{0}\) is a general point in a wall \(W\) of \(C\),  we write \(\cont{W}\) for the morphism from \(\mc{M}_{C}\) onto (the normalisation of) its image under the morphism \(f\colon \mc M_{C} \to \overline{\mc M_{\theta_0}}\).  Note that \(\overline{\mc{M}_{\theta_0}}\) is not necessarily irreducible.
 \end{definition}

 The morphism \(\cont{W}\) is determined by \(\lb{C}(\theta_{0})\).  Since \(\theta_0\) lies in the closure of \(C\), the line bundle \(\lb{C}(\theta_{0})\) is nef.  If \(\lb{C}(\theta_{0})\) is in fact ample then \(\cont{W}\) is an isomorphism by Lemma~\ref{lemma:GIT}.  Otherwise \(\lb{C}(\theta_{0})\) corresponds to a class on the boundary of the ample cone of \(\mc{M}_{C}\) and, since \(\theta_{0}\) is general in the wall \(W\),  this class lies in the interior of a facet of the ample cone.  Thus the corresponding contraction \(\cont{W}\) is primitive,  i.e.,  it cannot be further factored into birational morphisms between normal varieties.

 \begin{definition}
 The wall $W$ is said to be of type \0, \I, \II, \III\ as follows:
 \begin{itemize}
 \item type \0\ if \(\cont{W}\) is an isomorphism.
 \item type \I\  if \(\cont{W}\) contracts a curve to a point.
 \item type \II\  if \(\cont{W}\) contracts a divisor to a point.
 \item type \III\  if \(\cont{W}\) contracts a divisor to a curve.
 \end{itemize}
 \end{definition}

 \begin{remark}
 The terminology `type \I, \II, \III' is standard for 3-fold birational contractions,  see Wilson~\cite{Wilson:kahler}.
 Proposition~\ref{prop:notypeII} to follow asserts that type \II\ walls do not exist. 
 \end{remark}

 \subsection{The fibre of \protect\(\cont{W}\protect\)}
 \label{sec:fibreofcont}
 Hereafter we assume \(G\subset \SL(3,\C)\) is a finite Abelian subgroup.  
 Fix a chamber \(C\),  suppose \(\theta\in C\) and \(\theta_{0}\) is on a wall of \(C\).
 Then there are nonzero representations \(\Ro, \Rt\) of \(G\) such that 
 \begin{enumerate}
 \item[\one] \(R \cong \Ro \oplus \Rt\); and
 \item[\two] there exist \(\theta_{0}\)-semistable \(G\)-equivariant coherent sheaves \(S\) and \(Q\) that satisfy \(H^{0}(S) \cong\Ro\) and \(H^{0}(Q) \cong\Rt\) as \(\C[G]\)-modules,  such that \(\theta_{0}(S) = \theta_{0}(Q) = 0\).
 \end{enumerate}
 Assume $\theta_{0}$ is general in the wall so that it is not contained in two distinct hyperplanes in $\Theta$ defined by proper subrepresentations of $R$.
 Since \(G\) is Abelian, every irreducible representation is of multiplicity one in $R$ and therefore \(\Ro\) and \(\Rt\) are the only representations of \(G\) that can be embedded in \(R\) as proper subrepresentations with this property. 
 In particular, both $S$ and $Q$ above are $\theta_0$-stable.
 We assume without loss of generality that \(\theta(S) = \theta(\Ro) > 0\).

 \begin{lemma}
 \label{lemma:fibre}
 The fibre (with reduced structure) of the morphism $\cont{W}$ over the point $\cont{W}(F)$ is isomorphic to \(\mathbb{P}(\GExt_{\mc{O}_{\C^3}}^{1}(Q, S)^{\vee})\).
 \end{lemma}
 \begin{proof}
 The fibre $\cont{W}^{-1}(\cont{W}(F))$ consists of $G$-constellatiions with the same S-equivalence class as $F$,  where S-equivalence is considered with respect to $\theta_0$.
 Thus every $G$-constellation in the fibre is an extension of $Q$ by $S$.
 Conversely, every such extension is $\theta$-stable and the assertion is proved.
 \end{proof}

 \begin{proposition}
 \label{prop:notypeII}
 There are no walls of type \II.
 \end{proposition}
 \begin{proof}
 Assume \(\theta_{0}\) is a general point of a type \II\ wall of \(C\).  Suppose that a \(G\)-constellation \(F\) corresponds to a point on the contracted surface that lies on no toric curve,  and that \(S\subset F\) is a subsheaf with \(\theta_{0}(S) = 0\).
 Since \(\cont{W}\) contracts a surface to a point,  we have \(\dim \GExt^{1}_{\mc{O}_{\C^{3}}}(F/S, S) = 3\) by Lemma \ref{lemma:fibre}.  This contradicts the technical result that we postpone to Proposition~\ref{prop:rigidity}\one.
 \end{proof}

 \subsection{The unstable locus of a wall crossing}
 \label{sec:newunstablelocus}
 We use the same notation and consider the same situation
as in the previous subsection.
 The following lemma is straightforward.

 \begin{lemma}
 Let \(F\) be a \(\theta\)-stable \(G\)-constellation.
 Suppose there is a subsheaf \(S \subset F\) such that \(H^{0}(S) \cong \Ro\) as a \(\C[G]\)-module.
 Then,  we have the following.
 \begin{enumerate}
 \item \(S\) and \(Q = F/S\) are \(\theta_{0}\)-stable and \(\GHom_{\C^3}(S, Q) = 0\).
 \item If \(S' \subset F\) has the same property then \(S' = S\).
 \end{enumerate}
 \end{lemma}
 
 We now introduce a natural scheme structure on the following subset of \(\mc{M}_C\):
\[
 Z = \{\, F \in \mc M_{C} \st \exists\; S \subset F \text{ such that \(H^{0}(S) \cong \Ro\) as a \(\C[G]\)-module} \,\}.
 \]
 Let \(h\colon (\text{Schemes over $\C$})\to (\text{Sets})\) be the functor such that \(h(\mc T)\) is the set of \(f \in \Hom(\mc T, \mc M_{C})\) such that there exists a quotient \(\nu \colon(f\times \text{id}_{\C^3})^* \mc U \twoheadrightarrow \mc Q_{\mc{T}}\) flat over \(\mc T\) with \(H^0(\mc Q_t ) \cong R_2 \otimes k(t)\) for all \(t \in \mc T\).

 \begin{lemma}
 \label{lemma:subscheme}
 There exists a closed subscheme $Z$ of $\mc M_C$ representing the functor \(h\).
 \end{lemma}
 \begin{proof}
 Define a functor \(\widetilde h\colon (\text{Schemes over $\C$})\to (\text{Sets})\) by
 \[
 \widetilde h(\mc T)=\{ (f, \nu) | f \in \Hom(\mc T, \mc M_{C}) \text{ and $\nu$ is a quotient as above}\}.
 \]
 Then \(\widetilde h\) is represented by a closed subscheme $Z$ of Grothendieck's Quot-scheme
 \[
\Quot^{\dim \Rt}_{\uni{C}/(\mc{M}_{C} \times \C^{3})/\mc{M}_{C}}
 \]
 of quotients of restrictions of the universal sheaf \(\uni{C}\) to the fibres of \(\mc{M}_{C} \times \C^{3} \to \mc{M}_{C}\) (see \cite{Grothendieck:hilb}). 
The morphism \(\mc{M}_{C} \times \C^{3} \to \mc{M}_{C}\) is nonprojective, but the restriction to the support of the universal sheaf \(\Supp(\uni{C})\subset \mc{M}_{C} \times_{X} \C^{3}\) is projective,  so the above Quot scheme is well defined.  Now,  \(G\) acts on \(\uni{C}\) and acts equivariantly on  \(\mc{M}_{C}\times \C^{3}\) (we give \(\mc{M}_{C}\) the trivial action of \(G\)),  so \(G\) acts on the Quot scheme.  Observe that \({Z}\) is the union of several connected components of the \(G\)-fixed point locus and consists of quotients \(Q\) such that \(H^{0}(Q)\cong \Rt\) as a \(\C[G]\)-module.

 The forgetful map \((f, \nu) \mapsto f\) leads to a morphism \({Z} \to \mc{M}_{C}\).
 Part (2) of the above lemma shows that this morphism is injective as a map,  while part (1) asserts that the tangent map is also injective.  The structure morphism is therefore a closed embedding \({Z} \hookrightarrow \mc{M}_{C}\) and we obtain an isomorphism \(\widetilde h \to h\).  \end{proof}

 \begin{definition}
 The above subscheme $Z$ is the \emph{unstable locus} in \(\mc{M}_{C}\) defined by \(W\).
 \end{definition}

 \begin{remark}
 It will be shown in Propositions~\ref{prop:rigiddivisor},~\ref{prop:typeIflop} and~\ref{prop:typeIIIiso} that the unstable locus \(Z\) is always connected.
 \end{remark}

 Symmetrically,  there is a closed subscheme \(Z'\) of \(\mc{M}_{C'}\) such that
 \[
 Z' = \{\, F' \in \mc M_{C'} \st \exists\; S' \subset F' \text{ such that \(H^{0}(S')\cong \Rt\) as a \(\C[G]\)-module} \,\},
 \]
 Since both \(\mc{M}_{C} \setminus Z\) and \(\mc{M}_{C'}\setminus Z'\) parametrise \(G\)-constellations that are simultaneously \(\theta\)-stable and \(\theta'\)-stable, we have the following:

 \begin{lemma}
 \label{lemma:common}
 There is an isomorphism \(\mc{M}_{C} \setminus Z \cong \mc{M}_{C'} \setminus Z'\).  Furthermore,  the restriction of the universal \(G\)-constellation \(\taut{C}\) is isomorphic as a family of $G$-constellations to the restriction of \(\taut{C'}\).
 \end{lemma}

 A consequence of Lemma~\ref{lemma:subscheme} is the existence of a universal subsheaf \(\mc{S} \subset \taut{C}\vert_{Z}\) such that \(\mc{S}\) and \(\mc{Q} = (\taut{C}\vert_{Z})/\mc{S}\) are locally free \(\owe_{Z}\)-modules.  Similarly, there is a universal subsheaf \(\mc{S'} \subset \taut{C'}\vert_{Z'}\) and quotient \(\mc{Q}'\).


 \section{Crossing walls of type \0}
 \label{sec:type0}
 Throughout this section,  \(C\) and \(C'\) denote chambers separated by a wall \(W\) of type \0.  Recall that the moduli \(\mc{M}_{C}\) and \(\mc{M}_{C'}\) are isomorphic by Lemma~\ref{lemma:GIT}.  We investigate the unstable locus of the wall crossing and determine how the tautological bundles change as a parameter \(\theta\) crosses the wall \(W\).

 \subsection{Change of tautological bundles}
 \label{sec:type0taut}
 
 \begin{proposition}
 \label{prop:Cartier}
 For a wall of type \0,  the unstable locus \(Z = Z'\) is a Cartier divisor that we denote \(D\).
 \end{proposition}
 \begin{proof}
 It follows from Lemmas~\ref{lemma:GIT} and~\ref{lemma:common} that \(Z\) is set-theoretically identical to \(Z'\).  
 Assume that \(\codim Z \geq 2\).
 Then the universal homomorphism \(\taut{C} \to V \otimes \taut{C}\) and \(\taut{C}\) are uniquely determined by their restrictions to \(\mc{M}_{C}\setminus Z\).
 But this implies that the universal $G$-constellations on \(\mc{M}_{C}\) and \(\mc{M}_{C'}\) are isomorphic to each other,  contradicting the assumption that \(C\neq C'\).  Thus \(\codim Z = 1\).

 Let \(\mathfrak{Z} \subset Z\) be the union of the associated components of \(Z\) (i.e., the closures of the scheme-theoretic points determined by associated primes of $\owe_{Z}$) whose codimensions are greater than $1$.
 Let \(\mc K\) be the kernel of the surjection $\taut{C} \to \mc Q$,  the composite of the restriction to \(Z\) and the obvious surjection.
 Since \(Z\setminus \mathfrak{Z}\) is a Cartier divisor on \(\mc{M}_{C}\setminus \mathfrak{Z}\),  \(\mc{K}\) is torsion-free on \(\mc{M}_{C}\) and locally free on \(\mc{M}_{C} \setminus \mathfrak{Z}\).  An inverse elementary transformation (see Maruyama~\cite{Maruyama:avb}) establishes the following exact sequence of locally free sheaves on \(Z \setminus \mathfrak{Z}\):
 \begin{equation}
 \label{eqn:obstruction}
 0 \rightarrow \mc Q|_{(Z \setminus \mathfrak{Z})} \otimes \owe_{(\mc M_{C}\setminus \mathfrak{Z})}(-(Z \setminus \mathfrak{Z})) \rightarrow \mc K |_{(Z \setminus \mathfrak{Z})} \to \mc S|_{(Z \setminus \mathfrak{Z})} \rightarrow 0.
\end{equation}
 All of these sheaves are flat families of \(G\)-sheaves on \(\C^{3}\).

 Take an arbitrary closed point \(z\in Z\setminus \mathfrak{Z}\).
 Let \(\Spec{A} \subset \mc{M}_{C}\) be a zero-dimensional subscheme supported at \(z\) with the following property: the scheme theoretic intersection \(\Spec{A} \cap Z\) is defined by a nonzero element \(\varepsilon \in A\) satisfying \(\mathfrak{m} \, \varepsilon = 0\) for \(\mathfrak{m}\) the maximal ideal of A.
 If we restrict (\ref{eqn:obstruction}) to \(z\),  we obtain a class \(e\) in \(\GExt_{\owe_{\C^3}}^{1}(\mc{S}_{z}, \mc{Q}(-(Z \setminus \mathfrak{Z})_{z})\).  
 As is well known (see \cite[Proposition~7.14]{Buchweitz:smd}), \(e\) is the obstruction class for the Quot-scheme, i.e.,  the class \(e\) is the obstruction to lifting the quotient \(\taut{C} \otimes A/(\varepsilon) \twoheadrightarrow \mc{Q} \otimes A/(\varepsilon)\) to a quotient of \(\taut{C} \otimes A\), flat over \(A\).
 The definition of \(Z\) as a scheme implies that we do not have such a lift
and hence that $e \ne 0$.
 Then we can see that the nontrivial extension \(\mc{K}_{z}\) is \(\theta'\)-stable.
 It follows that \(\mc{K}|_{(Z \setminus \mathfrak{Z})}\) and hence \(\mc{K}|_{(\mc M_{C}\setminus \mathfrak{Z})}\) are families of \(\theta'\)-stable \(G\)-constellations.
 Therefore \(\mc{K}|_{(\mc{M}_{C} \setminus \mathfrak{Z})}\) is isomorphic to \(\taut{C'}\otimes L|_{(\mc{M}_{C} \setminus \mathfrak{Z})}\),  for some line bundle \(L\).
 As a result,  the double dual \(\mc{K}^{**}\) is isomorphic to \(\taut{C'}\otimes L\).
 
 Assume \(\mathfrak{Z}\) is nonempty and take a closed point \(\mathfrak{z}\in \mathfrak{Z}\).
If $\mathfrak{z}$ is not on an associated component of $Z$ of codimension-one then the dimension of the support of \(\mc{Q}\) at \(\mathfrak{z}\) is at most one.  Therefore \((\mc{K}^{**})_{\mathfrak{z}} \cong (\taut{C})_\mathfrak{z}\),  but \(\mc K^{**}\) is also isomorphic to \(\taut{C'}\otimes L\) and hence \((\mc{K}^{**})_{\mathfrak{z}}\) is both \(\theta\)- and \(\theta'\)-stable.  This contradicts the definition of \(Z\), so we may assume \(\mathfrak{z}\) lies on an embedded component of $Z$.
 Let \(D\subseteq Z\) denote the Cartier divisor of \(\mc{M}_{C}\) obtained from \(Z\) by removing the embedded components.  
 Then we can argue as above by considering the surjection \(\taut{C} \to \mc Q\vert_{D}\) with kernel $\mc{K}^{**}$.
 The assumption that \(\mathfrak{z}\) lies on an embedded component of \(Z\) implies that the obstruction $e$ must vanish,
contradicting the $\theta'$-stability of $(\mc K^{**})_{\mathfrak{z}}$.
Thus, \(\mathfrak{Z} = \emptyset\).
\end{proof}

 \begin{corollary}
 \label{coro:type0taut}
 \begin{enumerate}
 \item If \(\rho_{0}\subseteq \Ro\), then \(\mc{S}\cong \mc{Q}'\) and \(\mc{S}'(D) \! := \mc S'\otimes \owe_{\mc M_{C}}(D)\cong \mc{Q}\).  Moreover,  the tautological bundles $\taut{C}$ and $\taut{C'}$ are related by the following elementary transformation of vector bundles (\cite{Maruyama:avb}):
$$
\begin{CD}     @.      0     @.      0      @.  @. \\
@.      @AAA     @AAA      @.  @. \\
0 @>>> \mc S @>>> \taut{C} |_D  @>>> \mc S'(D) @>>> 0 \\
@.       @AAA      @AAA            @|         @.  \\
0 @>>> \taut{C'} @>>> \taut{C} @>>> \mc S'(D) @>>> 0 \\
@.        @AAA     @AAA       @. @. \\
@.        \taut{C}(-D) @= \taut{C}(-D)   @.  @. \\
@.        @AAA     @AAA       @. @. \\
@. 0 @. 0 @. @.
\end{CD}
$$
 \item Otherwise,  \(\rho_{0}\subseteq \Rt\),  in which case \(\mc S(D)\cong \mc{Q}'\) and \(\mc{S'}\cong \mc{Q}\),  and we replace \(\mc{S'}(D)\) by \(\mc{S'}\) and \(\taut{C'}\) by \(\taut{C'}(-D)\) in the above diagram.
 \end{enumerate}
 \end{corollary}
 \begin{proof}
 The proof of the previous proposition showed \(\mc{K}\cong \taut{C'}\otimes L\) and hence,  by the normalisation \(\taut{\rho_{0}}'\cong \owe_{\mc{M}_{C}}\),  that \(\mc{K}_{\rho_{0}}\cong L\).  By definition,  \(\mc{K}_{\rho_{0}}\) is the kernel of \(\taut{\rho_{0}}\to \mc{Q}_{\rho_{0}}\),  so \(L\cong \owe_{\mc{M}_{C}}\) for \(\rho_{0}\subseteq \Ro\).  This gives \(\mc{K} \cong \taut{C'}\) if \(\rho_{0} \subseteq \Ro\).  The proof that \(\mc{K} \cong \taut{C'}(-D)\) if \(\rho_{0} \subseteq \Rt\) is similar.
 \end{proof}

 \begin{corollary}
 \label{coro:taut0changes}
 The tautological line bundles \(\taut{\rho}\) and \(\taut{\rho}'\) are related as follows.
 \begin{enumerate}
 \item \(\rho_{0} \subseteq \Ro\implies \taut{\rho}' \cong \taut{\rho}\) for \(\rho \subseteq \Ro\) and \(\taut{\rho}' \cong \taut{\rho}(-D)\) for \(\rho \subseteq \Rt\).
 \item \(\rho_{0} \subseteq \Rt\implies\taut{\rho}' \cong \taut{\rho}(D)\) for \(\rho \subseteq \Ro\) and \(\taut{\rho}' \cong \taut{\rho}\) for \(\rho \subseteq \Rt\).
 \end{enumerate}
 \end{corollary}

 \subsection{The unstable locus for a wall of type \0} 

 Let \(D\) denote the divisor that forms the unstable locus of a type \0\ wall crossing. 
 \begin{proposition}
 \label{prop:compact}
 The unstable locus \(D\) of a type \0\ wall is compact.
 \end{proposition}
 \begin{proof}
 If \(D\) is noncompact then it contracts onto a curve \(\ell \subset X\).
 Take a nonzero point \(x \in \pi^{-1}(\ell)\subset \C^{3}\).
 The stabiliser \(H\) of \(x\) is nontrivial.
 The category of finite length \(G\)-sheaves supported by the \(G\)-orbit \(G\cdot x\) and the category of finite length \(H\)-sheaves supported at \(x\) are equivalent as in \cite[Lemma~8.1]{Bridgeland:mim}.
 Therefore \(G\)-constellations supported by \(G\cdot x\) correspond to \(H\)-constellations supported at \(x\).  
 Moreover,  the map \(\Theta_{G} \to \Theta_{H}\) of the stability-parameter spaces that is induced by the natural map \(\ind^{G}_{H}\colon R(H) \to R(G)\) is compatible with the above correspondence between \(G\)- and \(H\)-constellations.
 Thus, we can replace \(G\) by \(H\) and thereby assume \(G\) fixes a nonzero \(x \in \C^{3}\).
 Then \(G\) is a subgroup of \(\SL(2, \C) \times \text{Id} \subset \SL(3, \C)\).

 Let \(E\) be a simple \(G\)-constellation.
 Then \(\Supp(E)\) is of the form \(G\cdot x \times a\) for some \(x \in \C^{2}\) and \(a \in \C\).
 If \(t\) is a parameter of \(\C\) at \(a\) then multiplication by \(t\colon E \to E\)
is a nilpotent endomorphism of \(E\) as a \(G\)-sheaf on \(\C^{3}\).
 Since \(E\) is simple, this map must be zero.
 Thus \(E\) is essentially a \(G\)-constellation on \(\C^{2}\) and hence \(\mc{M}_{C}\) is the product of \(\C\) with the moduli space for the restriction of \(G\) to the subgroup \(\SL(2, \C)\) as above.
 In the two-dimensional case,  Kronheimer~\cite{Kronheimer:hkq} proved that the moduli space corresponding to a nongeneric parameter is singular,  so walls of type \0\ do not exist.
 This contradicts the assumption,  so \(D\) is compact.
 \end{proof}

 To proceed with our study of the type \0\ walls of a chamber \(C\) we require the technical results listed in \S\ref{sec:repsquiver}.

 \begin{lemma}
 \label{lemma:reduced}
 The divisor \(D\) is reduced.
 \end{lemma}
 \begin{proof}
 It suffices to show that \(D\) is generically reduced since $D$ is Cartier.
 In a neighbourhood of a point on a two-dimensional torus orbit of $D$, the divisor $D$ is defined by the equation $u_i^{\rho}=0$ for pairs \((\rho, i)\) such that \(\rho \subset \Ro\) and \(\rho\rho_i \not\subset \Ro\) (see \S\ref{sec:repsquiver} for the notation).
 Then Corollary~\ref{coro:oneis0} implies \(D\) is reduced.
 \end{proof}


 \begin{corollary}
 \label{coro:rigid}
 Suppose \(F\) is a \(\theta\)-stable \(G\)-constellation with a nontrivial proper subsheaf \(S\) such that \(\theta_{0}(S) = 0\).  Then either \(S\) or \(Q=F/S\) is rigid.
 \end{corollary}
 \begin{proof}
 Let \(D_{i}\) be an irreducible component of \(D\).  The morphism \(\cont{W}\) is injective,  so \(\dim \GExt_{\mc{O}_{\C^3}}^{1}(Q, S) = 1\) by Lemma~\ref{lemma:fibre}.  Proposition~\ref{prop:rigidity} asserts that either \(S\) or \(Q\) is rigid on the open torus of the toric stratum \(D_{i}\).  Since all such subsheaves (or quotients) are \(\theta_{0}\)-stable,  there is a morphism from \(D_{i}\) to the moduli of \(\theta_{0}\)-stable \(G\)-sheaves that is constant on the subtorus of \(D_{i}\).  It follows that all of the subsheaves (or quotients) on \(D_{i}\) are rigid.  Since \(D_{i}\) was arbitrary,  the result holds on \(D\).  
 \end{proof}

 \begin{remark}
 The corollary shows that either \(S\) or \(Q\) is a constant family on a connected component of \(D\).  Proposition~\ref{prop:rigiddivisor} asserts that \(D\) is connected.   Thus,  the unstable locus \(D\) of a type \0\ wall crossing parametrises either rigid subsheaves or rigid quotient sheaves.
 \end{remark}


 \section{The chamber structure via Fourier--Mukai}
 \label{sec:fourier-mukai}
 In this section the walls of the parameter space \(\Theta\) are written in terms of the \(K\)-theory of the moduli space \(Y = \mc{M}_{C}\) using the Fourier--Mukai transform introduced in \S\ref{sec:bkr}.  For any chamber \(C\) and any \(\theta\in C\),  we write the Fourier--Mukai transform as \(\FM{C} := \FM{\theta}\).

 \subsection{Understanding \protect\(\Theta\protect\) and \protect\(L_{C}(\theta)\protect\) via Fourier--Mukai}
 \label{sec:ThetaFM}
 The results of this section do not require that \(G\subset \SL(3,\C)\) is Abelian. Consider the topological filtration
 \[
 K(Y)_{\Q} = F^{0} \supset F^{1} \supset F^{2} \supset F^{3} = 0,
 \]
 where \(F^{i}\subset K(Y)_{\Q}\) denotes the subspace spanned by sheaves with support of codimension at least \(i\).  Consider also the dimension filtration
 \[
 0 = F_{-1} \subset F_{0} \subset F_{1} \subset F_{2} = K_{0}(Y) _{\Q}
 \]
 of the compactly supported \(K\)-theory of \(Y\), where \(F_{i}\subset K_{0}(Y)_{\Q}\) denotes the subspace spanned by sheaves with support of dimension at most \(i\).
 The pairing
 \begin{equation}
 \label{eqn:pairing}
 \chi(\alpha, \beta) = \sum_i (-1)^i \dim \Ext^i(\alpha, \beta)
 \end{equation}
 is defined for \(\alpha \in K(Y)\) and \(\beta \in K_{0}(Y)\).
 \begin{proposition}
 \label{prop:pairing}
 The pairing \(\chi\colon K(Y) \times K_0(Y) \to \Z\) is perfect.
 Moreover, we have \(F^1 = F_0^{\perp}\) and \(F^2 = F_1^{\perp}\) with respect to this pairing.
 \end{proposition}
 \begin{proof}
 There is a perfect pairing \(\chi^{G} \colon K^{G}(\C^3) \times K^G_0(\C^3) \to \Z\)
 as in \cite[\S9.2]{Bridgeland:mim}.
 The first statement follows from this via the Fourier--Mukai transform.
 Riemann--Roch gives
 \[
\chi(\alpha, \beta) = \int_{\tau^{-1}(0)} \ch(\alpha^{\vee}) \cap \ch_{\tau^{-1}(0)}^Y(\beta)
\cap \operatorname{Td}(Y)
 \]
 where \(\ch \colon K(Y)_{\Q} \to A(Y)_{\Q}\) and \(\ch_{\tau^{-1}(0)}^Y\colon K_0(Y)_{\Q} \to A_{\tau^{-1}(0)}(Y)_{\Q}\) denote the Chern character isomorphisms  (see Fulton~\cite[\S18.3]{Fulton:it}).
 It follows that the intersection pairing \(A(Y)_{\Q} \times A_{\tau^{-1}(0)}(Y)_{\Q} \to \Q\) is perfect.
 This implies the second statement.
 \end{proof}
 
 Write \(\fm{C}\colon K_{0}(Y)\to K_{0}^{G}(\C^{3})\) for the isomorphism (\ref{eqn:fmktheory}) determined by any parameter \(\theta\in C\) and write
 \[
 \fmstar{C}\colon K^{G}(\C^{3})\rightarrow K(Y)
 \]
 for its adjoint.  
 Recall that by definition \(\Theta = R^{\perp}\subset \Hom(R(G),\Q)\).
 If we identify \(R(G)\) with \(K^G_0(\C^{3})\) via the isomorphism sending  \(\rho\) to \(\rho\otimes\owe_{0}\),  then \(\Hom(R(G),\Z) \cong K^G(\C^{3})\) by the perfect pairing \(\chi^G\).
 In this paper, we fix this isomorphism and embed $\Theta$ in $K^G(\C^3)_{\Q}$.
 In other words, we regard \(\theta \in \Theta\) as the class \([\sum_{\rho} \theta(\rho)(\rho \otimes \owe_{\C^3})] \in K^G(\C^3)_{\Q}\).
 
 \begin{corollary}
 For each chamber \(C\subset \Theta\),  the Fourier--Mukai transform induces an isomorphism \(\fmstar{C}\colon \Theta \to F^{1}\).
 \end{corollary}
 \begin{proof}
 The Fourier--Mukai transform maps the class \([R \otimes \owe_0] \in K^G_0(\C^3)\) to the class of a point in \(K_0(Y)\).  
 This spans \(F_0\), so \(\fmstar{C}(\Theta) = F_0^{\perp}= F^{1}\).
 \end{proof}

 For a divisor \(D\) and a compact curve \(\ell\) on \(Y\),  the pairing from Proposition~\ref{prop:pairing} is simply \(\chi(\owe_D, \owe_{\ell}) = -\deg\left(\owe_{Y}(D)\vert_{\ell}\right)\).
 Thus we fix an isomorphism
 \begin{equation}
 \label{eqn:picisom}
 F^1/F^2 \overset{\sim}{\longrightarrow} \Pic(Y)_{\Q}
 \end{equation}
 such that the class \([\mc F]\) of a torsion sheaf \(\mc{F}\) is mapped to \(\det(\mc F)^{-1}\).  Recall that to \(\theta \in \Theta\) we associate a (fractional) line bundle \({L}_{C}(\theta)\) on \(Y\).  
 For the isomorphism \(\fmstar{C}\colon \Theta\to F^{1}\), we have
 \[
 \fmstar{C}(\theta)=\fmstar{C}\big{(}\textstyle{\sum_{\rho}}\theta(\rho)[\rho \otimes \owe_{\C^3}]\big{)} = \sum_{\rho} \theta(\rho) [\taut{\rho}^{-1}],
 \]
 whose determinant is \(L_{C}(\theta)^{-1}\).  The identification (\ref{eqn:picisom}) then gives
 \[
{L}_{C}(\theta) = (\fmstar{C}(\theta)\mod F^{2}) \in F^{1}/F^{2} \cong \Pic(Y)_{\Q}.
 \]
 Thus we have established the following commutative diagram:
 \[
 \begin{CD}
 \Theta @>{\fmstar{C}}>> F^1 \\
 @V{L_{C}}VV              @VV{\res}V  \\
 \Pic(Y)_{\Q}  @<{\sim}<{\det^{-1}}< F^1/F^2
 \end{CD}
 \]
 In particular,  for a curve \(\ell\) on \(Y\) we have  
 \begin{equation}
 \label{eqn:degLC}
 \deg \left(L_{C}(\theta)\vert_{\ell}\right) = \sum_{\rho} \chi(\taut{\rho}\vert_{\ell})\cdot \theta(\rho) = \theta(\fm{C}(\owe_{\ell})).
 \end{equation}

 \subsection{The walls of \protect\(\Theta\protect\) via the Fourier--Mukai transform}
 \label{sec:0viaFM}
 We now use the Fourier--Mukai transform to write the inequalities defining the walls of a chamber \(C\) in terms of the \(K\)-theory of \(\mc{M}_{C}\).  The inequalities determined by walls of type \I\ or \III\ follow immediately from the relation (\ref{eqn:degLC}):

 \begin{corollary}
 \label{coro:curvewalls}
 Given a wall \(W\subset C\) of type \I\ or \III,  let \(\ell\subset \mc{M}_{C}\) denote a curve contracted by \(\cont{W}\).  Then \(\theta(\fm{C}(\owe_{\ell})) = 0\) for all \(\theta\in W\),  and furthermore \(\theta(\fm{C}(\owe_{\ell})) > 0 \text{ for all }\theta\in C\)
 \end{corollary}
 \begin{proof}
 For \(\theta\in W\),  \( L_{C}(\theta)\) is the pull back of a line bundle via \(\cont{W}\) that contracts a curve \(\ell\subset \mc{M}_{C}\) to a point.  That is,  \(\theta\in W\) if and only if \(\deg \left(L_{C}(\theta)\vert_{\ell}\right) = 0\) which,  by (\ref{eqn:degLC}),  is equivalently \(\theta(\fm{C}(\owe_{\ell})) = 0\).  
 \end{proof}

 The analogous statement for walls of type \0\ (see Corollary~\ref{coro:rigidquot} below) is more involved and will occupy the rest of this subsection.  First we establish a useful vanishing result.

 \begin{lemma}
 \label{lemma:tau*}
 For \(\sigma, \rho \in \Irr(G)\) we have \(\RR \tau_*(\taut{\rho}\otimes \taut{\sigma}^{-1}) \cong \left(\pi_*({\owe}_{\C^{3}}) \otimes \sigma \otimes \rho^*\right)^G\). 
 In particular,  \(H^{i}(Y,\taut{\rho} \otimes \taut{\sigma}^{-1}) = 0\) for \(i > 0\).
 \end{lemma}
 \begin{proof}
 It is well known that
 \[
 R^{i}\tau_{*}(\taut{\rho} \otimes \taut{\sigma}^{-1}) \cong H^{i}(Y,\taut{\rho} \otimes \taut{\sigma}^{-1}) \cong \Ext^{i}(\owe_{Y},\taut{\rho} \otimes \taut{\sigma}^{-1}) \cong \Ext^{i}(\taut{\rho}^{-1},\taut{\sigma}^{-1}),
 \]
 which in turn is equal to \(\Hom^{i}_{D(Y)}(\taut{\rho}^{-1}, \taut{\sigma}^{-1})\).  The Fourier--Mukai transform \(\FM{C}\) induces the isomorphism 
 \[
 \Hom^{i}_{D(Y)}(\taut{\rho}^{-1}, \taut{\sigma}^{-1}) \cong \Hom^{i}_{D^G(\C^{3})}(\rho \otimes \owe_{\C^{3}}, \sigma \otimes \owe_{\C^{3}}),
 \]
 which equals \(\GExt^{i}(\rho \otimes \owe_{\C^{3}}, \sigma \otimes \owe_{\C^{3}})\cong \GExt^{i}(\owe_{\C^{3}}, \rho^{*}\otimes \sigma \otimes \owe_{\C^{3}})\).  Now
 \[
 \GExt^{i}(\owe_{\C^{3}}, \rho^{*}\otimes \sigma \otimes \owe_{\C^{3}}) \cong H^{i}(\C^{3},\rho^{*}\otimes \sigma \otimes \owe_{\C^{3}})^{G} = R^{i}\pi_{*}(\rho^{*}\otimes \sigma \otimes \owe_{\C^{3}})^{G}.
 \]
 This is zero unless \(i= 0\).  This completes the proof.
 \end{proof}

 \begin{proposition}
 \label{prop:rigiddivisor}
 Let \(D \subset \mc{M}_C\) be the unstable locus of a wall $W$ of type \0. 
 Write $R=R_1 \oplus R_2$,  with $R_1$ corresponding to the subsheaves
 determining $W$ as in \S\ref{sec:fibreofcont}.
 Then the following hold.
 \begin{enumerate}
 \item[\one] For \(\rho\subseteq R_2\) and \(\sigma\subseteq R_1\),  we have \(H^{i}(\taut{\rho} \otimes \taut{\sigma}^{-1}\vert_D) = 0\) for \(i=0,1,2\).
 \item[\two] \(D\) is connected.
 \item[\three] \(H^{1}(\owe_D) = H^{2}(\owe_D)=0\).
 \end{enumerate}
 \end{proposition}
 \begin{proof}
 By Corollary~\ref{coro:taut0changes}, \(\taut{\rho} \otimes \taut{\sigma}^{-1}\) becomes \(\taut{\rho} \otimes \taut{\sigma}^{-1}(-D)\) after crossing the wall.  Lemma~\ref{lemma:tau*} implies that
 \[
 \RR \tau_* \taut{\rho} \otimes \taut{\sigma}^{-1}(-D)
 = (\pi_*(\owe_M)\otimes \sigma \otimes \rho^*)^{G}
 = \RR \tau_* \taut{\rho} \otimes \taut{\sigma}^{-1}.
 \]
 The first assertion then follows from the long exact sequence in cohomology associated to the sequence 
 \begin{equation}
 \label{eqn:sestaut}
 0 \to \taut{\rho}\otimes \taut{\sigma}^{-1}(-D) \to \taut{\rho}\otimes \taut{\sigma}^{-1} \to \taut{\rho}\otimes \taut{\sigma}^{-1}\vert_{D} \to 0.
 \end{equation}
 
 For the second assertion,  take a connected component $D_1$ of $D$.
 By Corollary~\ref{coro:rigid}, there is either a rigid quotient $Q$ or a rigid subsheaf $S$ destabilising the \(G\)-constellations parametrised by $D_1$.
 We first assume the quotient is rigid.
 Fix \(\rho\subset R_2\) and consider \(\FM{C}^{i}(\taut{\rho}^{-1}\otimes \omega_{D_1})\).  Applying Serre duality to (\ref{eqn:FMi}) gives 
 \begin{equation}
 \label{eqn:FMC}
 \FM{C}^{i}(\taut{\rho}^{-1}\otimes \omega_{D_1}) \cong \bigoplus_{\sigma \in \Irr(G)} H^{2-i}(\taut{\rho}\otimes \taut{\sigma}^{-1}\vert_{D_1})^{*} \otimes \sigma. 
 \end{equation}
 The first assertion implies \(H^{i}(\taut{\rho} \otimes \taut{\sigma}^{-1}\vert_{D_1})\) may only be nonzero for \(\sigma\subset \Rt\).  Then \(\rho,\sigma\subset R_2\),  hence \(\taut{\rho}\vert_{D_1} \cong \taut{\sigma}\vert_{D_1}\) by the rigidity of $Q$,  giving
  \begin{equation}
 \label{eqn:FMCQ}
 \FM{C}^{i}(\taut{\rho}^{-1}\otimes \omega_{D_1}) \cong H^{2-i}(\owe_{D_1})^{*} \otimes Q. 
 \end{equation}
 In the case where the subsheaf is rigid, by using similar arguments or Lemma~\ref{lemma:symmetry}, we obtain
  \[
 \FM{C}^{i}(\taut{\sigma}^{-1}|_{D_1}) \cong H^{i}(\owe_{D_1})\otimes S.
 \]
 Now assume \(D\) is not connected and let \(D_{1}\) and \(D_{2}\) be distinct connected components.  
 First consider the case where both components parametrise rigid quotients, and write \(Q_{1}\) and \(Q_{2}\) for the rigid quotients on \(D_{1}\) and \(D_{2}\) respectively.  
 Since they give the same wall, their images \([Q_1]\) and \([Q_2]\) in the Grothendieck group \(K_{0}^{G}(\C^{3})\cong R(G)\) are the same.  
 It follows from (\ref{eqn:FMCQ}) that \(\fm{C}(\taut{\rho}^{-1}\otimes \omega_{D_{j}}) = [\FM{C}(\taut{\rho}^{-1}\otimes \omega_{D_{j}})] = \chi(\owe_{D_{j}})[Q_{j}]\in R(G)\) for \(j = 1,2\),  so that \(\fm{C}(\taut{\rho}^{-1}\otimes \omega_{D_1})\) and \(\fm{C}(\taut{\rho}^{-1}\otimes \omega_{D_2})\) are proportional in \(R(G)\).  
 This cannot happen since \([\taut{\rho}^{-1} \otimes \omega_{D_1}]\) and \([\taut{\rho}^{-1}\otimes \omega_{D_2}]\) are independent in \(K_{0}(Y)\) and \(\fm{C}\) is an isomorphism.
 The case where both parametrise rigid subsheaves also cannot occur,  the argument is similar.
 Finally,  assume we have a rigid quotient \(Q_1\) on \(D_1\) and a rigid subsheaf \(S_2\) on \(D_2\).
 In this case, we can argue in \(K_0(Y)/F_0 \overset{\sim}{\to} R(G)/\Z[R]\).
 The classes \([\taut{\rho} \otimes \omega_{D_{1}}]\) and \([\taut{\sigma}\vert_{D_{2}}]\) are still independent in the left hand side but \([Q_1] + [S_2] = 0\) in the right hand side, a contradiction.
 Thus \(D\) is connected and we obtain
 \[
 \fm{C}(\taut{\rho}^{-1}\otimes \omega_{D}) = \chi(\owe_{D})[Q]
 \]  
 in the rigid quotient case and
 \[
 \fm{C}(\taut{\sigma}^{-1}|_D) = \chi(\owe_D)[S]
 \]
 in the rigid subsheaf case.
 We show the final assertion in the rigid quotient case, the other case is similar
 (or follows from Lemma~\ref{lemma:symmetry}).
 We have
 \(H^{2}(\owe_{D}) = 0\) since \(X\) is a rational singularity and hence \(\FM{C}^{0}(\taut{\rho}^{-1}\otimes \omega_{D}) = 0\).  Thus we have the distinguished triangle 
 \[
 \FM{C}^{1}(\taut{\rho}^{-1}\otimes \omega_{D})[-1] \to \FM{C}(\taut{\rho}^{-1}\otimes \omega_{D}) \to \FM{C}^{2}(\taut{\rho}^{-1}\otimes \omega_{D})[-2]\stackrel{g}{\longrightarrow}
 \]
 in \(D_{0}^{G}(\C^{3})\),  where \(g \colon \FM{C}^{2}(\taut{\rho}^{-1}\otimes \omega_{D})[-2]\to \FM{C}^{1}(\taut{\rho}^{-1}\otimes \omega_{D})\).  By (\ref{eqn:FMCQ}) and the connectedness of \(D\),  this triangle is 
 \[
 H^{1}(\owe_{D})^{*} \otimes Q[-1] \to \FM{C}(\taut{\rho}^{-1}\otimes \omega_{D}) \to Q[-2] \stackrel{g}{\longrightarrow} H^{1}(\owe_{D})^{*} \otimes Q.
 \]
 The morphism \(g\) is given by an element in \(\GExt_{\owe_{\C^{3}}}^{2}(Q, Q) \otimes H^{1}(\owe_{D})^{*}\).  Serre duality and the rigidity of \(Q\) imply \(\GExt_{\owe_{\C^{3}}}^{2}(Q, Q) = \GExt_{\owe_{\C^{3}}}^{1}(Q, Q)^{\vee}=0\),  and hence \(g\) must be zero.  Thus \(\FM{C}(\taut{\rho}^{-1}\otimes \omega_D)\) is the direct sum of \(Q[-2]\) and \(H^{1}(\owe_{D})^{*} \otimes Q[-1]\).  Since \(D\) is connected and \(\Phi_{C}\) is an equivalence we have \(H^{1}(\owe_{D}) = H^{1}(\owe_{D})^{*} = 0\).
 \end{proof}

 \begin{corollary}
 \label{coro:rigidquot}
 The unstable locus \(D\) of a type \0\ wall parametrises either a fixed rigid quotient \(Q\) or a fixed rigid subsheaf \(S\).  
 Moreover:
 \begin{enumerate}
 \item
 In the rigid quotient case,  we have \(\FM{C}(\taut{\rho}^{-1} \otimes \omega_D) = Q[-2]\) for \(\rho\subseteq R_2\).  In particular, the inequality \(\theta([Q]) < 0\) for \(\theta \in C\) defining the wall can be written \(\theta(\fm{C}(\taut{\rho}^{-1} \otimes \omega_{D})) < 0\).
 \item
 In the rigid subsheaf case,  we have \(\FM{C}(\taut{\sigma}^{-1}\vert_D) = S\) for \(\sigma\subseteq R_1\).  In particular, the inequality \(\theta([S]) > 0\) for \(\theta \in C\) defining the wall can be written \(\theta(\fm{C}(\taut{\sigma}^{-1}\vert_{D})) > 0\).
\end{enumerate}
 \end{corollary}

 \subsection{The chamber structure via Fourier--Mukai}
 Corollaries~\ref{coro:curvewalls} and \ref{coro:rigidquot} enable us to write the inequalities defining a chamber \(C\) via the Fourier--Mukai transform as follows: \(\theta\in C\) if and only if
 \begin{itemize}
 \item \(\theta(\fm{C}(\owe_{\ell})) > 0\) for every exceptional curve \(\ell\).
 \item \(\theta(\fm{C}(\taut{\sigma}^{-1}\vert_{D})) > 0\),  where \(D\) is the unstable locus of a type \0\ wall that parametrises a rigid subsheaf $S$ and $\sigma$ is some (in fact any) representation in \(H^{0}(S)\).
 \item \(\theta(\fm{C}(\taut{\rho}^{-1}\otimes \omega_{D'})) < 0\),  where \(D'\) is the unstable locus of a type \0\ wall parametrising a rigid quotient \(Q\) and \(\rho\) is some (in fact any) representation in \(H^{0}(Q)\).
 \end{itemize}

 Unfortunately it is not clear how to find pairs \((D, \sigma)\) and \((D',\rho)\) before knowing the inequalities from type \0\ walls and their unstable loci.
 We now introduce additional inequalities that are satisfied by parameters \(\theta\) in the chamber \(C\) leading to the more streamlined statement of Theorem~\ref{thm:chamberstrong}.  Of course,  the extra inequalities are redundant.   

 \begin{lemma}
 \label{lemma:quotinequal}
 For every \(\theta \in C\),  every compact reduced divisor \(D\) and every irreducible representation \(\rho\), we have \(\theta(\fm{C}(\taut{\rho}^{-1} \otimes \omega_D)) < 0\).
 \end{lemma}
 \begin{proof}
 Since \(C\) is an open set, a linear function does not achieve a maximum value on \(C\), so it is sufficient to show that \(\theta (\fm{C}(\taut{\rho}^{-1} \otimes \omega_D)) \leq 0\).  By the symmetry of \(\Theta\) from Lemma~\ref{lemma:symmetry}\one\ we may assume \(\rho = \rho_{0}\), in which case (\ref{eqn:FMC}) gives \(\FM{C}^{i}(\omega_D) = \oplus_{\sigma} H^{2-i}(\taut{\sigma}^{-1}\vert_{D})^{*}\otimes \sigma\).  Lemma~\ref{lemma:tau*} gives \(H^{i}(\taut{\sigma}^{-1}) = 0\) for \(i > 0\) which,  when combined with the long exact sequence of cohomologies arising from (\ref{eqn:sestaut}),  proves \(H^{i}(\taut{\sigma}^{-1}\vert_D) \cong H^{i+1}(\taut{\sigma}^{-1}(-D))\) for \(i = 1,2\).  The vanishing of \(H^{3}\)'s gives \(\FM{C}^0(\omega_{D}) = 0\),  so
 \[
 \theta(\fm{C}(\omega_{D})) = \theta([\FM{C}^{2}(\omega_{D})]) - \theta([\FM{C}^{1}(\omega_{D})]).
 \]
 It is enough therefore to show \(\theta([\FM{C}^{2}(\omega_{D})]) \leq 0\) and \(\theta([\FM{C}^{1}(\omega_D)]) \geq 0\).  We may assume \(D\) is connected.  To prove that \(\theta([\FM{C}^{2}(\omega_{D})]) \leq 0\),  note that \(\FM{C}^{2}(\omega_{D}) = \oplus_{\sigma} H^{0}(\taut{\sigma}^{-1}\vert_{D})^{*}\otimes \sigma\).  Taking sufficiently many points \(y_{1}, \dots, y_{k} \in D\) shows that \(\FM{C}^{2}(\omega_{D})\) is a quotient of 
 \[
 \bigoplus_{\sigma} \bigoplus_{i} H^{0}(\taut{\sigma}^{-1} \otimes \owe_{y_i})^{*}\otimes \sigma \cong \bigoplus_{i} \FM{C}(\owe_{y_i}).
 \]
  For \(i = 1,\dots ,k\),   \(\FM{C}(\owe_{y_i})\) is \(\theta\)-stable and so \(\oplus_{i} \FM{C}(\owe_{y_i})\) is \(\theta\)-semistable.  Thus the quotient \(\FM{C}^{2}(\omega_{D})\) must satisfy \(\theta([\FM{C}^{2}(\omega_{D})]) \leq 0\) as required.

 To prove that \(\theta([\FM{C}^{1}(\omega_{D})]) \geq 0\) we first claim that \(\FM{C}^{1}(\omega_{E}) = 0\), where \(E\) is the maximal reduced compact divisor in \(\mc{M}_{C}\).  This is equivalent by (\ref{eqn:FMC}) to \(H^{1}(\taut{\sigma}^{-1}\vert_{E}) = 0\) for all \(\sigma\in \Irr(G)\).  Since \(H^{i}(\taut{\sigma}^{-1}) = 0\) for \(i > 0\),  we prove as above that \(H^{1}(\taut{\sigma}^{-1}\vert_E) \cong H^2(\taut{\sigma}^{-1}(-E))\).  This vanishes because the line bundle \(\owe_{Y}(-E)\) is generated by global sections in codimension-one (\(E\) coincides with the fibre of \(\pi(0) \in \C^{3}/G\) off a one-dimensional subset),  \(H^{2}(\taut{\sigma}^{-1}) = 0\) and we have no \(H^{3}\)'s.  This proves the claim.

 Let \(D\) be an arbitrary reduced compact divisor.
 Write \(E=D+D'\),  where \(E\) is as above and \(\ell'\!:= D \cap D'\) is a (complete intersection) curve.
 By applying \(\FM{C}\) to
 \[
 0 \to \omega_D \oplus \omega_{D'} \to \omega_E \to \omega_{\ell'} \to 0,
 \]
 we obtain \(\FM{C}^{0}(\omega_{\ell'}) \cong \FM{C}^{1}(\omega_D) \oplus \FM{C}^1(\omega_{D'})\). Choosing sufficiently many points \(y_{1}, \dots, y_{r}\) of \(\ell'\) enables us to embed \(\FM{C}^{0}(\omega_{\ell'}) = \oplus_{\sigma} H^{0}(\taut{\sigma} \otimes \omega_{\ell'})\) into
 \[
 \bigoplus_{\sigma} \bigoplus_{i} H^{0}(\taut{\sigma} \otimes \omega_{\ell'}\otimes \owe_{y_i})^{*}\otimes \sigma \cong \bigoplus_{i} \FM{C}(\owe_{y_i}).
 \]
 Here each \(\FM{C}(\owe_{y_i})\) is \(\theta\)-stable so \(\oplus_i \FM{C}(\owe_{y_i})\) is \(\theta\)-semistable.  Thus its submodule \(\FM{C}^{1}(\omega_{D})\) must satisfy \(\theta([\FM{C}^{1}(\omega_{D})]) \geq 0\) as required.
 \end{proof}

 \begin{corollary}
 \label{coro:quotinequal}
 For every \(\theta \in C\),  every compact reduced divisor \(D\) and every irreducible representation \(\rho\), we have \(\theta(\fm{C}(\taut{\rho}^{-1}\vert_D)) > 0\).
 \end{corollary}
 \begin{proof}
 This follows immediately from Lemma~\ref{lemma:symmetry}.
 \end{proof}

 \begin{theorem}
 \label{thm:chamberstrong}
 Let \(C \subset \Theta\) be a chamber.  Then \(\theta\in C\) if and only if
 \begin{itemize}
 \item for every exceptional curve \(\ell\) we have \(\theta(\fm{C}(\owe_{\ell})) > 0\).
 \item for every compact reduced divisor \(D\) and irreducible representation \(\rho\) we have 
 \[
 \theta(\fm{C}(\taut{\rho}^{-1}\otimes \omega_{D})) < 0\quad\text{and}\quad\theta(\fm{C}(\taut{\rho}^{-1}\vert_D)) > 0.
 \]
 \end{itemize}
 \end{theorem}


 \section{Crossing walls of type \I\ and \III}
 This section repeats the analysis of \S\ref{sec:type0}  for walls of type \I\ and \III.  That is,  for walls of each type we determine the unstable locus and calculate the changes in the moduli spaces and the tautological bundles that occur as the parameter \(\theta\in \Theta\) passes through the wall.  

\subsection{Crossing a wall of type \I\ induces a flop}

 \begin{proposition}
 \label{prop:typeIflop}
 Let \(C\) and \(C'\) be chambers separated by a wall \(W\) of type \I.  Then the unstable locus \(Z\subset \mc{M}_{C}\) of the wall crossing is the curve \(\ell\) contracted by \(\cont{W}\).  Furthermore, \(\mc{M}_{C}\) is a flop of \(\mc{M}_{C'}\).
 \end{proposition}
 \begin{proof} We break the proof up into several steps.

 \para \textsc{Step 1:} We begin by showing that \(\ell\) is a connected component of \(Z\).  The wall \(W\) gives rise to a defining inequality \(\theta(\Ro) > 0\) of the chamber \(C\),  for some subrepresentation \(\Ro\subset R\) as in \S\ref{sec:fibreofcont}.  On the other hand,  Corollary~\ref{coro:curvewalls} asserts that the same defining inequality may be written as \(\theta(\fm{C}(\owe_{\ell})) > 0\).  Hence there exists \(q\in\Q_{>0}\) such that \(\theta(\fm{C}(\owe_{\ell})) = q\, \theta(\Ro)\) for all \(\theta\in \Theta\).  Since \(\Theta = R^{\perp}\subset R(G)\cong K_{0}^{G}(\C^{3})\),  there exists \(k\in \Q\) such that \(\fm{C}(\owe_{\ell}) = q\, \Ro - k\, R\).  In this equality, \(\fm{C}(\owe_{\ell})\) is integral,  every coefficient of \(R\) is 1 and \(\Ro\) has both 0 and 1 as its coefficients.  It follows that \(k,q\in \Z\).  Since \(\fm{C}([\text{pt}]) = R\), we have \(\fm{C}(\owe_{\ell}(k)) = q\Ro\).  The line bundle \(\owe_{\ell}(k)\) is not divisible in \(K_0(Y)\),  so \(q=1\), i.e., 
 \[
 \fm{C}(\owe_{\ell}(k)) = \Ro.
 \]
 Then \(\deg(\taut{\rho} \otimes \owe_{\ell}(k))\) is $0$ or $-1$, therefore
 \begin{equation}
 \label{eqn:subtypeI}
S:= \FM{C}(\owe_{\ell}(k))
 \end{equation} 
 is a $G$-sheaf.
 Rigidity of \(\owe_{\ell}(k)\) implies \(S\) is rigid.
 In particular,  \([S]\) is an isolated point in the moduli space of \(\theta_0\)-stable \(G\)-equivariant coherent sheaves on \(\C^3\).  For \(y\in \mc{M}_{C}\),  write \(F_y\) for the corresponding \(G\)-constellation.  The equivalence \(\FM{C}\) gives \(\Hom_{\owe_{\mc M_C}}(\owe_{\ell}(k), \owe_y) \cong \GHom_{\owe_{\C^3}}(S, F_y)\), so
 \[
 y \in \mc M_C \text{ is on } \ell \iff F_y \text{ contains }S\text{ as a subsheaf.}
 \]
 Consider the \(\theta_{0}\)-unstable locus \(Z \subset \mc{M}_{C}\).
 We have a morphism from \(Z\) to the moduli of \(\theta_0\)-stable \(G\)-equivariant coherent sheaves associating \(y\in Z\) to the destabilising subsheaf of \(F_y\).
 The above shows that the fibre over \([S]\) is exactly \(\ell\).
 Since \([S]\) is an isolated point, the curve \(\ell\) is a connected component of \(Z\).  
 Note also that the quotient of \(F_y\) by \(S\) for each \(y \in \ell\) is the rigid \(G\)-sheaf
\[
 Q := \FM{C}^1(\owe_{\ell}(k-1)) = \FM{C}(\owe_{\ell}(k-1))[1].
 \]

 \para \textsc{Step 2:} 
 By symmetry,  the curve \(\ell'\subset \mc{M}_{C'}\) contracted to a point in \(\overline{\mc{M}_{\theta_{0}}}\) is a connected component of the unstable locus \(Z'\) in \(\mc{M}_{C'}\).  
 Set \(Z_{c} := Z\setminus \ell\) and \(Z_{c}' := Z\setminus \ell'\).  
 The isomorphism \(\mc{M}_{C}\setminus Z\cong \mc{M}_{C'}\setminus Z'\) of Lemma~\ref{lemma:common} ensures that the ample line bundle \(L_{C'}(\theta')\) is the proper transform of \(L_C(\theta')\) in the neighbourhood \(\mc{M}_{C}\setminus Z_{c}\) of \(\ell\).  
 The first assertion of Lemma~\ref{lemma:GIT} states that \(L_C(\theta')^{-1}\) is \(f\)-ample.    Since \(f\) contracts only the curve \(\ell\),  we see that \(\mc{M}_{C'}\setminus Z_{c}'\) is an \(L_C(\theta')\)-flop of \(\mc{M}_{C}\setminus Z_{c}\).  
 The question of whether \(\mc{M}_{C'}\) is a flop of \(\mc{M}_{C}\) is local over \(\overline{\mc{M}_{\theta_0}}\),  so the second assertion is proven.

 \para \textsc{Step 3:} 
 Since \(\mc{M}_{C}\dashrightarrow \mc{M}_{C'}\) is a flop of \(\ell\),  \(Z_{c}\subset \mc{M}_{C}\) is set-theoretically identical to \(Z_{c}'\subset \mc{M}_{C'}\).  
 Assume \(Z_{c}\) is nonempty.
 The proofs of Propositions~\ref{prop:Cartier},~\ref{prop:compact} and Lemma~\ref{lemma:reduced} show that \(Z_{c}\) is a compact and reduced Cartier divisor.  
 Furthermore, take a connected component \(D_1\) of \(Z_c\).  Then the argument that proves Corollary~\ref{coro:rigid} shows that either we have a rigid subsheaf \(S'\) or a rigid quotient \(Q'\) that is common to the \(G\)-constellations parametrised by \(D_1\).

 \para \textsc{Step 4:}
 Let \(\mc{M}\) denote the common blow up \(\mc{M}_C\) along \(\ell\) and \(\mc{M}_{C'}\) along \(\ell'\).
 The elementary transformation of the tautological bundle \(\taut{} = \taut{C}\) along the pullback of \(Z\) to \(\mc{M}\) gives a family of \(\theta'\)-stable \(G\)-constellations for some \(\theta'\in C'\),  and symmetrically for \(\taut{}' = \taut{C'}\).  
 Using this relation between the pullbacks of the tautological bundles to \(\mc{M}\) in conjunction with Lemma~\ref{lemma:tau*} establishes the analogue of Proposition~\ref{prop:rigiddivisor}\one\ (with \(Z_{c}\) replacing \(D\)).

 \para \textsc{Step 5:}
 To derive a contradiction we argue as in the proof of Proposition~\ref{prop:rigiddivisor}\two.  Apply the Fourier--Mukai transform so that in \(K_0^G(\C^3)\) we have either
 \[
 \fm{C}(\taut{\sigma}^{-1}\vert_{D_{1}}) = [S'] = [S]\quad\text{or}\quad\fm{C}(\taut{\rho} \otimes \omega_{D_{1}}) = [Q'] = [Q].
 \]
 But this is impossible because \(\fm{C}(\owe_{\ell}(k))= [S]\),  \(\fm{C}(\owe_{\ell}(k-1))= -[Q]\) and \(\fm{C}\) is an isomorphism.
 Thus \(Z_{c}\) is empty and \(Z = \ell\) as required.
 \end{proof}

 \begin{warning}
 The previous proposition is a particular case of a result of Thaddeus~\cite[Theorem~3.3]{Thaddeus:git}.  Note however that the assumption on the codimension of the unstable locus (denoted there \(X^{\pm}\)) is missing from the statement of \cite[Theorem~3.3]{Thaddeus:git};  it appears only in the statement of the simplest case \cite[Proposition~1.6]{Thaddeus:git}.  
 \end{warning}

 \begin{corollary}
 \label{coro:propertransform}
 The tautological bundle \({\taut{\rho}}'\) is the proper transform of \(\taut{\rho}\) and we have:
 \begin{description}
 \item[Case \plusone] \(\rho_0\subseteq \Rt \iff \{\, \deg(\taut{\rho}|_{\ell}) \,|\, \rho \in \Irr(G)\,\}=\{\,0, 1\,\}\).
 \item[Case \minusone] \(\rho_0\subseteq \Ro \iff \{\, \deg(\taut{\rho}|_{\ell}) \,|\, \rho \in \Irr(G)\,\}=\{\,0, -1\,\}\).
 \end{description}
 \end{corollary}
 \begin{proof}
 In the proposition we saw that \(\deg(\taut{\rho} \otimes \owe_{\ell}(k))\) is $0$ or $-1$.
 Taking \(\rho = \rho_0\) shows that \(k=-1\) or 0.
 If \(k=-1\),  Case \plusone\ occurs.
 Otherwise, Case \minusone\ occurs.
 \end{proof}

\subsection{Walls of type \III\ and change of tautological bundles}
 Suppose \(C\) and \(C'\) are chambers separated by a wall \(W\) of type \III\ so that both \(\cont{W}\colon \mc{M}_{C}\to Y_{0}\) and \(\cont{W}'\colon \mc{M}_{C'}\to Y_{0}\) contract a divisor to a curve.  The universal property of blowing up
together with the primitivity of the contractions
ensures that both  \(\mc{M}_{C}\) and \(\mc{M}_{C'}\) can be characterised as the  blow-up the curve that forms the singular
locus of $Y_0$.  In particular,   \(\mc{M}_{C}\) is isomorphic to \(\mc{M}_{C'}\).

 \begin{proposition} 
 \label{prop:typeIIIiso}
 Let \(C\) and \(C'\) be chambers separated by a wall \(W\) of type \III.  Then \(\mc{M}_{C}\) is isomorphic to \(\mc{M}_{C'}\) and the unstable locus \(Z\subset \mc{M}_{C}\) of the wall crossing is the divisor \(D\) that is contracted onto a curve by \(\cont{W}\). 
 \end{proposition}
 \begin{proof}
 We need only prove the second assertion.  Write \(\cont{W}\colon \mc M_{C} \to Y_0\) for the contraction determined by the wall, \(B \subset Y_0\) for the image of the divisor \(D\) and \(\ell\) for a fibre of \(D \to B\).
 We argue as in Proposition~\ref{prop:typeIflop}.
 Thus, there is an integer $k$ such that \(S:= \FM{C}(\owe_{\ell}(k))\) is a \(G\)-sheaf satisfying \([S] = \Ro\), and \(B\) is a connected component of the moduli space of \(\theta_0\)-stable \(G\)-sheaves \(E\) with \(H^0(E) \cong R_1\).
 The divisor \(D\) is the pull-back of \(B\) and is therefore a connected component of \(Z\).
 The proof that \(Z_{c} := Z\setminus D\) is empty is similar to \textsc{Steps 3-5} of Proposition~\ref{prop:typeIflop},  though in \textsc{Step 4} we need not perform a blow-up since we already have \(\mc{M}_{C}\cong \mc{M}_{C'}\).
 \end{proof}

 \begin{corollary}
 \label{coro:changeIII}
 The tautological bundles \(\taut{\rho}\) and \(\taut{\rho}'\) are related as follows:
 \begin{description}
 \item[Case \minusone] \(\rho_0\subseteq \Ro\iff \{\, \deg(\taut{\rho}|_{\ell}) \,|\, \rho \in \Irr(G)\,\}=\{\,0, -1\,\}\),  in which case
 \[
{\taut{\rho}}'=
 \begin{cases}
 \taut{\rho}(-D) & \text{if $\deg(\taut{\rho}|_{\ell})=-1$}\\
 \taut{\rho} & \text{if $\deg(\taut{\rho}|_{\ell})=0$.}
 \end{cases}
 \]
 \item[Case \plusone] \(\rho_0\subseteq \Rt \iff \{\, \deg(\taut{\rho}|_{\ell}) \,|\, \rho \in \Irr(G)\,\}=\{\,0, 1\,\}\), in which case
 \[
 {\taut{\rho}}'=
 \begin{cases}
 \taut{\rho}& \text{if $\deg(\taut{\rho}|_{\ell})=0$}\\
 \taut{\rho}(D) & \text{if $\deg(\taut{\rho}|_{\ell})=1$.}
 \end{cases}
 \]
 \end{description}
 \end{corollary}
 \begin{proof}
 This is proved in the same way as in Corollary~\ref{coro:type0taut}.
 \end{proof}

 \section{Derived equivalences induced by wall crossing}
 \label{sec:twists}
 This section considers the composite \(\FM{C'}^{-1}\circ \FM{C}\colon D(\mc{M}_{C})\to D(\mc{M}_{C'})\) of the Fourier--Mukai transforms determined by adjacent chambers \(C\) and \(C'\).  For chambers separated by a wall of type \0\ or \III\ we determine the explicit self-equivalence of \(D(\mc{M}_{C})\).  For chambers separated by a type \I\ wall we exhibit the equivalence of derived categories between flops.

\subsection{Walls of type \0\ and Seidel--Thomas twists}

 Let \(Y\) be an \(n\)-dimensional nonsingular quasiprojective variety and let \(D_{c}(Y)\) denote the derived category of compactly supported coherent sheaves on \(Y\).

 \begin{definition}\cite{Seidel:bga}
 An object \(\mc{E} \in D_{c}(Y)\) is \emph{spherical} if \(\Hom^{k}_{D(Y)}(\mc{E},\mc{E}) = 0\) unless \(k\) is \(0\) or \(n\), \(\Hom^{0}_{D(Y)}(\mc{E},\mc{E}) \cong \Hom^{n}_{D(Y)}(\mc{E},\mc{E}) \cong \C\) and \(\mc{E} \otimes \omega_{Y} \cong \mc{E}\).
 The \emph{twist} \(T_{\mc{E}}\) along a spherical object \(\mc{E}\) is defined via the distinguished triangle
 \[
 \RR\Hom_{\owe_{Y}}(\mc{E}, \alpha) \ltensor_{\C} \mc{E} \overset{\ev}{\longrightarrow} \alpha \longrightarrow T_{\mc E}(\alpha)
 \]
 for any \(\alpha \in D(Y)\),  where ev is the evaluation morphism.  The twist \(T_{\mc{E}}\) is an exact self-equivalence with inverse \(T'_{\mc{E}}\) defined via the distinguished triangle
 \[
 T'_{\mc E}(\alpha) \longrightarrow \alpha  \overset{\text{ev}}{\longrightarrow} \RR\Hom_{\owe_{Y}}(\alpha,\mc{E})^{\vee} \ltensor_{\C} \mc{E}\quad \mbox{ for any } \alpha \in D(Y).
 \]
 \end{definition}

 Before establishing that a type \0\ wall crossing induces a twist on the derived category,  we pause to review the twist on the Grothendieck groups \(K(Y)\) and \(K_{c}(Y)\) of the corresponding derived categories that we require in \S\ref{sec:maintheorem}.  
 The pairing defined in (\ref{eqn:pairing}) determines \(\chi(\alpha, \beta)\) for \(\alpha, \beta \in K_{c}(Y)\).
 If \(\dim Y\) is odd and \(\omega_{Y}\) is trivial then \(\chi(-,\!-\!)\) is skew symmetric on \(K_{c}(Y)\).
 The following is immediate from the definition.

 \begin{lemma}
 \label{lemma:Ktheorytwist}
 The twist \(T_{\mc{E}}\) acts on \(K(Y)\) and \(K_{c}(Y)\) as
 \[
 (T_{\mc{E}})_{*}(y) = y - \chi\left([\mc{E}], y\right) [\mc{E}].
 \]
 This clearly acts as the identity on \([\mc{E}]^{\perp}\).  Moreover, if the dimension of \(Y\) is odd then
\(\chi((T_{\mc{E}})_{*}(y), [\mc{E}]) = \chi(y, [\mc{E}])\), so the sets
 \[
 [\mc{E}]^{\pm}\!:=\{\, y \in K(Y)\st \pm \chi(y, [\mc E]) > 0 \,\}
 \]
 are invariant under \(T_{\mc{E}}\).
 \end{lemma}

 We now consider the self-equivalence induced by crossing a type \0\ wall.  
 Let \(C, C'\) be adjacent chambers separated by a wall \(W\) of type \0\ and \(\taut{\rho}, \taut{\rho}'\) be the tautological bundles on \(\mc{M}_{C}, \mc{M}_{C'}\) respectively.  Note that \(Y=\mc{M}_{C}\cong\mc{M}_{C'}\).  Hereafter,  the notation \(f(?)\) denotes the functor \([?\mapsto f(?)]\).

 \begin{proposition}
 \label{prop:twist}
 Suppose \(W\) is defined by the equation \(\theta(\fm{C}(\taut{\rho}^{-1}\otimes \omega_{D})) = 0\) (so \(D\) parametrises a rigid quotient \(Q\) with \(\rho\in Q\)).  Then either
 \begin{enumerate}
 \item \label{twist:1} \(\taut{\rho}\vert_{D}\not\cong \owe_{D}\),  in which case \(\FM{C'}^{-1} \circ \FM{C}\cong T_{\taut{\rho}^{-1} \otimes \omega_D}\)\emph{;} or
 \item \(\taut{\rho}\vert_{D}\cong \owe_{D}\),  in which case \(\FM{C'}^{-1} \circ \FM{C}\cong \owe_{Y}(-D)\otimes T_{\omega_D}(?)\),
 \end{enumerate}
 where \(\cong\) between functors denotes natural isomorphism.  Otherwise,  \(W\) is defined by the equation \(\theta(\fm{C}(\taut{\sigma}^{-1}\vert_{D})) = 0\) (so \(D\) parametrises a rigid subsheaf \(S\) and \(\sigma\in S\)).  Then either
 \begin{enumerate}
 \item[3.] \(\taut{\sigma}\vert_{D}\not\cong\owe_{D}\),  in which case \(\FM{C'}^{-1} \circ \FM{C}\cong T'_{\taut{\sigma}^{-1}\vert_{D}}\)\emph{;} or
 \item[4.] \(\taut{\sigma}\vert_{D}\cong\owe_{D}\),  in which case \(\FM{C'}^{-1} \circ \FM{C}\cong T'_{\omega_D}(\owe_{Y}(D)\otimes \,?)\).
 \end{enumerate}
 \end{proposition}

 \begin{proof}
 For part~(\ref{twist:1}),  set \(\mc{E} = \taut{\rho}^{-1}\otimes \omega_{D}\). 
 Since an object in \(D(\mc{M}_C)\) is represented by a complex of sheaves that are direct sums of \(\taut{\sigma}^{-1} \cong \FM{C}^{-1}(\sigma \otimes \owe_{\C^3})\) for \(\sigma \in \Irr(G)\),  we need only establish functorial isomorphims
 \[
 \FM{C'}^{-1} \circ \FM{C}(\taut{\sigma}^{-1}) \cong T_{\mc{E}}(\taut{\sigma}^{-1}).
 \]
 The left hand side is \({\taut{\sigma}'}^{-1}\) which,  by Corollary~\ref{coro:taut0changes},  is isomorphic to \(\taut{\sigma}^{-1}\) for  \(\sigma\not\subset Q\) and to \(\taut{\sigma}^{-1}(D)\) for \(\sigma\subset Q\).  To compute the right hand side, write  \begin{equation}
 \label{eqn:RRtwist1}
 \RR\Hom_{\owe_{Y}}(\mc{E},\taut{\sigma}^{-1}) \cong \RR\Gamma(\taut{\rho}\otimes\taut{\sigma}^{-1}\otimes\owe_{D}[-1]). 
 \end{equation}
 For \(\sigma\not\subset Q\),  this is zero by Proposition~\ref{prop:rigiddivisor}\one,  and hence \(T_{\mc{E}}(\taut{\sigma}^{-1})\cong\taut{\sigma}^{-1}\) as required.  Otherwise,  \(\sigma\subset Q\).  Then \(\taut{\sigma}\vert_{D}\cong \taut{\rho}\vert_{D}\) and hence (\ref{eqn:RRtwist1}) is \(\RR\Gamma(\owe_{D}[-1]) \cong \C[-1]\) by Proposition~\ref{prop:rigiddivisor}\three.  The twist \(T_{\mc{E}}(\taut{\sigma}^{-1})\) is therefore the cone over the map
 \[
 \C[-1]\otimes \mc{E}\longrightarrow \taut{\sigma}^{-1}.
 \]
 This is a sheaf given by a class in \(\Ext^{1}_{\owe_{Y}}(\mc{E},\taut{\sigma}^{-1})\).  Since \(\taut{\sigma}\vert_{D}\cong \taut{\rho}\vert_{D}\), this is simply \(\Ext^{1}_{\owe_{Y}}(\taut{\sigma}^{-1}\otimes \omega_{D},\taut{\sigma}^{-1})\cong \Ext^{1}_{\owe_{Y}}(\omega_{D},\owe_{Y}) \cong H^{2}(\omega_{D})^{\vee}\cong H^{0}(\owe_{D})\cong \C\).  Clearly \(\taut{\sigma}^{-1}(D)\) is one such (nontrivial) class,  so \(T_{\mc{E}}(\taut{\sigma}^{-1}) \cong \taut{\sigma}^{-1}(D)\).
 
 For the functoriality,   consider the diagram
 \[
 \begin{CD}
D(\mc M_{C}) @>f>> D(\mc M_{C'}) \\
@VVV                            @VVV   \\
D\big{(}\mc{M}_{C} \setminus \tau^{-1}(0)\big{)} @>\sim>> D\big{(}\mc M_{C'} \setminus (\tau')^{-1}(0)\big{)},
 \end{CD}
 \]
 where \(\tau,\tau'\) denote the contractions from the moduli \(\mc{M}_{C},\mc{M}_{C'}\) to \(\C^{3}/G\).  The vertical arrows are restrictions and the lower horizontal arrow is the equivalence induced by \(\mc{M}_{C} \smallsetminus \tau^{-1}(0) \cong \mc{M}_{C'} \smallsetminus (\tau')^{-1}(0)\).  The restrictions of \(\FM{C}^{-1},\FM{C'}^{-1}\) to \(D^{G}(\C^{3}\setminus \{0\})\) are Fourier--Mukai transforms determined by \(\mc{R}\vert_{\mc{M}_{C} \smallsetminus \tau^{-1}(0)}\cong\mc{R}'\vert_{\mc{M}_{C'} \smallsetminus (\tau')^{-1}(0)}\),  so the diagram commutes if \(f = \FM{C'}^{-1} \circ \FM{C}\).  On the other hand,  \(\mc{E}\) is supported on \(D\subset \tau^{-1}(0)\) so the restriction of the twist \(T_{\mc{E}}\) to \(D\big{(}\mc M_{C} \setminus \tau^{-1}(0))\) is the identity.  Thus the diagram also commutes if \(f = T_{\mc{E}}\).  Moreover, the restriction maps
 \[
 \Hom(\mc{F}_{1},\mc{F}_{2}) \to \Hom({\mc{F}_{1}}|_{\mc{M}_{C} \setminus \tau^{-1}(0)},{\mc{F}_{2}}|_{\mc M_{C} \setminus \tau^{-1}(0)})
 \]
  are injective for locally free sheaves \(\mc{F}_{1},\mc{F}_{2}\) on \(\mc{M}_{C}\).  This establishes the functoriality and hence completes the proof of part (1).  Part (2) is similar.   Parts (3) and (4) follow from (1) and (2) respectively by exchanging the roles of \(C\) and \(C'\).
 \end{proof}

 \subsection{Walls of type \I\ and Bondal--Orlov derived equivalences}

 Let \(C, C'\) denote chambers separated by a wall of type \I.  Then \(\mc{M}_{C} \dashrightarrow \mc{M}_{C'}\) is the classical flop with respect to a \((-1, -1)\)-curve \(\ell \subset \mc{M}_{C}\).  If we denote by \(p\colon \widetilde{\mc{M}} \to \mc{M}_{C}\) the blowup of \(\mc{M}_{C}\) along \(\ell\),  then the exceptional locus \(E\) is isomorphic to \(\mathbb{P}^{1} \times \mathbb{P}^{1}\) and we have another contraction \(q\colon \widetilde{\mc{M}} \to \mc{M}_{C'}\).  In this situation,  Bondal and Orlov~\cite{Bondal:sod} observed that there is an equivalence of derived categories
 \[
 \RR q_{*} \LL p^{*}\colon D(\mc{M}_{C}) \to D(\mc{M}_{C'})
 \]
with inverse \(\RR p_{*} (\owe(E) \otimes \LL q^{*}(?))\colon D(\mc{M}_{C'}) \to D(\mc{M}_{C})\).

 \begin{proposition}
 \label{prop:bondalorlov}
 When crossing a wall of type \I,  the induced equivalence of derived categories is determined by the cases \plusminusone\ of Corollary~\ref{coro:propertransform} as follows:
 \begin{enumerate}
 \item If \plusone\ holds we have \(\FM{C'}^{-1} \circ \FM{C} \cong \RR q_{*} \LL p^{*}\).
 \item Otherwise \minusone\ holds and \(\FM{C'}^{-1} \circ \FM{C} \cong \RR q_{*} (\owe(E) \otimes \LL p^{*}(?))\).
 \end{enumerate}
 \end{proposition}
  \begin{proof}
 For part (1),  we establish functorial isomorphisms 
 \[
 \FM{C'}^{-1} \circ \FM{C}(\taut{\sigma}^{-1}) \cong \RR q_{*} \LL p^{*}(\taut{\sigma}^{-1}).
 \]
 The left hand side is isomorphic to \({\taut{\rho}'}^{-1}\).  It follows from Case \plusone\ from Corollary~\ref{coro:propertransform} that the right hand side is the proper transform of \(\taut{\sigma}^{-1}\), i.e.,  \({\taut{\rho}'}^{-1}\).  The functoriality is proved as in Proposition~\ref{prop:twist}.  Part (2) is similar.
 \end{proof}

 \subsection{Self-equivalences induced by walls of type \III}
 Finally,  assume that chambers \(C\) and \(C'\) are separated by a wall of type \III.  
 Let \(Y=\mc{M}_{C}\cong \mc{M}_{C'}\) denote the moduli and \(Y \to Y_0\) the contraction determined by the wall.
 In this case, the unstable locus \(D \subset Y\) is the divisor contracted to a curve in \(Y_{0}\).
 Let \(p, q\) be the first and the second projections of the fibre product \(Y \times_{Y_0} Y\).

 \begin{proposition}
 \label{prop:horjaszendroi}
 When crossing a wall of type \III,  the induced self-equivalence of the derived category is determined by the cases \plusminusone\ of Corollary~\ref{coro:changeIII} as follows:
 \begin{enumerate}
 \item If \minusone\ holds we have \(\FM{C'}^{-1} \circ \FM{C}\cong\RR q_{*} \LL p^{*}\).
 \item If \plusone\ holds then \(\FM{C'}^{-1} \circ \FM{C}\cong\RR q_{*} (\omega_{Y\times_{Y_0}Y} \otimes \LL p^{*}(?))\).
 \end{enumerate}
 \end{proposition}
 \begin{proof}
 Similar to that for walls of type \I,  except that the unstable locus \(D\) need not be compact so instead we consider the restrictions of \(\FM{C}^{-1}, \FM{C'}^{-1}\) to \(D^{G}(\C^{3}\setminus \pi^{-1}(\tau(D)))\).
 \end{proof}

 \begin{remark}
 Set \(E=D, Z=\cont{W}(D)\) and \(\mc{E}=\owe_D\). Then we have
 \[
 \owe_{Y\times Y}\cong \widetilde{\mc E}_R \otimes p^* \owe_Y(-D),
 \quad
 \omega_{Y \times Y} \cong \widetilde{\mc E}_L \otimes q^* \owe_Y(D),
 \]
 where \(\widetilde{\mc E}_R\) and \(\widetilde{\mc E}_L\) denote the kernels (or correspondences) of the EZ transforms introduced by Horja~\cite{Horja:dca}.
 Thus, our derived equivalences are EZ transforms composed with tensor products by line bundles \(\owe_Y(\pm D)\).
 \end{remark}



 \section{Proof of the main theorem}
 \label{sec:maintheorem}
 In this section we prove Theorem~\ref{thm:main}.
 Since every projective crepant resolution is obtained by successive flops from \(\ghilb\),  it is enough to show that if a chamber \(C\) defines \(\mathcal{M}_C \cong Y\) then for any flop \(Y'\) of \(Y\) there is a chamber \(C'\) such that \(\mathcal{M}_{C'} \cong Y'\).  It is not true in general that such \(C'\) is found by a single wall crossing from \(C\),  see Example~\ref{ex:3stepflop}.
 Thus we consider the union of the closures \(\overline{C}\) of all chambers \(C\) with \(\mc{M}_{C} \cong Y\) and try to find a desired wall of type \I\ on the boundary of this union.
 To understand the relation with the ample cone clearly, we transform each \(\overline{C}\) via the isomorphism \(\fm{C}^{*}\colon \Theta\to F^{1}\).
 It turns out that adjacent chambers in \(\Theta\) need not be transformed to adjacent cones in \(F^{1}\),  but tensoring by some line bundle \(\owe_{Y}(D)\) fills the gap.

 To make the argument precise we first introduce an action of a subgroup of \(\Pic(Y)\) on \(F^1\) as follows.  
 The Picard group of \(Y\) acts on \(K(Y)\) by tensor product, preserving the filtration \(F^{\bullet}\) and hence \(\Pic(Y)\) acts on \(F^1\).
 Note that this action is trivial on \(F^1/F^2 \cong \Pic(Y)_{\Q}\)
 because $[\mathcal F \otimes L] - [\mathcal F] \in F^2$ for
 a torsion sheaf $\mathcal F$ and a line bundle $L$.
 Let \(S_{1}, \dots, S_{b}\) denote the compact irreducible surfaces on \(Y\) and \(\Pic^c(Y) \subset \Pic(Y)\) the subgroup generated by the \(\owe_Y(S_i)\).
 For \(\xi \in F^1\), \(\owe_Y(S_i)\) acts as
 \[
 \owe_Y(S_i) \otimes \xi = \xi + [\owe_{S_i} \otimes \xi],
 \]
 where the equality holds in \(K(Y)\).  Indeed, \(\owe_Y(S_i) \otimes \xi = [\owe_{Y}]\otimes \xi + [\owe_{S_i}(S_{i})] \otimes \xi\),  and the difference \(([\owe_{S_{i}}(S_{i})] - [\owe_{S_{i}}])\otimes \xi\) lies in \(F^{3} = 0\).

 We define \(\Amp'(Y) \subset F^1/F^2\) to be the set of nef line bundles that do not contract a divisor to a point, i.e.,  \(\Amp'(Y)\) is obtained from the nef cone \(\overline{\Amp(Y)}\) by removing the walls that induce birational morphisms of type \II.  Clearly
 \[
\Amp(Y) \subset \Amp'(Y) \subset \overline{\Amp(Y)}.
 \]

 \begin{lemma}
 \label{lemma:linindep}
 Suppose \(\res(\xi)\in \Amp'(Y)\) for \(\xi\in F^{1}\).  Then \([\owe_{S_1} \otimes \xi], \dots, [\owe_{S_b} \otimes \xi]\) are linearly independent in $F^2$.
 \end{lemma}
 \begin{proof}
 Put \(L = \res(\xi)\) and take a general section \(s\) of \(L^n\) for suitable \(n > 0\).
 The section \(s\) determines a nonsingular surface \(S\).
 Our assumption that \(L \in \Amp'(Y)\) ensures that \(S \cap S_i\) are nonempty curves that are contracted in \(\tau\).
 This implies that the intersection matrix \((c_1(L). S_i. S_j)_{i,j}\) is negative definite and proves the assertion.
 \end{proof}

 \begin{proposition}
 For every \(\xi \in \res^{-1}(\Amp'(Y))\) there is a neigbourhood \(N(\xi)\) of \(\xi\) such that only finitely many pairs \((L, C)\) satisfy \(N(\xi) \cap L \otimes \fmstar{C}(\overline{C}) \neq \emptyset\), for \(L\in \Pic^{c}(Y)\) and \(C\subset \Theta\) a chamber.  Moreover
 \begin{equation}
 \label{eqn:union}
 \res^{-1}(\Amp'(Y)) \subset \bigcup_{L \in \Pic^c(Y)}
 \bigcup_{\mathcal{M}_C \cong Y} L \otimes \fmstar{C}(\overline{C}).
 \end{equation}
 \end{proposition}
 \begin{proof}
 Since $\fmstar{C}(\overline{C})$ is bounded in the fibre direction of the map \(\res\colon F^{1}\to F^{1}/F^{2}\) by Theorem~\ref{thm:chamberstrong},  and since the generators of \(\Pic^{c}(Y)\) act independently in the fibre direction of \(\res\) by Lemma~\ref{lemma:linindep},  we obtain the first assertion of the proposition.
 It follows that the intersection of the right hand side with the left hand side of (\ref{eqn:union}) is closed in the left hand side. 
Therefore, to establish (\ref{eqn:union}) it is sufficent to show that the interior of the left hand side is contained in the right hand side.  That is,  we now show
 \[
\res^{-1}(\Amp(Y)) \subset \bigcup_{L \in \Pic^c(Y)}
\bigcup_{\mathcal{M}_C \cong Y} L \otimes \fmstar{C}(\overline{C}).
 \]
 If this is not true, then again by the first assertion, there is a boundary point $L \otimes \fmstar{C}(\theta_0)$ of the right hand side that is contained in the left hand side, where $L \in \Pic^c(Y)$, $C$ is a chamber with $\mathcal{M}_C \cong Y$ and $\theta_0 \in \overline{C}$.
 We may assume that $\theta_0$ is in the interior of a wall $W$ of $C$.
 Since $L_C(\theta_0)=\res(\fmstar{C}(\theta_0))$ is ample,  $W$ is a wall of type \0\ determined by a rigid subsheaf or quotient sheaf $F$.
 Let $D$ be the compact divisor that forms the unstable locus for this wall and let \(C'\) be the chamber satisfying \(\overline{C'} \cap \overline{C} = W\).
 Then, by Proposition \ref{prop:twist} and Lemma \ref{lemma:Ktheorytwist}, we have either
 \begin{equation}
 \label{eqn:F1conesa}
 \fmstar{C}(\overline{C}) \cap \fmstar{C'}(\overline{C'}) =\fmstar{C}(W)(=\fmstar{C'}(W))
 \end{equation}
 or
 \begin{equation}
 \label{eqn:F1conesb}
 \fmstar{C}(\overline{C}) \cap (\owe_Y(\pm D) \otimes \fmstar{C'}(\overline{C'}))
 =\fmstar{C}(W)(=\owe_Y(\pm D) \otimes \fmstar{C'}(W)).
 \end{equation}
 If (\ref{eqn:F1conesa}) holds,  then \(L \otimes \fmstar{C}(\theta_{0})\in L\otimes \fmstar{C}(W) =  L\otimes\fmstar{C}(\overline{C}) \cap L\otimes\fmstar{C'}(\overline{C'})\). This contradicts the assumption that \(L \otimes \fmstar{C}(\theta_{0})\) is a boundary point of the right hand side of (\ref{eqn:union}).  The contradiction for the (\ref{eqn:F1conesb}) case is similar.
 \end{proof}
 
 \begin{corollary}
 Let $B \subset \overline{\Amp(Y)}$ be a wall of the nef cone giving rise to a birational contraction of type \I.  Then there is a chamber $C$ defining \(\mc{M}_{C}\cong Y\) and there is a type \I\ wall \(W\) of \(\mc{C}\) such that \(\res(\fmstar{C}(W)) \subset B\).  The same statement holds with \III\ replacing \I.
\end{corollary}

The main theorem follows from this and Proposition~\ref{prop:typeIflop}.


 \section{The moduli space of \protect\(G\protect\)-clusters}
 \label{sec:ghilb}
 Throughout this section we let \(C\) (rather than \(C_{0}\) as in \S\ref{sec:intro}) denote the chamber defining \(\ghilb = \mc{M}_{C}\).  In this case,  the inequalities defined by walls of type \0\ can be written in a particularly simple form,  see Theorem~\ref{thm:ghilbwalls}.  We also establish that none of the inequalities \(\theta(\fm{C}(\owe_{\ell})) > 0\) determined by \((-1,-1)\)-curves \(\ell\subset \ghilb\) that induce flops are redundant.  As a consequence,   the flop of any such curve is achieved by crossing a wall of \(C\);  this is false in general,  see Example~\ref{ex:3stepflop}.  The main point is that \(\ghilb\) exhibits nice properties,  even compared with other crepant resolutions \(\mc{M}_{C'}\to\C^{3}/G\).

 \subsection{How to calculate \protect\(\ghilb\protect\)}
 \label{sec:htc}
 Let \(G\subset \SL(3,\C)\) be a finite Abelian subgroup of order \(r = \vert G\vert\).  
 Choose coordinates \(x,y,z\) on \(\C^{3}\) to diagonalise the action of \(G\),  write \(\widetilde{M} \cong \Z^3\) for the lattice of Laurent monomials in \(x,y,z\) and \(\widetilde{N}\) for the dual lattice with basis \(e_{1},e_{2},e_{3}\).  
 Write \(g = \mbox{diag}(\varepsilon^{\alpha_{1}},\varepsilon^{\alpha_{2}},\varepsilon^{\alpha_{3}})\in G\) with \(0\leq \alpha_{j} < r\),  where \(\varepsilon\) is a primitive \(r\)th root of unity.  
 To each \(g\in G\) we associate the vector \(v_{g} = \frac{1}{r}(\alpha_{1},\alpha_{2},\alpha_{3})\).  
 Write \(N := \widetilde{N} + \sum_{g\in G} \; \Z\cdot v_{g}\) and \(M := \Hom(N,\Z)\) for the dual lattice of \(G\)-invariant Laurent monomials.  
 The affine toric variety \(U_{\sigma} = \Spec\C[\sigma^{\vee}\cap M]\) defined by the cone \(\sigma = \sum \R_{\geq 0} e_{i}\) in \(N_{\R} := N\otimes_{\Z}\R\) is the quotient \(\C^3/G\).  
 The \emph{junior simplex} \(\Delta \subset N_{\R}\) is the triangle with vertices \(e_{1}, e_{2}, e_{3}\),  containing lattice points \(\frac{1}{r}(\alpha_{1}, \alpha_{2}, \alpha_{3})\) with \(\frac{1}{r}(\alpha_{1} + \alpha_{2} + \alpha_{3}) = 1\).   
 Crepant toric resolutions \(X_{\Sigma}\rightarrow\C^{3}/G\) are determined by basic triangulations \(\Sigma\) of \(\Delta\);  we identify triangulations of \(\Delta\) with fans,  so vertices,  lines and triangles of \(\Sigma\) define torus-invariant surfaces, curves and points respectively of \(X_{\Sigma}\).  

 Nakamura~\cite{Nakamura:ago} proved that \(\ghilb = X_{\Sigma}\) for one such triangulation \(\Sigma\).  
 Craw and Reid~\cite{Craw:htc} introduced an alternative construction of \(\Sigma\) by first partitioning \(\Delta\) into \emph{regular triangles}.  
 The fan \(\Sigma\) is obtained by drawing \(r-1\) tesselating lines parallel to the sides of each regular triangle of side \(r\).  
 There are two types of regular triangle of side \(r\) (see \cite[Figure~10]{Craw:htc}):
 \begin{enumerate}
 \item[(a)] a \emph{corner triangle}.  By permuting \(x,y,z\) if necessary,  we may assume the sides of this triangle are cut out by the ratios \(x^{d}\! : y^{b}\), \(y^{e}\! : \!x^{a}\) and \(z^{f}\! : \! y^{c}\) for some \(d - a = e - b - c = f = r\);  or 
 \item[(b)] the (unique) \emph{meeting of champions triangle} whose sides are cut out by the ratios \(x^{d}\! : y^{b}\), \(x^{a}\! : \!z^{f}\) and \(y^{e}\! : \! z^{c}\) for some \(d - a = e - b = f - c = r\).
 \end{enumerate}
 Moreover,  the lines of the regular tesselations of the regular triangles are cut out by the \(G\)-invariant ratios of monomials
 \begin{eqnarray}
     x^{d-i}\!:\!y^{b+i}z^{i};\quad y^{e-j}\!:\!z^{j}x^{a+j};\quad z^{f-k}\!:\!x^{k}y^{c+k}\quad \mbox{in Case (a),} \label{eqn:prop3.2a}\\
     x^{d-i}\!:\!y^{b+i}z^{i},\quad y^{e-j}\!:\!z^{c+j}x^{j},\quad z^{f-k}\!:\!x^{a+k}y^{k}\quad \mbox{in Case (b),}\label{eqn:prop3.2b}
 \end{eqnarray}
 for \(i,j,k = 0,\dots ,r-1\).   
 The edges of a basic triangle \(\tau\in \Sigma\) defining the affine chart \(U_{\tau} \subset X_{\Sigma}\) are cut out by the ratios from (\ref{eqn:prop3.2a}) or (\ref{eqn:prop3.2b}) if \(i,j,k\) satisfy either \(i+j+k = r-1\),  in which case \(\tau\) is an \emph{up} triangle,  or \(i+j+k=r+1\),  in which case \(\tau\) is \emph{down}.  
 An up triangle \(\tau\) in Case (a) defines \(U_{\tau}\cong \C^{3} = \Spec \C[\xi,\eta,\zeta]\) with
 \begin{equation}
 \label{eqn:coordsupa}
 \xi = x^{d-i}/y^{b+i}z^{i},\quad \eta = y^{e-j}/z^{j}x^{a+j}, \quad \zeta = z^{f-k}/x^{k}y^{c+k}.
 \end{equation}
 A down triangle \(\tau\) in Case (a) defines \(U_{\tau}\cong\C^{3} = \Spec \C[\lambda,\mu,\nu]\) such that
  \begin{equation}
 \label{eqn:coordsdowna}
 \lambda = y^{b+i}z^{i}/x^{d-i},\quad \mu = z^{j}x^{a+j}/y^{e-j}, \quad \nu = x^{k}y^{c+k}/z^{f-k}.
 \end{equation}
 Case (b) is similar.  
 See \cite{Craw:htc} for more details.

 \subsection{Walls of the chamber defining \protect\(\ghilb\protect\)}
 Let \(C\) denote the chamber defining \(\mc{M}_{C} = \ghilb\).   
 That is,  \(C\) is the unique chamber containing the cone \(\Theta_{+}\) introduced in \S\ref{sec:chamberghilb}.
 We now show that type \0\ walls of \(C\) defined by rigid subsheaves are related to the marking of torus-invariant surfaces in \(\ghilb\) with characters of \(G\) according to the recipe of Reid~\cite{Reid:mc}.

 First we review Reid's recipe.  
 Lines in the fan \(\Sigma\) are cut out by the \(G\)-invariant ratios listed in (\ref{eqn:prop3.2a}) and (\ref{eqn:prop3.2b}).  
 The monomials in each ratio lie in the same character space of the \(G\)-action,  and we mark the line with the common character.  
 As for the (interior) vertices \(v\) in \(\Sigma\),  there are essentially two cases:  
 \begin{itemize}
 \item If three straight lines in \(\Sigma\) pass through \(v\) then \(v\) defines a del Pezzo surface of degree 6.  
 In this case,  \(v\) is marked with a pair of characters \(\rho\), \(\rho'\) determined uniquely from the characters marking the three straight lines meeting at \(v\). 
 \item Otherwise,  \(v\) is marked with a single character \(\rho\) determined uniquely from the characters marking the lines meeting at \(v\).
 \end{itemize}
 In this way,  every character of \(G\) marks either a unique interior vertex in \(\Sigma\) or a line in \(\Sigma\) possibly passing through several vertices.   
 Compact torus-invariant surfaces \(D\) in \(\ghilb\) are the toric strata defined by the interior vertices,  and we say that a character \(\rho\) marks \(D\) when \(\rho\) marks the corresponding vertex \(v\).  
 Complete proofs and worked examples of Reid's recipe are given by Craw~\cite{Craw:emc}. 

 The reader unfamiliar with Reid's recipe may take the following lemma as the definition of the marking of a divisor \(D\) with character(s).

 \begin{lemma}
 \label{lemma:socle}
 Let \(D\) be a compact  torus-invariant divisor \(D\) in \(\ghilb\).  A character \(\rho\) marks \(D\) if and only if \(\rho\) lies in the socle of every torus-invariant \(G\)-cluster in \(D\).
 \end{lemma}
 \begin{proof} 
 The proof of Theorem~6.1 from Craw~\cite{Craw:emc} establishes that if \(\rho\) marks \(D\) then \(\rho\) lies in the socle of every torus-invariant \(G\)-cluster in \(D\).   
 For the converse,  we prove the case where \(D\) is a del Pezzo of degree 6, the other cases are similar (and easier!).  
 The divisor \(D\) is the toric stratum of a vertex \(v\) inside a regular triangle of side \(r \geq 3\),  and the torus-invariant points of \(D\) are the strata defined by the six triangles in \(\Sigma\) with vertex \(v\).  
 If \(v\) lies in a regular corner triangle then the ratios (\ref{eqn:prop3.2a}) cutting out the lines through \(v\) satisfy \(i+j+k=r\) for some \(1\leq i,j,k\leq r-1\).  
 Write \(\owe_{Z}\) for the \(G\)-cluster corresponding to the toric stratum of the triangle \(\tau\) with vertex \(v\) and two edges cut out by \(x^{d-i}\! : \! y^{b+i}z^{i}\) and \(z^{f-k}\! : \! x^{k}y^{c+k}\).  
 The method from \cite[\S5]{Craw:emc} shows that the socle of \(\owe_{Z}\) consists of the monomials \(y^{e-j}z^{i}\) and \(x^{k}y^{e-j}\) whose characters \(\rho\) and \(\rho'\) mark \(D\),  together with the monomials 
 \[
 m_{1} = x^{d-i-1}y^{c+k}, \; m_{2} = x^{d-i-1}z^{j}, \; m_{3} = x^{a+j}z^{f-k-1}, \; m_{4} = y^{b+i}z^{f-k-1}
 \]
 whose characters we denote \(\rho_{1},\dots ,\rho_{4}\) respectively.  
 We claim that \(\rho_{1},\rho_{2}\notin\socle(\owe_{Z'})\),  where \(\owe_{Z'}\) is the \(G\)-cluster defined by the other triangle with vertex \(v\) and one edge cut out by \(x^{d-i}\! : \! y^{b+i}z^{i}\).  
 Indeed,  the calculation of \(\owe_{Z'}\) gives \(x\cdot m_{1},x\cdot m_{2}\in \socle(\owe_{Z'})\).
 Hence \(m_{1},m_{2}\in \owe_{Z'}\smallsetminus \socle(\owe_{Z'})\) so that neither \(\rho_{1}\) nor \(\rho_{2}\) lie in \(\socle(\owe_{Z'})\).  
 Similarly,  if  \(\owe_{Z''}\) is the \(G\)-cluster defined by the triangle with vertex \(v\) that shares the edge cut out by \(z^{f-k}\! : \! x^{k}y^{c+k}\),  then  \(z\cdot m_{3},z\cdot m_{4}\in \socle(\owe_{Z''})\) so \(\rho_{3},\rho_{4}\notin \socle(\owe_{Z''})\) as before.   
 This completes the proof when \(v\) lies in a regular corner triangle,  the meeting of champions case is similar.
  \end{proof}

 \begin{corollary}
 \label{coro:irred}
 A divisor \(D\) in \(\ghilb\) parametrising a rigid subsheaf defining a wall of type \0\ in \(C\) is necessarily irreducible.
 \end{corollary}
 \begin{proof}
 Write \(S\) for the rigid subsheaf.
 If \(\rho\in \socle(S)\) then \(\rho\) lies in the socle of every torus-invariant \(G\)-cluster \(\owe_{Z}\) in \(D\).  
 Assume \(D\) is not irreducible and let \(D_{1}, D_{2}\) be irreducible components.  Then \(\rho\) lies in the socle of every torus-fixed point of both \(D_{1}\) and \(D_{2}\).  
 Lemma~\ref{lemma:socle} implies that \(\rho\) marks both \(D_{1}\) and \(D_{2}\),  a contradiction.
 \end{proof}

 \begin{proposition}
 \label{prop:ghilbsubs}
 \(\theta([S]) = \theta(\fm{C}(\taut{\rho}^{-1}\vert_{D})) = 0\) is a type \0\ wall of \(C\) if and only if \(D\) is irreducible and \(S = \owe_{0}\otimes \rho\) for some \(\rho\) marking \(D\).  In this case,  the wall is \(\theta(\rho) = 0\) and \(\FM{C}^{-1}(\owe_{0}\otimes \rho) = \taut{\rho}^{-1}\vert_{D}\).
 \end{proposition}
 \begin{proof}
 We claim that if \(\rho\) marks a (compact irreducible) divisor \(D\) then \(\theta(\rho) = 0\) is a type \0\ wall of \(C\) with unstable locus \(D\).  Indeed,  \(\rho\subseteq S\) by Lemma~\ref{lemma:socle} so that \(\rho\) is a subsheaf of a \(G\)-cluster and hence \(\theta(\rho)\geq 0\) for all \(\theta\in\overline{C}\).  Thus \(\rho\) lies in the cone dual to \(\overline{C}\),  so \(\overline{C}\cap \rho^{\perp}\) is a face of \(\overline{C}\).  The cone \(\Theta_{+}\) from (\ref{eqn:Theta+}) lies in \(C\) and \(\overline{\Theta_{+}}\cap \rho^{\perp}\) is a codimension-one face of \(\overline{\Theta_{+}}\), so \(W = \{\theta\in \overline{C}\st \theta(\rho)=0\} = \overline{C}\cap \rho^{\perp}\) is a wall of \(C\).  The unstable locus is a divisor containing \(D\) and hence,  by the previous corollary,  coincides with \(D\) as required.

 Suppose \(S = \FM{C}(\taut{\rho}^{-1}\vert_{D})\) is a rigid subsheaf \(S\) giving a wall of type \0.  Then \(D\) is irreducible by the previous corollary.  We may assume that \(\rho\) lies in the socle of \(S\),  so \(\rho\) marks \(D\) by Lemma~\ref{lemma:socle}.  The claim asserts that \(\theta(\rho) = 0\) is a wall of type \0\ with unstable locus \(D\).  Corollary~\ref{coro:rigidquot} gives \(\fm{C}(\taut{\rho}^{-1}\vert_{D}) = \owe_{0}\otimes \rho\) and hence \(S =  \owe_{0}\otimes \rho\) as required.  Conversely,  the claim shows that \(0 = \theta(\rho) = \theta([S])\) is a wall of type \0.

 To prove the final statement,  note that \(\FM{C}(\taut{\rho}^{-1}\vert_{D}) = \owe_{0}\otimes \rho\) for \(S = \owe_{0}\otimes \rho\) by Corollary~\ref{coro:rigidquot}.  The result follows by applying \(\FM{C}^{-1}\).
 \end{proof}

 \begin{remark}
 This result shows that for every compact irreducible divisor \(D\) and for each \(\rho\) marking \(D\),  if \(\taut{\sigma}\vert_{D} \cong \taut{\rho}\vert_{D}\) for some \(\sigma\in \Irr(G)\) then \(\sigma = \rho\). 
 \end{remark}
 
 \begin{theorem}
 \label{thm:ghilbwalls}
 A parameter \(\theta\) lies in the chamber \(C\) defining \(\ghilb\) if and only if
 \begin{itemize}
 \item for every exceptional curve \(\ell\) we have \(\theta(\fm{C}(\owe_{\ell})) > 0\).
 \item for every compact irreducible divisor \(D\) and for every character \(\rho\) marking \(D\) we have \(\theta(\fm{C}(\taut{\rho}^{-1}\vert_{D})) = \theta(\rho) > 0\).
 \item for every compact divisor \(D'\) (parametrising a rigid quotient) we have \(\theta(\fm{C}(\omega_{D'})) < 0\).
 \end{itemize}
 \end{theorem}
 \begin{proof}
 Compare Theorem~\ref{thm:chamberstrong} and Proposition~\ref{prop:ghilbsubs}. For the rigid quotients,  note that \(\rho_{0}\in Q\) because \(G\)-clusters are quotients \(\owe_{\C^{3}}/\mc{I}\).  Clearly \(\taut{\rho_{0}}^{-1}\) is trivial on \(\ghilb\), so the inequality is simply \(\theta(\fm{C}(\omega_{D'})) < 0\). 
  \end{proof}

 \begin{example}
 \label{ex:11.1.2.8}
 For the cyclic quotient singularity \(\frac{1}{11}(1,2,8)\),  write \(\rho_{0},\dots ,\rho_{10}\) for the irreducible representations and write parameters \(\theta\) as \((\theta_{0},\dots ,\theta_{10})\) where \(\theta_{i} = \theta(\rho_{i})\). 
 The toric fan defining \(\ghilb\) together with its marking according to Reid's recipe appears in \cite[Figure~4(b)]{Craw:emc}.  
 The vertices are marked with \(\rho_{3}, \rho_{4},\rho_{7},\rho_{9},\rho_{10}\) and we denote by \(D_{3}, D_{4}, D_{7}, D_{9}, D_{10}\) the corresponding compact irreducible surfaces.  The remaining nontrivial characters mark the lines.  
 
 The three \((-1,-1)\)-curves \(\ell\subset \ghilb\) are marked with \(\rho_{1},\rho_{5},\rho_{6}\) and we label the inequalities \(\theta(\fm{C}(\owe_{\ell})) > 0\) arising from the type \I\ walls accordingly:
 \[
 f_{1} := \theta_1 + \theta_3 + \theta_9 > 0, \quad f_5 := \theta_5 + \theta_7 > 0, \quad f_6 := \theta_6 > 0.
 \]
 There is a unique Hirzebruch surface \(\mathbb{F}_{4}\subset \ghilb\).  The pair of toric curves that are fibres of the morphism \(\mathbb{F}_{4}\to \mathbb{P}^{1}\) are marked with \(\rho_{8}\) so we label the inequality \(\theta(\fm{C}(\owe_{\ell})) > 0\) arising from the type \III\ wall accordingly:
 \[
 f_8 := \theta_8 + \theta_9 + \theta_{10} > 0.
 \]
 It is also convenient to introduce the inequality \(f_2 := \theta_2 + \theta_3 + \theta_7 + \theta_{10} + 2\theta_4 > 0\) that is of the form \(\theta(\fm{C}(\owe_{\ell})) > 0\) for the \((1,-3)\)-curve \(\ell\) in \(\mathbb{P}^{2}\subset \ghilb\) marked with \(\rho_{2}\).  We will soon see that this inequality is redundant (it must be since there are no type \II\ walls!).   For the type \0\ walls,  the inequalities 
 \[
 \theta_{3} > 0, \quad \theta_{4} > 0, \quad \theta_{7} > 0, \quad \theta_{9} > 0, \quad \theta_{10} > 0.
 \]
 arising from the walls whose unstable locus \(D_{i}\) parametrises rigid subsheaves are determined by the marking of surfaces according to Theorem~\ref{thm:ghilbwalls}.  

 Finally,  to calculate the compact divisors \(D'\) parametrising rigid quotients \(Q\),  recall that \(\rho_{i}\subset Q\) if and only if \(\taut{\rho_{i}}\vert_{D'} \cong \owe_{D'}\).  Thus,  we require only the divisors on which tautological bundles are trivial.  We write the appropriate divisor \(D'\) (written as a linear combination of the irreducible divisors \(D_{i}\)) next to the corresponding inequality:
 \begin{eqnarray*}
 D_{4}  & & f_2 > \theta_4  \\
 D_{9}  & & f_1 + f_5 + f_8 > \theta_9  \\
 D_{4}+D_{7} & & f_2 + f_5 + f_6 > \theta_4 + \theta_7 \\
 D_{3} + D_{4} & & f_2 + f_8 + f_6 > \theta_3 + \theta_4 \\
 D_{3} + D_{4} + D_{7} & & f_1 + f_2 + f_5 + f_6 > \theta_3 + \theta_4 + \theta_7 \\
 D_{4} + D_{7}+ D_{10} & & f_2 + f_5 + f_6 + f_8 > \theta_4 + \theta_7 + \theta_{10} \\
 D_{3} + D_{4} + D_{7} + D_{9} + D_{10} & & f_1 + f_2 + f_5 + f_6 + f_8 > \theta_3 + \theta_4 + \theta_7 + \theta_9 + \theta_{10}
 \end{eqnarray*}
 Note that since \(\theta_{4} > 0\),  the inequality determined by \(D' = D_{4}\) ensures that \(f_{2} > 0\) is redundant.  Notice also that the cone \(\Theta_{+} = \{\theta\in \Theta\st \theta_{i} > 0 \text{ for }i > 0\}\) is a proper subcone in the chamber \(C\) cut out by all of these inequalities.  
 \end{example}

 \subsection{Flopping a single curve in \protect\(\ghilb\protect\)}
 
 Fix a \((-1,-1)\)-curve \(\ellzero\) in \(\ghilb\) that induces a flop.  We now show that the inequality \(\theta(\fm{C}(\owe_{\ellzero})) > 0\) determined by \(\ellzero\) according to Theorem~\ref{thm:ghilbwalls} is not redundant.  It follows that the equation \(\theta(\fm{C}(\owe_{\ellzero})) = 0\) defines a type \I\ wall of \(C\) and hence,  by Proposition~\ref{prop:typeIflop},  the flop of \(\ellzero\) is induced by crossing the wall.

 First, note that \(\ellzero\) is the toric stratum of a line in the triangulation \(\Sigma\) lying strictly inside a regular triangle.

 \begin{proposition}
 \label{prop:deg0or1}
 Every tautological bundle \(\taut{\rho}\) has degree \(0\) or \(1\) on \(\ellzero\).
 \end{proposition}
 \begin{proof}
 We prove the case where the line lies inside a regular corner triangle,  the meeting of champions case is similar.  To compute \(\degree \taut{\rho}\vert_{\ellzero}\),  compare the \(G\)-clusters \(\owe_{Z}\) and \(\owe_{Z'}\) corresponding to the toric strata defined by the triangles \(\tau\) and \(\tau'\) that intersect along the line.  One of these triangles is \emph{down},  say \(\tau\),  so \(U_{\tau} \cong \C^{3}\) has coordinates listed in (\ref{eqn:coordsdowna}).  It follows that the subscheme \(Z\subset \C^{3}\) corresponding to \((\lambda,\mu,\nu) = (0,0,0)\) in \(U_{\tau}\) has defining ideal
 \[
 I = \langle x^{d-i+1},y^{e-j+1},z^{f-k+1},y^{b+i}z^{i},x^{a+j}z^{j},x^{k}y^{c+k},xyz\rangle.
 \]
 The set \(\Gamma\) of monomials in \(\C[x,y,z]\smallsetminus I\) forms a basis of \(H^{0}(Z,\owe_{Z})\) so by \S\ref{sec:chamberghilb},  the generator of \(\taut{\rho}\) over \(U_{\tau}\) is the (unique) monomial in \(\Gamma\) lying in the \(\rho\)-character space.  To compute the set \(\Gamma'\) of monomials in the complement of the ideal \(I'\) defining \(Z'\subset \C^{3}\),  we assume \(\tau\cap\tau'\) is cut out by \(x^{d-i}\! : \!y^{b+i}z^{i}\).  Then \(\tau'\) is \emph{up} and \(U_{\tau'} \cong \C^{3}\) has coordinates
 \[
 \xi = x^{d-i}/y^{b+i}z^{i}, \quad\eta = y^{e-j+1}/z^{j-1}x^{a+j-1},\quad \zeta = z^{f-k+1}/x^{k-1}y^{c+k-1},
 \]
 hence \(I' = \langle x^{d-i},y^{e-j+1},z^{f-k+1},y^{b+i+1}z^{i+1},x^{a+j}z^{j},x^{k}y^{c+k},xyz\rangle\).  Let \(\gamma\) denote the set obtained from \(\Gamma\) by removing the monomials divisible by \(x^{d-i}\).  This coincides with the set obtained from \(\Gamma'\) by removing the monomials divisible by \(y^{b+i}z^{i}\).  We have \(\degree \taut{\rho}\vert_{\ellzero} = 0\) whenever \(\gamma\) contains a monomial in the \(\rho\)-character space.  The bijection \(\Gamma\smallsetminus \gamma \rightarrow \Gamma'\smallsetminus \gamma\) defined by \(m\mapsto\lambda\cdot m = (y^{b+i}z^{i}/x^{d-i})\cdot m\) determines the transition function of \(\taut{\rho}\) from \(U_{\tau}\) to \(U_{\tau'}\),  hence \(\degree \taut{\rho}\vert_{\ellzero} = 1\) for \(\rho\in \Gamma\smallsetminus \gamma\). This proves the result when the line \(\tau\cap\tau'\) is cut out by \(x^{d-i}\! : \!y^{b+i}z^{i}\),  the cases where \(\tau\cap\tau'\) is cut out by \(y^{e-j}\!:\!z^{j}x^{a+j}\) and \(z^{f-k}\!:\!x^{k}y^{c+k}\) are similar.
 \end{proof}

 \begin{remark}
 Corollary~\ref{coro:propertransform} asserts that Proposition~\ref{prop:deg0or1} is a necessary condition for \(\ellzero\) to determine a wall of type \I.
 \end{remark}

 \begin{lemma}
 For every compact invariant irreducible surface \(D\) there exists a character \(\sigma\) of \(G\) such that \(\taut{\sigma}\vert_{\ellzero}\) is trivial and \(\taut{\sigma}\vert_{D}\) is nontrivial.
 \end{lemma}
 \begin{proof}
 We prove the case where \(D\) is a del Pezzo surface of degree 6,  the other cases are similar.  Write \(\taut{i}\), \(\taut{j}\), \(\taut{k}\),  \(\taut{l}\),  \(\taut{m}\) for the tautological bundles satisfying \(\taut{l}\otimes \taut{m} = \taut{i}\otimes \taut{j}\otimes \taut{k}\) from Theorem~6.1 of \cite{Craw:emc}.  If \(\taut{i}\), \(\taut{j}\), \(\taut{k}\) are all nontrivial on \(\ellzero\) then they have degree 1 on \(\ellzero\) by the previous proposition.  We deduce from the above relation that either \(\taut{l}\) or \(\taut{m}\) has degree at least 2 on \(\ellzero\).  This contradicts the previous proposition,  hence at least one of \(\taut{i}\), \(\taut{j}\), \(\taut{k}\) is trivial on \(\ellzero\).
 \end{proof}

 \para For the curve \(\ellzero\) as above,  let \(\theta_{0}\in\Theta\) be such that \(\theta_{0}(\rho) = 0\) if \(\taut{\rho}\vert_{\ellzero}\) is nontrivial,  and \(\theta_{0}(\rho) > 0\) if \(\taut{\rho}\vert_{\ellzero}\) is trivial and \(\rho \neq \rho_{0}\).  Then \(\theta_{0}\in  \overline{\Theta_{+}} \subset \overline{C}\) and \(\deg \lb{C}(\theta_{0})\vert_{\ellzero} = 0\).

 \begin{lemma}
 The parameter \(\theta_{0}\) satisfies the following:
 \begin{itemize}
 \item for every exceptional curve \(\ell\) we have \(\theta_{0}(\fm{C}(\owe_{\ell})) \geq 0\).
 \item for every compact irreducible divisor \(D\) and for every character \(\rho\) marking \(D\) we have \(\theta_{0}(\fm{C}(\taut{\rho}^{-1}\vert_{D})) = \theta_{0}(\rho) \geq 0\).
 \item for every compact divisor \(D'\) (parametrising a rigid quotient) we have \(\theta_{0}(\fm{C}(\omega_{D'})) < 0\).
 \end{itemize}
 \end{lemma}
 \begin{proof}
 Since \(\theta_{0} \in \overline{C}\),  it is enough to show \(\theta_{0}(\fm{C}(\omega_{D'}))<0\) for every compact divisor \(D'\) parametrising a rigid quotient \(Q\).  Corollary~\ref{coro:rigidquot} gives \(Q = \FM{C}(\omega_{D'})[2]\),  so it is enough to show \(\theta_{0}([Q]) < 0\) or,  equivalently,  \(\sum_{\rho \notin Q} \theta_{0}(\rho) > 0\).  We already have \(\sum_{\rho \notin Q} \theta_{0}(\rho) \geq 0\) because \(\rho_{0}\in Q\) and \(\theta_{0}\in \overline{\Theta_{+}}\).  By the above lemma,  there exists \(\sigma\) such that \(\taut{\sigma}\vert_{\ellzero}\) is trivial and \(\taut{\sigma}\vert_{D'}\) is nontrivial,  hence \(\sigma\notin Q\).  The definition of \(\theta_{0}\) then gives \(\sum_{\rho \notin Q} \theta_{0}(\rho) > 0\) as required.
 \end{proof}

 \begin{lemma}
 \label{lemma:eta}
 There exists \(\eta \in \overline{C}\) satisfying \(\lb{C}(\eta) \cong \lb{C}(\theta_{0})\) such that
 \begin{itemize}
 \item for every exceptional curve \(\ell\) we have \(\eta(\fm{C}(\owe_{\ell})) \geq 0\).
 \item for every compact irreducible divisor \(D\) and for every character \(\rho\) marking \(D\) we have \(\eta(\fm{C}(\taut{\rho}^{-1}\vert_{D})) = \eta(\rho) > 0\).
 \item for every compact divisor \(D'\) (parametrising a rigid quotient) we have \(\eta(\fm{C}(\omega_{D'})) < 0\).
 \end{itemize}
 \end{lemma}
 \begin{proof}
 Let \(\epsilon > 0\) be sufficiently small.  Since the sheaves \(\owe_{D}\) defined by compact irreducible divisors \(D\) form a basis of \(F_{2}/F_{1} = (F^{2})^{\vee}\),  it is possible to choose \(\eta\) such that \(\lb{C}(\eta) \cong \lb{C}(\theta_0)\) and \(\eta(\fm{C}(\owe_{D})) = \theta_{0}(\fm{C}(\owe_{D})) + \epsilon\) for every compact irreducible surface \(D\).   In other words,  we perturb \(\theta_{0}\) in the fibre direction of the map \(F^{1}\to F^{1}/F^{2}\).  Then we have
 \begin{eqnarray*}
 \eta(\fm{C}(\taut{\rho}^{-1}\vert_{D})) = \theta_{0}(\fm{C}(\taut{\rho}^{-1}\vert_{D})) + \epsilon & > & 0 \\
 \eta(\fm{C}(\omega_{D'})) = \theta_{0}(\fm{C}(\omega_{D'})) + n \cdot \epsilon & < & 0,
 \end{eqnarray*}
 where \(n\) is the number of irreducible components of \(D'\).  The inequalities in the statement of the lemma hold so \(\eta \in \overline{C}\).
 \end{proof}

 The parameter \(\eta \in \overline{C}\) of Lemma~\ref{lemma:eta} satisfies \(\eta(\fm{C}(\owe_{\ellzero})) = \deg \lb{C}(\eta)\vert_{\ellzero} = 0\),  but lies strictly off the walls of type \0.
 We now construct a parameter that also lies off the walls of type \I\ and \III\ defined by curves \(\ell\neq \ellzero\).   
 Indeed,  consider a nef line bundle \(L\) such that the degree of \(L\) is zero on \(\ellzero\) and positive on the other curves.  
 There exists \(\xi \in \Theta\) (which need not be in \(\overline{C}\)) such that \(\lb{C}(\xi) \cong L\).  
 Choose a sufficiently large \(t \in \mathbb{Q}\) to ensure that the parameter \(\theta \! := \xi + t\eta\) lies in \(\overline{C}\) but does not lie on walls defined by rigid subsheaves or quotients.  
 Since \(\lb{C}(\xi)\) has positive degree on all curves \(\ell\neq \ellzero\),  we have 
 \[
 \theta(\fm{C}(\owe_{\ell})) = \degree \lb{C}(\theta)\vert_{\ell} > \degree \lb{C}(\eta)\vert_{\ell} = \eta(\fm{C}(\owe_{\ell}))\geq 0,
 \]
 so \(\theta\) lies strictly off walls defined by curves \(\ell\neq \ellzero\).  Moreover,
 \[
 \theta(\fm{C}(\owe_{\ellzero})) = \degree \lb{C}(\theta)\vert_{\ellzero} = \degree L\vert_{\ellzero} + t\cdot\degree \lb{C}(\eta)\vert_{\ellzero} = t\cdot \eta(\fm{C}(\owe_{\ellzero})) = 0,
 \]
 so \(\theta\) lies on the wall of type \I\ determined by \(\ellzero\).  In particular,  \(\theta\) lies in \(\overline{C}\) and satisfies all of the strict inequalities defining \(C\) except \(\theta(\fm{C}(\owe_{\ellzero})) > 0\).  This proves the main result of this section:

 \begin{theorem}
 \label{thm:singleflop}
 The flop of any single torus-invariant curve \(\ellzero\) in \(\ghilb\) is achieved by crossing a wall of the chamber \(C\).
 \end{theorem}

 It is often possible to realise a pair of flops from \(\ghilb\) by consecutive wall crossings,  though we do not know if this is true for any finite Abelian subgroup of \(\SL(3,\C)\).  However,  the next example illustrates a sequence of three classical flops from \(\ghilb\) that cannot be realised with only three wall crossings.

\begin{example}
 \label{ex:3stepflop}
 Consider the action of \(G = \Z/6\times \Z/2\) with generators the diagonal matrices \(g_{(1,0)} = \diag(\varepsilon,\varepsilon,\varepsilon^{4})\) and \(g_{(0,1)} = \diag(-1,1,-1)\),  where \(\varepsilon\) is a primitive \(6^{\text{th}}\) root of unity.  The quotient \(\C^{3}/G\) is the singularity of type \(\frac{1}{6}(1,1,4)\oplus \frac{1}{2}(1,0,1)\).  The fan defining \(\ghilb\) is shown in Figure~\ref{fig:3stepflop}.
 \begin{figure}[!ht]
 \centering
 \psset{unit=1cm}
 \begin{pspicture}(5.5,5)
 \put(0,0.2){
 \psline{-}(0,0)(5.5,0)\psline{-}(5.5,0)(2.75,4.763)\psline{-}(2.75,4.763)(0,0)
 \psdots(0,0)(5.5,0)(2.75,4.763)
 \rput(3.15,4.763){\(e_{1}\)}
 \rput(5.85,0){\(e_{2}\)}
 \rput(-0.35,0){\(e_{3}\)}
 \psdots(2.75,0)(1.375,2.3815)(4.125,2.3815)(1.375,0.794)(4.125,0.794)(2.75,1.5877)(2.75,3.175)
 \psset{linecolor=gray}
 \psline(0,0)(2.75,1.5877)
 \psline(5.5,0)(2.75,1.5877)
 \psline(2.75,4.763)(2.75,1.5877)
 \psline(2.75,0)(1.375,0.794)
 \psline(2.75,0)(4.125,0.794)
 \psline(4.125,0.794)(1.375,0.794)
 \psline(1.375,2.3815)(2.75,3.175)
 \psline(1.375,2.3815)(1.375,0.794)
 \psline(1.375,0.794)(2.75,3.175)
 \psline(4.125,2.3815)(2.75,3.175)
 \psline(4.125,0.794)(2.75,3.175)
 \psline(4.125,2.3815)(4.125,0.794)
 \psset{linecolor=white}
 \rput(2.75,2.3814){\rput(-0.2,-0.2){\pspolygon[fillcolor=white,fillstyle=solid](0,0)(0.4,0)(0.4,0.4)(0,0.4)\rput(0.2,0.2){\(\ell\)}}}
 \rput(2.0625,2){\rput(-0.2,-0.2){\pspolygon[fillcolor=white,fillstyle=solid](0,0)(0.5,0)(0.5,0.5)(0,0.5)\rput(0.2,0.25){\(\ell_{1}\)}}}
 \rput(3.4375,2){\rput(-0.2,-0.2){\pspolygon[fillcolor=white,fillstyle=solid](0,0)(0.5,0)(0.5,0.5)(0,0.5)\rput(0.2,0.25){\(\ell_{2}\)}}}
 }
 \end{pspicture}
 \caption{The fan $\Sigma$ of $\ghilb$ for $G = \frac{1}{6}(1,1,4)\oplus \frac{1}{2}(1,0,1)$}
 \label{fig:3stepflop}
 \end{figure}
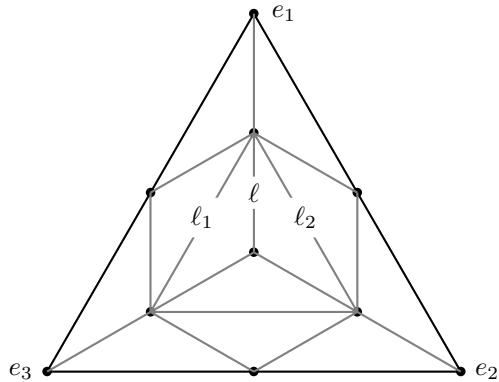 
 
 The toric strata labelled \(\ell\), \(\ell_{1}\), \(\ell_{2}\) in Figure~\ref{fig:3stepflop} determine torus-invariant curves in \(\ghilb\).  Reid's recipe marks the lines from \(e_{i}\) to \(\frac{1}{3}(1,1,1)\) with \(\rho_{(2,0)}\),  so the exceptional \(\mathbb{P}^{2}\) (the toric stratum of the ray through \(\frac{1}{3}(1,1,1)\)) is marked with \(\rho_{(2,0)}\otimes \rho_{(2,0)} = \rho_{(4,0)}\).  Set \(\sigma = \rho_{(4,0)}\).  The line bundle \(\taut{\sigma}\) is generated by \(y^{4}\) at one torus-invariant point of \(\ell\) and by \(z^{4}\) at the other.  Since \(\ell\) is cut out by the ratio \(y^{2}\! : \! z^{2}\), we have \(\deg(\taut{\sigma}\vert_{\ell}) = 2\).  Similarly,  \(\deg(\taut{\sigma}\vert_{\ell_{1}}) = \deg(\taut{\sigma}\vert_{\ell_{2}}) = 0\).  

 Let \(Y'\) be the crepant resolution obtained by flopping \(\ghilb\) along \(\ell_{1}\) then \(\ell_{2}\),  and let \(\ell'\subset Y'\) denote the proper transform of \(\ell\).  The flops of both \(\ell_{1}\) and \(\ell_{2}\) are achieved by crossing a wall of the corresponding chamber and hence the tautological bundles \(\taut{\rho}'\) on the moduli space \(Y'\) are the proper transforms of those on \(\ghilb\).  The generators of \(\taut{\sigma}'\) at the torus-invariant points of \(\ell'\) are \(y^{4}\) and \(z^{4}\) so, as above,  we have \(\deg(\taut{\sigma}'\vert_{\ell'}) = 2\).  It follows from Corollary~\ref{coro:propertransform} that the \((-1,-1)\)-curve \(\ell'\) induces a flop but it does not determine a wall of type \I, i.e., the inequality \(\theta(\fm{C'}(\owe_{\ell'}))\geq 0\) is redundant. 
 \end{example}

 \begin{remark}
 Proposition~\ref{prop:deg0or1} also holds if \(\ell_{0}\) is a torus-invariant \((0,-2)\)-curve.  As a result, the proof of Theorem~\ref{thm:singleflop} establishes that every birational contraction of type \III\ from \(\ghilb\) can be achieved by passing into a wall of the chamber \(C\).
 \end{remark}


 \section{Representations of the McKay quiver}
 \label{sec:repsquiver}
 
 In this section, we prove technical lemmas required in \S\ref{sec:classify-walls}. 
 Following Sardo Infirri~\cite{Sardo-Infirri:ros},  we regard points of \(\mc{M}_{C}\) as representations of the McKay quiver rather than as \(G\)-constellations.
 We are mainly concerned with representations defined by points in two-dimensional torus-orbits.

 \subsection{The universal representation}
 Recall the construction of the moduli space from \S\ref{sec:construction}.
 Assume \(n = 3\) and \(G\subset \SL(V)\) is Abelian.
 Decompose \(V^{*} = \rho_1 \oplus \rho_2 \oplus \rho_3\) as a representation of \(G\).
 An element \(B \in \Hom_{\C[G]}(R, V \otimes R)\) decomposes as 
 \[
 B = \oplus_{\rho\in \Irr(G)}(b_1^{\rho}, b_2^{\rho}, b_3^{\rho}) = (b_i^{\rho})_{i, \rho},
 \]
 where \(b_i^{\rho}\colon R_{\rho} \to R_{\rho\rho_i}\) is a linear map.
 The commutativity condition \(B \wedge B = 0\) becomes \(b_j^{\rho\rho_i} b_i^{\rho}  = b_i^{\rho\rho_j} b_j^{\rho}\) for every \(1\leq i,j\leq 3\) and \(\rho\in \Irr(G)\). These are the defining equations of \(\mc{N} \subset \Hom_{\C[G]}(R, V \otimes R)\).

 We regard \((b_i^{\rho})\) as a representation of the \emph{McKay quiver}.
 This by definition is the quiver with vertex set \(\Irr(G)\) and arrows \(a_i^{\rho}\) from \(\rho\) to \(\rho\rho_i\), one for each \(i=1, 2, 3\) and \(\rho \in \Irr(G)\).
 Since \(\dim \taut{\rho} = 1\),  regard each \(b_i^{\rho}\) as a complex number and define a map that sends the arrow \(a_{i}^{\rho}\) to \(b_i^{\rho}\in \C\).  This is a representation of the McKay quiver with dimension vector \((1,\dots,1)\) (see \cite{Sardo-Infirri:ros}).
 In this way, \(\mc{M}_{C}\) is the moduli space of equivalence classes of \(\theta\)-stable representations of the McKay quiver with dimension vector \((1,\dots,1)\) satisfying the commutativity relations.

 The tautological bundle \(\taut{}\) on \(\mc{M}_{C}\) determines the universal representation of the McKay quiver which we write as \((u_i^{\rho})\),  where \(u_i^{\rho}\colon \taut{\rho} \to \taut{\rho\rho_i}\) a homomorphism of line bundles.
 Let \(x, y, z\) be the coordinates on \(V = \C^3\) lying in \(\rho_1, \rho_2, \rho_3 \subset V^{*}\) respectively.
 Then \(xyz\) is \(G\)-invariant and we regard it as a function on \(\mc{M}_{C}\).

 \begin{lemma}
 For any \(\rho \in \Irr(G)\),  \(u_3^{\rho\rho_1\rho_2}u_2^{\rho\rho_1}u_1^{\rho}\colon \taut{\rho} \to \taut{\rho}\) is multiplication by \(xyz\).
 \end{lemma}
 \begin{proof}
 This follows from the definition of the morphism \(\mc{M}_{\theta} \to X = \C^{3}/G\) given in Proposition~\ref{prop:maptox}.
 \end{proof}

 Let \(T\) denote the big torus of the toric variety \(\mc{M}_{C}\).  Consider an affine open set \(U_{\tau}\) corresponding to a cone \(\tau\).  
 Since \(\mc{M}_{C}\) is a crepant resolution, we have coordinates \(\xi_1, \xi_2, \xi_3\) on \(U_{\tau}\) such that \(\xi_1\xi_2\xi_3 = xyz\).
 Recall that a two dimensional \(T\)-orbit of \(U_{\tau}\) is precisely the locus \((\xi_{i} = 0)\),  for some \(i = 1,2,3\).  This gives the following result:
 
 \begin{corollary}
 \label{coro:oneis0}
 On a two-dimensional $T$-orbit, exactly one of \(u_1^{\rho}, u_2^{\rho\rho_1},  u_3^{\rho\rho_1\rho_2}\) is zero.
 Similarly, exactly one of \(u_1^{\rho}, u_3^{\rho\rho_1},  u_2^{\rho\rho_1\rho_3}\) is zero.
 Moreover, the zero has multiplicity one along the orbit.
 \end{corollary}

 \subsection{The McKay quiver embedded in a torus}
 It is convenient to consider the \emph{universal cover} of the McKay quiver.  (It is not actually a simply connected graph,  but it does enable one to ``unwrap'' the McKay quiver across the whole plane just like the universal cover of the two-torus.)
 We now consider the character group \(G^{*}\) rather than the set \(\Irr(G)\).  
 By taking \(\rho_1, \rho_2, \rho_3\) as generators of \(G^{*}\) we have an exact sequence
 \begin{equation}
 \label{eqn:lattice}
0 \to M \to \Z^3 \to G^{*} \to 0.
 \end{equation}
 Since \(G\subset \SL(V)\),  \(M\) contains \((1,1,1)\).
 If we put \(M'= M/\Z(1,1,1)\) and \(H = \Z^3/\Z(1,1,1)\) then we have an exact sequence
 \[
0 \to M' \to H \to G^{*} \to 0.
 \]
 Denote by \(f_1, f_2, f_3 \in H\) the images of \((1,0,0)\), \((0,1,0)\), \((0,0,1)\in M\) respectively.

 \begin{definition}
 The \emph{universal cover} of the McKay quiver is the quiver with vertex set \(H\) and arrows from \(h \in H\) to \(h + f_i\) for \(i=1,2,3\).  
 \end{definition}

 Alternatively,  the universal cover may be regarded as follows.  Identify \(\Z^{3}\) with the lattice of Laurent monomials in \(x,y,z\),  so \(M\) is the sublattice of \(G\)-invariant Laurent monomials.  Then \(H\) becomes the lattice of Laurent monomials in \(x,y,z\) modulo \(xyz\) and arrows correspond to `multiply by \(x,y\) or \(z\)'.  In the \(\ghilb\) case it is convenient to choose monomials in \(\C[x,y,z]/xyz\) (see \cite[\S7]{Reid:mc} for a picture),  but in general we do not make this choice.  

 The McKay quiver is the quotient of this universal cover by the action of \(M'\).
 By embedding the universal cover into \(H_{\R} = H \otimes \R\),   the McKay quiver becomes embedded in the two dimensional real torus \(H_{\R}/M'\).

 Let \(F\) be a \(G\)-constellation and write \(B = (b_i^{\rho})\) for the corresponding representation.
 Define a graph \(\Gamma_F\) whose vertex set is \(G^{*}\) and whose edges are \(a_i^{\rho}\) with \(b_i^{\rho} \ne 0\),  forgetting the orientation.
 This graph can be embedded into \(H_{\R}/M'\): the edge corresponding to \(a_i^{\rho}\) is the image of the line segment connecting a lift \(\widetilde{\rho}\in H_{\R}\) and \(\widetilde{\rho} + f_i\) (any lift will do since we work modulo \(M'\)).

 The commutativity condition ensures that if \(a_i^{\rho}\) and \(a_j^{\rho\rho_i}\) are edges of \(\Gamma_F\),  then so are \(a_j^{\rho}\) and \(a_i^{\rho\rho_j}\) for \(\rho \in G^*\) and \(i \ne j\).
 Let \(C(\rho\, ; i, j)\) denote the cycle formed by arrows \(a_i^{\rho}\), \(a_j^{\rho\rho_i}\), \((a_i^{\rho\rho_j})^{-1}\) and \((a_j^{\rho})^{-1}\) that occur in any single commutativity equation.
 This is the image under \(H_{\R}\to H_{\R}/M'\) of a parallelogram (or a \emph{diamond}) in \(H_{\R}\).
 We denote by \(D(\rho \, ; i, j) \subset H_{\R}/M'\) the image of the closed domain surrounded by a diamond mapped to \(C(\rho\, ; i, j)\).  Let \(\Cycle(F)\) denote the set of such cycles in \(\Gamma_F\) satisfying \(i<j\).
 Assume \(F\) is not on the big torus \(T\) of \(\mc{M}_C\).
 In this case, if \(C(\rho\, ; i,j) \in \Cycle(F)\) then the \emph{diagonal} \(a_k^{\rho \rho_i \rho_j}\) of the diamond is not contained in \(\Gamma_F\),  where \(\{i, j, k\}=\{1,2,3\}\).

 \begin{example}
 \label{ex:unicover}
For the cyclic quotient singularity of type \(\frac{1}{6}(1,2,3)\),  the fan \(\Sigma\) defining \(\ghilb = X_{\Sigma}\) is shown in Craw~\cite[Figure~5.6(a)]{Craw:thesis}.  Consider a point \(F\) on the two-dimensional torus orbit lying in the unique compact toric surface in \(\ghilb\).  The covering of \(H_{\R}\) by lifts of the graph \(\Gamma_{F}\subset H_{\R}/M'\) is shown in Figure~\ref{fig:unicover}:  the vertices corresponding to the trivial representation are marked with circles;  the arrows \(f_{1}\) are drawn horizontally,  the \(f_{2}\) are vertical and the \(f_{3}\) point down and to the left.  Only the arrows \(a_{i}^{\rho}\) for which \(b_{i}^{\rho}\neq 0\) are illustrated.

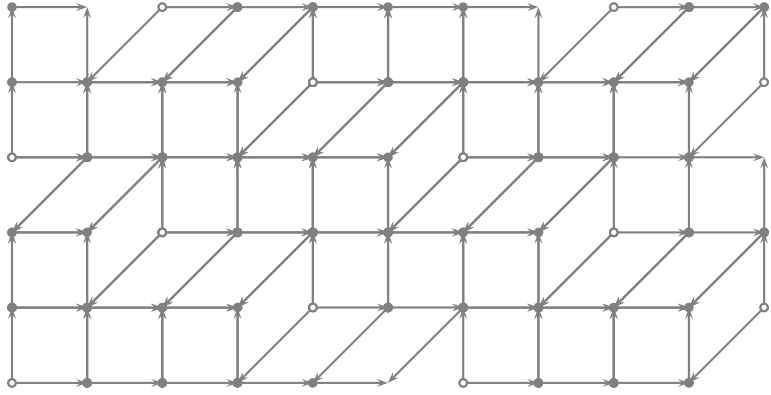
\begin{figure}[!ht]
 \centering
 \psset{unit=1cm}
 \begin{pspicture}(10,5.5)
 \put(0,0.2){
 \psset{linecolor=gray}
  \rput(4,1){\pstextpath[c]{\psline{o->}(0,0)(0,1)}{}
 \pstextpath[c]{\psline{o->}(0,0)(-1,-1)}{}
 \pstextpath[c]{\psline{*->}(0,1)(-1,0)}{}
 \pstextpath[c]{\psline{*-*}(-1,-1)(-1,0)}{}
}
 \rput(3,2){\pstextpath[c]{\psline{o->}(0,0)(1,0)}{}
 \pstextpath[c]{\psline{o->}(0,0)(-1,-1)}{}
 \pstextpath[c]{\psline{*->}(1,0)(0,-1)}{}
 \pstextpath[c]{\psline{*->}(-1,-1)(0,-1)}{}
} 
 \rput(2,2){\pstextpath[c]{\psline{o->}(0,0)(1,0)}{}
 \pstextpath[c]{\psline{o->}(0,0)(-1,-1)}{}
 \pstextpath[c]{\psline{*->}(1,0)(0,-1)}{}
 \pstextpath[c]{\psline{*->}(-1,-1)(0,-1)}{}
}
 \rput(2,0){\pstextpath[c]{\psline{o->}(0,0)(1,0)}{}
 \pstextpath[c]{\psline{o->}(0,0)(0,1)}{}
 \pstextpath[c]{\psline{*->}(1,0)(1,1)}{}
 \pstextpath[c]{\psline{*->}(0,1)(1,1)}{}
}
 \rput(1,0){\pstextpath[c]{\psline{o->}(0,0)(1,0)}{}
 \pstextpath[c]{\psline{o->}(0,0)(0,1)}{}
 \pstextpath[c]{\psline{*->}(1,0)(1,1)}{}
 \pstextpath[c]{\psline{*->}(0,1)(1,1)}{}
}
 \rput(0,0){\pstextpath[c]{\psline{o->}(0,0)(1,0)}{}
 \pstextpath[c]{\psline{o->}(0,0)(0,1)}{}
 \pstextpath[c]{\psline{*->}(1,0)(1,1)}{}
 \pstextpath[c]{\psline{*->}(0,1)(1,1)}{}
}
  \rput(8,2){\pstextpath[c]{\psline{o->}(0,0)(0,1)}{}
 \pstextpath[c]{\psline{o->}(0,0)(-1,-1)}{}
 \pstextpath[c]{\psline{*->}(0,1)(-1,0)}{}
 \pstextpath[c]{\psline{*-*}(-1,-1)(-1,0)}{}
}
 \rput(7,3){\pstextpath[c]{\psline{o->}(0,0)(1,0)}{}
 \pstextpath[c]{\psline{o->}(0,0)(-1,-1)}{}
 \pstextpath[c]{\psline{*->}(1,0)(0,-1)}{}
 \pstextpath[c]{\psline{*->}(-1,-1)(0,-1)}{}
} 
 \rput(6,3){\pstextpath[c]{\psline{o->}(0,0)(1,0)}{}
 \pstextpath[c]{\psline{o->}(0,0)(-1,-1)}{}
 \pstextpath[c]{\psline{*->}(1,0)(0,-1)}{}
 \pstextpath[c]{\psline{*->}(-1,-1)(0,-1)}{}
}
 \rput(6,1){\pstextpath[c]{\psline{o->}(0,0)(1,0)}{}
 \pstextpath[c]{\psline{o->}(0,0)(0,1)}{}
 \pstextpath[c]{\psline{*->}(1,0)(1,1)}{}
 \pstextpath[c]{\psline{*->}(0,1)(1,1)}{}
}
 \rput(5,1){\pstextpath[c]{\psline{o->}(0,0)(1,0)}{}
 \pstextpath[c]{\psline{o->}(0,0)(0,1)}{}
 \pstextpath[c]{\psline{*->}(1,0)(1,1)}{}
 \pstextpath[c]{\psline{*->}(0,1)(1,1)}{}
}
 \rput(4,1){\pstextpath[c]{\psline{o->}(0,0)(1,0)}{}
 \pstextpath[c]{\psline{o->}(0,0)(0,1)}{}
 \pstextpath[c]{\psline{*->}(1,0)(1,1)}{}
 \pstextpath[c]{\psline{*->}(0,1)(1,1)}{}
}
  \rput(6,3){\pstextpath[c]{\psline{o->}(0,0)(0,1)}{}
 \pstextpath[c]{\psline{o->}(0,0)(-1,-1)}{}
 \pstextpath[c]{\psline{*->}(0,1)(-1,0)}{}
 \pstextpath[c]{\psline{*-*}(-1,-1)(-1,0)}{}
}
 \rput(5,4){\pstextpath[c]{\psline{o->}(0,0)(1,0)}{}
 \pstextpath[c]{\psline{o->}(0,0)(-1,-1)}{}
 \pstextpath[c]{\psline{*->}(1,0)(0,-1)}{}
 \pstextpath[c]{\psline{*->}(-1,-1)(0,-1)}{}
} 
 \rput(4,4){\pstextpath[c]{\psline{o->}(0,0)(1,0)}{}
 \pstextpath[c]{\psline{o->}(0,0)(-1,-1)}{}
 \pstextpath[c]{\psline{*->}(1,0)(0,-1)}{}
 \pstextpath[c]{\psline{*->}(-1,-1)(0,-1)}{}
}
 \rput(4,2){\pstextpath[c]{\psline{o->}(0,0)(1,0)}{}
 \pstextpath[c]{\psline{o->}(0,0)(0,1)}{}
 \pstextpath[c]{\psline{*->}(1,0)(1,1)}{}
 \pstextpath[c]{\psline{*->}(0,1)(1,1)}{}
}
 \rput(3,2){\pstextpath[c]{\psline{o->}(0,0)(1,0)}{}
 \pstextpath[c]{\psline{o->}(0,0)(0,1)}{}
 \pstextpath[c]{\psline{*->}(1,0)(1,1)}{}
 \pstextpath[c]{\psline{*->}(0,1)(1,1)}{}
}
 \rput(2,2){\pstextpath[c]{\psline{o->}(0,0)(1,0)}{}
 \pstextpath[c]{\psline{o->}(0,0)(0,1)}{}
 \pstextpath[c]{\psline{*->}(1,0)(1,1)}{}
 \pstextpath[c]{\psline{*->}(0,1)(1,1)}{}
}
  \rput(4,4){\pstextpath[c]{\psline{o->}(0,0)(0,1)}{}
 \pstextpath[c]{\psline{o->}(0,0)(-1,-1)}{}
 \pstextpath[c]{\psline{*->}(0,1)(-1,0)}{}
 \pstextpath[c]{\psline{*-*}(-1,-1)(-1,0)}{}
}
 \rput(3,5){\pstextpath[c]{\psline{o->}(0,0)(1,0)}{}
 \pstextpath[c]{\psline{o->}(0,0)(-1,-1)}{}
 \pstextpath[c]{\psline{*->}(1,0)(0,-1)}{}
 \pstextpath[c]{\psline{*->}(-1,-1)(0,-1)}{}
} 
 \rput(2,5){\pstextpath[c]{\psline{o->}(0,0)(1,0)}{}
 \pstextpath[c]{\psline{o->}(0,0)(-1,-1)}{}
 \pstextpath[c]{\psline{*->}(1,0)(0,-1)}{}
 \pstextpath[c]{\psline{*->}(-1,-1)(0,-1)}{}
}
 \rput(2,3){\pstextpath[c]{\psline{o->}(0,0)(1,0)}{}
 \pstextpath[c]{\psline{o->}(0,0)(0,1)}{}
 \pstextpath[c]{\psline{*->}(1,0)(1,1)}{}
 \pstextpath[c]{\psline{*->}(0,1)(1,1)}{}
}
 \rput(1,3){\pstextpath[c]{\psline{o->}(0,0)(1,0)}{}
 \pstextpath[c]{\psline{o->}(0,0)(0,1)}{}
 \pstextpath[c]{\psline{*->}(1,0)(1,1)}{}
 \pstextpath[c]{\psline{*->}(0,1)(1,1)}{}
}
 \rput(0,3){\pstextpath[c]{\psline{o->}(0,0)(1,0)}{}
 \pstextpath[c]{\psline{o->}(0,0)(0,1)}{}
 \pstextpath[c]{\psline{*->}(1,0)(1,1)}{}
 \pstextpath[c]{\psline{*->}(0,1)(1,1)}{}
}
 \rput(6,4){\pstextpath[c]{\psline{o->}(0,0)(1,0)}{}
 \pstextpath[c]{\psline{o->}(0,0)(0,1)}{}
 \pstextpath[c]{\psline{*->}(1,0)(1,1)}{}
 \pstextpath[c]{\psline{*->}(0,1)(1,1)}{}
}
 \rput(5,4){\pstextpath[c]{\psline{o->}(0,0)(1,0)}{}
 \pstextpath[c]{\psline{o->}(0,0)(0,1)}{}
 \pstextpath[c]{\psline{*->}(1,0)(1,1)}{}
 \pstextpath[c]{\psline{*->}(0,1)(1,1)}{}
}
 \rput(4,4){\pstextpath[c]{\psline{o->}(0,0)(1,0)}{}
 \pstextpath[c]{\psline{o->}(0,0)(0,1)}{}
 \pstextpath[c]{\psline{*->}(1,0)(1,1)}{}
 \pstextpath[c]{\psline{*->}(0,1)(1,1)}{}
}
  \rput(10,1){\pstextpath[c]{\psline{o->}(0,0)(0,1)}{}
 \pstextpath[c]{\psline{o->}(0,0)(-1,-1)}{}
 \pstextpath[c]{\psline{*->}(0,1)(-1,0)}{}
 \pstextpath[c]{\psline{*-*}(-1,-1)(-1,0)}{}
}
 \rput(9,2){\pstextpath[c]{\psline{o->}(0,0)(1,0)}{}
 \pstextpath[c]{\psline{o->}(0,0)(-1,-1)}{}
 \pstextpath[c]{\psline{*->}(1,0)(0,-1)}{}
 \pstextpath[c]{\psline{*->}(-1,-1)(0,-1)}{}
} 
 \rput(8,2){\pstextpath[c]{\psline{o->}(0,0)(1,0)}{}
 \pstextpath[c]{\psline{o->}(0,0)(-1,-1)}{}
 \pstextpath[c]{\psline{*->}(1,0)(0,-1)}{}
 \pstextpath[c]{\psline{*->}(-1,-1)(0,-1)}{}
}
 \rput(8,0){\pstextpath[c]{\psline{o->}(0,0)(1,0)}{}
 \pstextpath[c]{\psline{o->}(0,0)(0,1)}{}
 \pstextpath[c]{\psline{*->}(1,0)(1,1)}{}
 \pstextpath[c]{\psline{*->}(0,1)(1,1)}{}
}
 \rput(7,0){\pstextpath[c]{\psline{o->}(0,0)(1,0)}{}
 \pstextpath[c]{\psline{o->}(0,0)(0,1)}{}
 \pstextpath[c]{\psline{*->}(1,0)(1,1)}{}
 \pstextpath[c]{\psline{*->}(0,1)(1,1)}{}
}
 \rput(6,0){\pstextpath[c]{\psline{o->}(0,0)(1,0)}{}
 \pstextpath[c]{\psline{o->}(0,0)(0,1)}{}
 \pstextpath[c]{\psline{*->}(1,0)(1,1)}{}
 \pstextpath[c]{\psline{*->}(0,1)(1,1)}{}
}
  \rput(10,4){\pstextpath[c]{\psline{o->}(0,0)(0,1)}{}
 \pstextpath[c]{\psline{o->}(0,0)(-1,-1)}{}
 \pstextpath[c]{\psline{*->}(0,1)(-1,0)}{}
 \pstextpath[c]{\psline{*-*}(-1,-1)(-1,0)}{}
}
 \rput(9,5){\pstextpath[c]{\psline{o->}(0,0)(1,0)}{}
 \pstextpath[c]{\psline{o->}(0,0)(-1,-1)}{}
 \pstextpath[c]{\psline{*->}(1,0)(0,-1)}{}
 \pstextpath[c]{\psline{*->}(-1,-1)(0,-1)}{}
} 
 \rput(8,5){\pstextpath[c]{\psline{o->}(0,0)(1,0)}{}
 \pstextpath[c]{\psline{o->}(0,0)(-1,-1)}{}
 \pstextpath[c]{\psline{*->}(1,0)(0,-1)}{}
 \pstextpath[c]{\psline{*->}(-1,-1)(0,-1)}{}
}
 \rput(8,3){\pstextpath[c]{\psline{o->}(0,0)(1,0)}{}
 \pstextpath[c]{\psline{o->}(0,0)(0,1)}{}
 \pstextpath[c]{\psline{*->}(1,0)(1,1)}{}
 \pstextpath[c]{\psline{*->}(0,1)(1,1)}{}
}
 \rput(7,3){\pstextpath[c]{\psline{o->}(0,0)(1,0)}{}
 \pstextpath[c]{\psline{o->}(0,0)(0,1)}{}
 \pstextpath[c]{\psline{*->}(1,0)(1,1)}{}
 \pstextpath[c]{\psline{*->}(0,1)(1,1)}{}
}
 \rput(6,3){\pstextpath[c]{\psline{o->}(0,0)(1,0)}{}
 \pstextpath[c]{\psline{o->}(0,0)(0,1)}{}
 \pstextpath[c]{\psline{*->}(1,0)(1,1)}{}
 \pstextpath[c]{\psline{*->}(0,1)(1,1)}{}
}
  \rput(2,2){\pstextpath[c]{\psline{o->}(0,0)(0,1)}{}
 \pstextpath[c]{\psline{o->}(0,0)(-1,-1)}{}
 \pstextpath[c]{\psline{*->}(0,1)(-1,0)}{}
 \pstextpath[c]{\psline{*-*}(-1,-1)(-1,0)}{}
}
 \rput(1,3){\pstextpath[c]{\psline{o->}(0,0)(1,0)}{}
 \pstextpath[c]{\psline{o->}(0,0)(-1,-1)}{}
 \pstextpath[c]{\psline{*->}(1,0)(0,-1)}{}
 \pstextpath[c]{\psline{*->}(-1,-1)(0,-1)}{}
} 
 \rput(0,3){\pstextpath[c]{\psline{o->}(0,0)(1,0)}{}
 \pstextpath[c]{\psline{*->}(1,0)(0,-1)}{}
}
 \rput(0,1){\pstextpath[c]{\psline{o->}(0,0)(1,0)}{}
 \pstextpath[c]{\psline{o->}(0,0)(0,1)}{}
 \pstextpath[c]{\psline{*->}(1,0)(1,1)}{}
 \pstextpath[c]{\psline{*->}(0,1)(1,1)}{}
}
 \rput(-1,1){
 \pstextpath[c]{\psline{*->}(1,0)(1,1)}{}
}
\rput(0,4){\pstextpath[c]{\psline{*->}(0,0)(0,1)}{}
\pstextpath[c]{\psline{*->}(0,1)(1,1)}{}
\pstextpath[c]{\psline{*->}(1,0)(1,1)}{}
 }
 \rput(3,0){\pstextpath[c]{\psline{*->}(0,0)(1,0)}{}
\pstextpath[c]{\psline{*->}(1,0)(2,0)}{}
\pstextpath[c]{\psline{*->}(2,1)(1,0)}{}
\pstextpath[c]{\psline{*->}(3,1)(2,0)}{}
 }
\rput(9,2){\pstextpath[c]{\psline{*->}(0,0)(0,1)}{}
\pstextpath[c]{\psline{*->}(0,1)(1,1)}{}
\pstextpath[c]{\psline{*->}(1,0)(1,1)}{}
 }
}
 \end{pspicture}
 \caption{Lifts of diamonds to $H_{\R}$}
 \label{fig:unicover}
 \end{figure} 

 \end{example}

 \begin{lemma}
 \label{lemma:domains}
 Assume \(F\) corresponds to a point on a two-dimensional $T$-orbit in $\mc M_{C}$. Then
 \[
H_{\R}/M'= \bigcup_{C(\rho\, ;i,j) \in \Cycle(F)} D(\rho\, ; i, j),
 \]
 where the domains intersect only along their common boundaries.  
 \end{lemma}
 \begin{proof}
 Let \(P\) be a point in \(H_{\R}/M'\).
 Since every representation occurs in \(F\),  \(P\) is contained in a closed domain surrounded by a triangle formed by \(a_1^{\rho}, a_2^{\rho\rho_1},  a_3^{\rho\rho_1\rho_2}\) or \(a_1^{\rho}, a_3^{\rho\rho_1},  a_2^{\rho\rho_1\rho_3}\) for some \(\rho\).
 Assume the former case, the latter is similar.
 Corollary~\ref{coro:oneis0} asserts that exactly one of \(b_1^{\rho}, b_2^{\rho\rho_1},  b_3^{\rho\rho_1\rho_2}\) is \(0\). Assume the last one is zero. 
 Then the commutativity condition \(a_2^{\rho\rho_1}a_1^{\rho}=a_1^{\rho\rho_2}a_2^{\rho}\) implies that \(C(\rho\,;1,2) \in \Cycle(F)\).
 We have \(P \in D(\rho\,;1,2)\) and therefore \(P\) is contained in the right hand side. The other cases are similar.
 \end{proof}

 \subsection{Rigidity result}

 A \(G\)-sheaf \(E\) is said to be \emph{simple} if \(\dim\GHom_{\owe_{\C^3}}(E,E)=1\).

 \begin{proposition}
 \label{prop:rigidity}
 Assume \(F\) corresponds to a point on a two-dimensional \(T\)-orbit in \(\mc{M}_{C}\).
 Take a \(G\)-subsheaf \(S\) of \(F\) with quotient \(Q\).
 \begin{enumerate}
 \item[\one] If both \(S\) and \(Q\) are simple then \(\dim \GExt^1_{\owe_{\C^3}}(Q, S)\leq 2\).
 \item[\two] If in addition \(\dim \GExt^1_{\owe_{\C^3}}(Q,S)=1\), then either $S$ or $Q$ is a rigid \(G\)-sheaf.
 \end{enumerate}
 \end{proposition}
 \begin{proof}
 Let \(\Gamma_S\) be a subgraph of $\Gamma_F$ whose vertex set consists of representations in $S$ and whose edges are those of $\Gamma_F$ that join two vertices in \(\Gamma_{S}\).
 We define $\Gamma_Q$ similarly.
 The assumption that $S$ and $Q$ are simple implies that both $\Gamma_S$ and $\Gamma_Q$ are connected.
 Since $S$ is a subsheaf, if $a_i^{\rho}$ is an edge of \(\Gamma_{F}\) that lies in neither $\Gamma_S$ nor $\Gamma_Q$, then $a_i^{\rho}$ is an arrow from a vertex in $\Gamma_Q$ to one in $\Gamma_S$.
 Denote by \(\Dom{S}\) the union of \(\Gamma_S\) with diamonds \(D(\rho\,;i,j)\) such that \(C(\rho\,;i,j) \in \Cycle(F)\) is contained in \(\Gamma_S\).
 Similarly,  \(\Dom{Q}\) denotes the domain determined by \(\Gamma_Q\).
 Let \(\Band(Q,S)\) be the complement of \(\Dom{S} \cup \Dom{Q}\).
 Lemma~\ref{lemma:domains} implies that if \(\Band(Q,S) \cap D(\rho\,;i,j)\) is nonempty then it is the diamond \(D(\rho\,;i,j)\) with two edges removed.  Moreover,  if the removed edges are not parallel then they must be adjacent edges whose orientations (as arrows) are different.
 The remaining edges of \(D(\rho\,;i,j)\) intersected with \(\Band(Q,S)\) join this set to an adjacent pair of such diamonds and hence form a connected component of \(\Band(Q,S)\) that is homeomorphic to \(S^1 \times (0,1)\).

 Now consider \(\GExt^1_{\owe_{\C^3}}(Q,S)\).
 Fix representations 
\[
 B_{11}\colon \Ro \to V \otimes \Ro\quad\text{and}\quad B_{22}\colon \Rt \to V \otimes \Rt
 \]
 that give rise to the \(G\)-sheaves \(S\) and \(Q\) respectively, for \(R = \Ro\oplus\Rt\).
 An extension class in \(\GExt^1_{\owe_{\C^3}}(Q,S)\) is given by a map \(B\colon [\Ro \oplus \Rt] \to V \otimes [\Ro \oplus \Rt]\) (up to automorphisms) satisfying the commutativity condition, such that the $(1,1)$- and $(2,2)$-entries are $B_{11}$ and $B_{22}$, and furthermore $B_{21}\colon \Ro \to V \otimes \Rt$ is the zero map.  Thus
 \[
\GExt^1_{\owe_{\C^3}}(Q,S)\cong\frac{
\left\{ B_{12}:R_2 \to V \otimes R_1\,\left|\,
B \wedge B=0 \text{ for } B=\begin{pmatrix} B_{11} & B_{12} \\ 0 & B_{22}\end{pmatrix}
\right.\right\}}
{\left\{ B_{22}A-AB_{11} \,\left|\,
A \in \GHom_{\C}(R_2, R_1) \right.\right\}}.
 \]
 In our case, the denominator is zero so we consider the numerator.  
 Each $B_{12}$ is determined by the values of arrows dividing $\overline{\Band(Q,S)}$ into diamonds.
 The commutativity condition for cycles \(C(\rho\,; i, j) \in \Cycle(F)\) contained in $\overline{\Band(Q,S)}$ implies that the value of such an arrow is determined by any of the two adjacent ones (together with $B_{11}$ and \(B_{22}\)).
 Therefore, $\dim \GExt^1_{\owe_{\C^3}}(Q,S)$ is the number of connected components of $\Band(Q,S)$, which we denote by $n$.
 Since the complement of $\Band(Q,S)$ in the torus $H_{\R}/M'$ has exactly two connected components $\Dom{S}$ and $\Dom{Q}$, we obtain $n \leq 2$. 
 Thus we obtain the first assertion.

 Assume $n=1$ so that $\Band(Q,S) \cong S^1 \times (0,1)$ is connected.
 Then either $\Dom{S}$ or $\Dom{Q}$ must be simply connected.
 We show if $\Dom{S}$ is simply connected, then $S$ is rigid, the other case is similar.
 The condition $B_{11} \wedge B_{11}=0$ implies that the values of deformations of $B_{11}$ at arrows not in $\Gamma_S$ must be zero.
 Hence it is enough to show that up to the action of $\Paut_{\C[G]}(R_1)$, we can assume $b_i^{\rho}=1$ if $a_i^{\rho} \subset \Gamma_S$.
 Fix $\sigma_0 \subset R_1$ and take an arbitrary $\sigma \subset R_1$.
 Since $\Gamma_S$ is connected, we can take a sequence of arrows
 \[
\gamma=\{\,\alpha_1^{\varepsilon_1}, \alpha_2^{\varepsilon_2}, \dots, 
\alpha_p^{\varepsilon_p}\,\}
 \]
 such that $\alpha_k$ is an arrow of the McKay quiver supported by $\Gamma_S$,  $\varepsilon_i= \pm 1$ determines the direction of arrows and such that $\gamma$ gives a connected path from $\sigma_0$ to $\sigma$.
 Put
 \[
\lambda_{\sigma} = \prod_{k=1}^{p} \beta_k^{-\varepsilon_k}
 \]
 where $\beta_k=b_i^{\rho}$ is the value of $B_{11}$ at the arrow $\alpha_k=a_i^{\rho}$.
 Since $\Dom{S}$ is simply connected, the first homology of $\Gamma_S$ is generated by cycles $C(\rho\,;i,j)$ contained in $\Gamma_S$.
 Then the commutativity condition shows that $\lambda_{\sigma}$ is independent of the choice of path $\gamma$.
 Put \(\lambda=(\lambda_{\sigma})_{\sigma \subset R_1} \in \Aut'_{\C[G]}(R_1)\).
 Then, $\lambda \cdot B_{11}$ has the desired property.
\end{proof}

 \begin{example}
 Consider once again Example~\ref{ex:unicover} and Figure~\ref{fig:unicover}. The quotient \(Q = \rho_0 + \rho_1\) is rigid and \(\Band(Q,S)\) is connected.  On the other hand,   the quotient \(Q = \rho_0 + \rho_1 + \rho_3\) is not rigid and \(\Band(Q, S)\) has two connected components.
We leave the details to the reader.
 \end{example}

 \medskip

 \noindent \textsc{Department of Mathematics,  University of Utah,\\
 155 South 1400 East, Salt Lake City,UT 84112, USA. \\
 E-mail:} \texttt{craw@math.utah.edu}
 
 \medskip

 \noindent \textsc{Department of Civil Engineering Systems, Kyoto University,\\
 Kyoto 606-8501, Japan\\
 E-mail:} \texttt{akira@kusm.kyoto-u.ac.jp}

 \end{document}